\documentclass[reqno,12pt]{amsart}
\usepackage{latexsym,amsmath,amsfonts,amscd,amssymb}
\renewcommand\today{March 11, 2003}

\theoremstyle{plain}  
\newtheorem{theorem}{Theorem}[section]
\newtheorem*{theorem*}{Theorem}
\newtheorem*{theoremA}{Theorem A}

\newtheorem{corollary}[theorem]{Corollary}
\newtheorem{lemma}[theorem]{Lemma}
\newtheorem{proposition}[theorem]{Proposition}

\newtheorem{tech-lemma}[theorem]{Technical Lemma}
\theoremstyle{definition}
\newtheorem{definition}[theorem]{Definition}
\theoremstyle{remark}

\newtheorem{remark}[theorem]{Remark}

\newtheorem*{claim*}{Claim}
\numberwithin{equation}{section}
\renewcommand{\leq}{\leqslant}
\renewcommand{\le}{\leqslant}
\renewcommand{\geq}{\geqslant}
\renewcommand{\ge}{\geqslant}

\newcommand{\lto}{\longrightarrow}
\newcommand{\R}{\mathbb{R}}

\newcommand{\C}{\mathbb{C}}
\newcommand{\HH}{\mathbb{H}}

\newcommand{\into}{\hookrightarrow}

\newcommand{\vol}{\mathrm{vol}}

\newcommand{\balpha}{{\boldsymbol{\alpha}}}
\newcommand{\btau}{{\boldsymbol{\tau}}}
\newcommand{\PU}{\mathrm{PU}}
\newcommand{\PGL}{\mathrm{PGL}} 
\newcommand{\SU}{\mathrm{SU}}
\newcommand{\SL}{\mathrm{SL}}
\newcommand{\U}{\mathrm{U}}
\newcommand{\GL}{\mathrm{GL}}
\newcommand{\GCD}{\mathrm{GCD}}

\DeclareMathOperator{\Sym}{Sym}

\DeclareMathOperator{\rk}{rk}
\DeclareMathOperator{\im}{im}
\DeclareMathOperator{\coker}{coker}
\DeclareMathOperator{\Hom}{Hom}

\DeclareMathOperator{\Ext}{Ext}
\DeclareMathOperator{\Id}{Id}
\DeclareMathOperator{\Quot}{Quot}

 \newcommand{\be}{\begin{equation}}
 \newcommand{\ee}{\end{equation}}
 
 \newfont{\bfc}{cmbsy10 scaled 1200}  
 \newfont{\gl}{eufm10 scaled \magstep1}  

 \newcommand{\dP}{\mathbb{P}}

 \newcommand{\cO}{{\mathcal O}}

 \newcommand{\lra}{\longrightarrow}
 \newcommand{\kahler}{K\"{a}hler}
 \newcommand{\ps}{{p^\ast}}
 \newcommand{\qs}{{q^\ast}}
 
 \newcommand{\xp}{{{X\times\dP^1}}}

 \newcommand{\he}{Hermitian--Einstein}

\begin{document}
\begin{titlepage}
  
\noindent {\Large\textbf{Moduli spaces of holomorphic triples 
over compact Riemann surfaces}} 
\bigskip

\noindent 
\textbf{Steven B. Bradlow}\footnotemark[1]$^{,}$\footnotemark[2] \\
Department of Mathematics, \\
University of Illinois, \\
Urbana,
IL 61801, 
USA \\
E-mail: \texttt{bradlow@math.uiuc.edu}
\medskip

\noindent 
\textbf{Oscar Garc{\'\i}a--Prada}\footnotemark[1]$^{,}
$\footnotemark[3]$^{,}$\footnotemark[5]$^{,}$\footnotemark[6]  \\
Instituto de Matem\'aticas y F\'{\i}sica Fundamental,\\
Consejo Superior de Investigaciones Cient\'{\i}ficas, \\
Serrano 113 bis,\\
28006 Madrid, Spain \\
E-mail: \texttt{oscar.garcia-prada@uam.es}

\medskip
\noindent 
\textbf{Peter B. Gothen}\footnotemark[1]$^{,}$\footnotemark[4]$^{,}$\footnotemark[5]   \\
Departamento de Matem{\'a}tica Pura, \\
Faculdade de Ci{\^e}ncias,
Universidade do Porto, \\
Rua do Campo Alegre 687, 4169-007 Porto,
Portugal \\
E-mail: \texttt{pbgothen@fc.up.pt}

\footnotetext[1]{
Members of VBAC (Vector Bundles on
Algebraic Curves), which is partially supported by EAGER (EC FP5
Contract no.\ HPRN-CT-2000-00099) and by EDGE (EC FP5 Contract no.\ 
HPRN-CT-2000-00101).}
\footnotetext[2]{Partially supported by 
the National Science Foundation under grant DMS-0072073 }
\footnotetext[3]{Partially supported by 
the Ministerio de Ciencia y Tecnolog\'{\i}a (Spain) under grant BFM2000-0024}
\footnotetext[4]{
Partially supported by the
Funda{\c c}{\~a}o para a Ci{\^e}ncia e a Tecnologia (Portugal) through
the Centro de Matem{\'a}tica da Universidade do Porto and through
grant no.\ SFRH/BPD/1606/2000.}
\footnotetext[5]{
Partially supported by the Portugal/Spain bilateral Programme Acciones
Integradas, 
grant nos.\ HP2000-0015 and AI-01/24}
\footnotetext[6]{
Partially supported by a  British EPSRC grant (October-December 2001)}

\bigskip
\noindent
\textbf{\today}
\vfill

\subsection*{Abstract.}
A holomorphic triple over a compact Riemann surface consists of two
holomorphic vector bundles and a holomorphic map between them.  After
fixing the topological types of the bundles and a real parameter,
there exist moduli spaces of stable holomorphic triples. In this paper
we study non-emptiness, irreducibility, smoothness, and birational
descriptions of these moduli spaces for a certain range of the
parameter. Our results have important applications to the study of the moduli
space of representations of the fundamental group of the surface into
unitary Lie groups of indefinite signature
(\cite{bradlow-garcia-gothen:2001,bradlow-garcia-gothen:2002}).
Another application, that we study in this paper, is to the existence
of stable bundles on the product of the surface by the complex
projective line.

\vfill

\newpage

\end{titlepage}

\addtocounter{page}{1}

\section{Introduction} \label{sec:introduction}

Let $X$ be a closed Riemann surface of genus $g \geq 2$.  The theory
of holomorphic triples has its origins
\cite{garcia-prada:1994,bradlow-garcia-prada:1996} in the search for
solutions to certain gauge theoretic equations on $X$, obtained by
dimensional reduction of the Hermitian--Einstein equation in 4
dimensions.  More precisely, solutions to the Hermitian--Einstein
equation on $X \times \mathbb{P}^1$ which are invariant under the
standard action of $\SU(2)$ on $\mathbb{P}^1$ correspond to solutions
to the so-called \emph{vortex equations} on $X$.  The
Hitchin--Kobayashi correspondence states that a solution to the
Hermitian--Einstein equation on $X \times \mathbb{P}^1$ gives rise to
a stable holomorphic bundle and that, conversely, any stable
holomorphic bundle admits a Hermitian--Einstein metric.  The
counterpart on $X$ states that there is a Hitchin--Kobayashi
correspondence between solutions to the coupled vortex equations and
\emph{stable holomorphic triples}. A holomorphic triple consists of a
pair of holomorphic vector bundles, $E_1$ and $E_2$, over $X$ and a
holomorphic map $\phi\colon E_2 \to E_1$ between them.  An important
feature of the stability condition for triples is that it depends on a
real parameter $\alpha$, corresponding to the fact that there is a
real parameter in the  vortex equations; thus
one is led to the concept of $\alpha$-stability of a holomorphic triple.
This parallels the fact that when studying Hermitian--Einstein metrics
and stable bundles on $X \times \mathbb{P}^1$ it is necessary to
choose a polarization on this complex surface.  We note that, as usual, there
are corresponding concepts of $\alpha$-polystable and
$\alpha$-semistable triples (see Section~\ref{sec:stable-triples}
below for precise definitions).

It was shown in \cite{bradlow-garcia-prada:1996} (see also
\cite{garcia-prada:1994}) that projective moduli spaces for
holomorphic triples exist.  (Later a direct 
 construction was given by
Schmitt using geometric in\-variant theory \cite{Schmitt}.)  Since the
stability condition depends on the real parameter $\alpha$, so do the
moduli spaces.  Fixing the topological invariants $n_i = \rk(E_i)$ and
$d_i = \deg(E_i)$, we denote the moduli space of $\alpha$-polystable
triples with the given invariants by
\begin{displaymath}
    \mathcal{N}_\alpha
  = \mathcal{N}_\alpha(n_1,n_2,d_1,d_2)\ ,
\end{displaymath}
and the moduli space of $\alpha$-stable triples by $\mathcal{N}^s_\alpha 
\subseteq \mathcal{N}_\alpha$.
In this paper we address the questions of smoothness, non-emptiness
and irreducibility of these moduli spaces.

Before describing our results in more detail, we explain our
motivation, which comes from the problem of determining the connected
components of the moduli space of representations of the fundamental
group of $X$ in $\PU(p,q)$.  A detailed study of this moduli space
appears in a companion paper \cite{bradlow-garcia-gothen:2002}; in the
following we briefly outline the main ideas.  The first point to
notice is that we may as well study the
connected components of the moduli space of projectively flat
$\U(p,q)$ bundles on $X$.  This moduli space can be divided into
disjoint closed subspaces $\mathcal{M}(a,b)$ indexed by a pair of
integers $(a,b)$, the Chern classes obtained from a reduction of
structure group to the maximal compact subgroup $\U(p) \times \U(q)$.
The values of $(a,b)$ are bounded by the Milnor--Wood type inequality
\begin{displaymath}
  \left\lvert\frac{aq-bp}{p+q}\right\rvert \leq \min\{p,q\}(g-1)\ .
\end{displaymath}
For each allowed value of $(a,b)$ one  expects the space
$\mathcal{M}(a,b)$ to be non-empty and connected, thus forming a
connected component of the moduli space.

By the work of Hitchin \cite{hitchin:1987,hitchin:1992}, Donaldson
\cite{donaldson:1987}, Simpson
\cite{simpson:1988,simpson:1992,simpson:1994a,simpson:1994b} and
Corlette \cite{corlette:1988}, the moduli spaces $\mathcal{M}(a,b)$
are homeomorphic to moduli spaces of so-called \emph{$\U(p,q)$-Higgs bundles}
on $X$: these are pairs $(E,\Phi)$, where $E$ is a holomorphic vector
bundle which decomposes as a direct sum $E = V \oplus W$ and the Higgs
field $\Phi \colon E \to E \otimes K$ is of the form
\begin{displaymath}
  \phi =
  \begin{pmatrix}
    0 & \beta \\
    \gamma & 0
  \end{pmatrix}
\end{displaymath}
with respect to the direct sum decomposition of $E$.  Here $K$ is the
canonical line bundle of $X$ and the invariants $a$ and $b$ appear as
the degrees of $V$ and $W$ respectively.  The $L^2$-norm of the Higgs
field gives us a Bott-Morse function on the moduli space (cf.\ Hitchin
\cite{hitchin:1987,hitchin:1992}).  Thus, connectedness of the spaces
$\mathcal{M}(a,b)$ will be a consequence of connectedness of the
corresponding subspaces of local minima.  In the case of flat
$\U(2,2)$-bundles, it was shown in \cite{gothen:2001} that the local
minima are represented by Higgs bundles for which either $\beta$ or
$\gamma$ vanishes.  One of the main results in
\cite{bradlow-garcia-gothen:2002} is that this is true in general.
The crucial observation is now that there is a bijective
correspondence between $\U(p,q)$-Higgs bundles $(E,\Phi)$ with $\beta
= 0$ or $\gamma = 0$ and holomorphic triples: if, say, $\gamma = 0$,
we obtain a holomorphic triple $T = (E_1,E_2,\phi)$ by setting $E_1 =
V \otimes K$, $E_2 = W$ and $\phi = \beta$.  It turns out that
$(E,\Phi)$ is (poly)stable as a $\U(p,q)$-Higgs bundle if and only if
the corresponding holomorphic triple $T$ is $\alpha$-(poly)stable for
$\alpha = 2g-2$.  It follows that the subspace of local minima on
$\mathcal{M}(a,b)$ is isomorphic to a moduli space of
$(2g-2)$-polystable holomorphic triples.  Thus the results of the
present paper imply results on non-emptiness and connectedness of the
moduli spaces $\mathcal{M}(a,b)$.  We refer the reader to
\cite{bradlow-garcia-gothen:2002} for the precise statements.

We now return to our main subject of study, the holomorphic triples.
In order for $\mathcal{N}_\alpha$ to be non-empty, one must have
$\alpha\geq \alpha_m$ with $\alpha_m=d_1/n_1-d_2/n_2\geq 0$. In the
case $n_1\neq n_2$ there is also a finite upper bound $\alpha_M$.
When the parameter $\alpha$ varies, the nature of the
$\alpha$-stability condition only changes for a discrete number of
so-called \emph{critical values} of $\alpha$ (see
Section~\ref{sec:stable-triples} for the precise statements).  We can
now state our main results.

\begin{theoremA}
\begin{enumerate}
\item[$(1)$] A triple $T=(E_1,E_2,\phi)$ of type $(n_1,n_2,d_1,d_2)$
is $\alpha_m$-polystable if and only if 
$\phi=0$ and $E_1$ and $E_2$ are polystable. We thus have 
$$
\mathcal{N}_{\alpha_m}(n_1,n_2,d_1,d_2)\cong M(n_1,d_1) \times M(n_2,d_2).
$$
where $M(n,d)$ denotes the moduli space of polystable bundles of rank 
$n$ and degree $d$.  In particular, 
$\mathcal{N}_{\alpha_m}(n_1,n_2,d_1,d_2)$ is non-empty and irreducible. 

\item[$(2)$] If $\alpha>\alpha_m $ is any value such that
$2g-2\leq\alpha$ (and $\alpha<\alpha_M$ if $n_1\ne n_2$) then
the moduli space $\mathcal{N}^s_\alpha(n_1,n_2,d_1,d_2)$ is non-empty, 
irreducible, and smooth of dimension 
$(g-1)(n_1^2 + n_2^2 - n_1 n_2) - n_1 d_2 + n_2 d_1 + 1$.  Moreover:

\quad $\bullet$ If $n_1=n_2 =n$ then the moduli space 
$\mathcal{N}^s_{\alpha}(n,n,d_1,d_2)$ is 
birationally  equivalent to a 
$\mathbb{P}^N$-fibration  over
$M^s(n,d_2)\times \mathrm{Sym}^{d_1-d_2}(X)$, where 
$M^s(n,d_2)$ denotes the subspace of stable bundles of type $(n,d_2)$,
$\mathrm{Sym}^{d_1-d_2}(X)$ is the symmetric product, and the fiber dimension is
$N=n(d_1-d_2)-1$.

\quad $\bullet$ If $n_1>n_2$ then the moduli space 
$\mathcal{N}^s_\alpha(n_1,n_2,d_1,d_2)$ is birationally 
equivalent to a $\mathbb{P}^N$-fibration over 
$M^s(n_1-n_2,d_1-d_2) \times M^s(n_2,d_2)$, where the fiber dimension is
$N=n_2d_1-n_1d_2+n_2(n_1-n_2)(g-1)-1$.  

\quad $\bullet$ If $n_1<n_2$ then the moduli space 
$\mathcal{N}^s_\alpha(n_1,n_2,d_1,d_2)$ is birationally 
equivalent to a $\mathbb{P}^N$-fibration over 
$M^s(n_2-n_1,d_2-d_1) \times M^s(n_1,d_1)$, where 
the fiber dimension is $N=n_2d_1-n_1d_2+n_1(n_2-n_1)(g-1)-1$. 

\item[$(3)$] If $n_1\ne n_2$ then the moduli space 
$\mathcal{N}_{\alpha_M}(n_1,n_2,d_1,d_2)$ is 
non-empty and irreducible. Moreover 
\begin{displaymath}
\mathcal{N}_{\alpha_M}(n_1,n_2,d_1,d_2)\cong
\begin{cases}
M(n_2,d_2) \times M(n_1-n_2,d_1-d_2) &\text{if $n_1>n_2$}\\ 
M(n_1,d_1) \times M(n_2-n_1,d_2-d_1) &\text{if $n_1<n_2$.} 
\end{cases}
\end{displaymath}
\end{enumerate}
\end{theoremA}

Our strategy for studying the moduli spaces is similar in spirit to
the one used by Thaddeus \cite{thaddeus:1994}: basically it consists
in obtaining a good understanding of the moduli space for a particular
(large) value of $\alpha$ and then keeping track of how the moduli
space changes as $\alpha$ varies.  In the following we explain this in
more detail and outline the contents of the paper.

After recalling the basic facts about holomorphic triples in
Section~\ref{sec:stable-triples}, we go on to study extensions and
deformations of triples in Section~\ref{sec:extensions-of-triples}.
Here we show that the quasi-projective variety $\mathcal{N}_\alpha^s
\subseteq \mathcal{N}_\alpha$ corresponding to $\alpha$-stable triples
is smooth for all values of $\alpha$ greater than or equal to $2g-2$
(Theorem~\ref{thm:smoothdim}).

In Sections \ref{sec:chi-bounds} and
\ref{sec:crossing-critical-values} we examine how the moduli spaces
differ for values of $\alpha$ on opposite sides of a critical value.
If $\mathcal{N}_{\alpha_c^{\pm}}$ denote the moduli spaces for values
of $\alpha$ above and below a critical value $\alpha_c$, we denote the
loci along which they differ by $\mathcal{S}_{\alpha_c^{\pm}}$
respectively.  Our main result (Theorem~\ref{thm:codim}) is that for
all $\alpha\ge 2g-2$ the codimension of $\mathcal{S}_{\alpha_c^{\pm}}$
is strictly positive.  It follows that the number of irreducible
components of the spaces $\mathcal{N}^s_{\alpha}$ are the same for all
$\alpha$ satisfying $\alpha\ge 2g-2$ and $\alpha_m<\alpha <\alpha_M$.
In order to estimate the codimension of the
$\mathcal{S}_{\alpha_c^{\pm}}$ we need to estimate the dimension of
certain spaces of extensions of triples.  It is notable that this
requires us to consider objects more general than triples, namely the
\emph{holomorphic chains} studied in \cite{alvarez-garcia-prada:2001}.
The rather technical details are in Section~\ref{sec:chi-bounds}: the
main result is Proposition~\ref{prop:H-alpha-stable} which is then
used to deduce the key
Proposition~\ref{prop:hyper-1-vanishing-criterion}.

Next we turn to the question of understanding the moduli spaces
$\mathcal{N}_\alpha$ for large values of the parameter $\alpha$.
After obtaining some preliminary results in
Section~\ref{sect:special-alpha}, we consider the case of triples with
$n_1 \neq n_2$ in Section~\ref{sec:n_1-neq-n_2}.  
Let $\mathcal{N}_L$ denote the moduli space of $\alpha$-polystable
triples for $\alpha$ between $\alpha_M$ and the largest critical value
smaller than $\alpha_M$.  We show that
this `large $\alpha$' moduli space is birationally equivalent to a
$\mathbb{P}^N$-fibration over a product of moduli spaces of stable
bundles (Theorem~\ref{thm:largealpha}).  Combining this fact with our
codimension estimates we obtain our main results on non-emptiness and
irreducibility of the moduli spaces $\mathcal{N}_{\alpha}$ and
$\mathcal{N}^s_{\alpha}$; these appear as
Theorem~\ref{thm:irreducibility-moduli-stable-triples} and
Corollary~\ref{cor:gcd=1}.

In Section~\ref{sec:n_1=n_2} we obtain analogous results in the case
when $n_1 = n_2$.  Even though there is no upper limit to $\alpha$ in
this case, the moduli spaces do stabilize for $\alpha$ sufficiently
large (Theorem~\ref{thm:stabilization}) and hence it makes sense to
consider the large $\alpha$ moduli space $\mathcal{N}_L$ also in this
case.  The birational description of $\mathcal{N}_L$ is given in
Theorem~\ref{thm:moduli-n1=n2}, while the main results on
non-emptiness and irreducibility are in
Theorem~\ref{thm:irreducibility-moduli-stable-triples-n1=n2}. 

Finally, in Section~\ref {sec:reduction}, we go
back to the origins of the theory of holomorphic triples and apply
our results on moduli of triples to deduce the existence of
$\SU(2)$-invariant Hermitian--Einstein metrics on complex vector
bundles on $X \times \mathbb{P}^1$; equivalently, our results imply
the existence of stable vector bundles on $X \times \mathbb{P}^1$.

This paper and its companion \cite{bradlow-garcia-gothen:2002} form a
substantially revised version of the preprint
\cite{bradlow-garcia-gothen:2002:preprint}. The main results proved in
this paper were announced in the note
\cite{bradlow-garcia-gothen:2001}.  In that note we claim (without
proof) that for $\alpha \geq 2g-2$, the moduli spaces
$\mathcal{N}_\alpha$ are irreducible without imposing the conditions
in $(2)$ or $(3)$ of Theorem~A.  This is a reasonable conjecture,
which we hope to come back to in a future publication.

\subsection*{Acknowledgements.} We thank the Mathematics Departments
of the University of Illinois at Urbana-Champaign, the University
Aut{\'o}noma of Madrid and the University of Aarhus, the Department of
Pure Mathematics of the University of Porto, the Mathematical Sciences
Research Institute of Berkeley, the Mathematical Institute of the
University of Oxford, and the Erwin Schr\"odinger International
Institute for Mathematical Physics in Vienna for their hospitality 
during various stages of
this research.  We thank Ron Donagi, Tom\'as G\'omez,
Rafael Hern\'andez, Nigel Hitchin, Alastair King, 
Vicente Mu\~noz, Peter Newstead, and S. Ramanan for  many
insights and patient explanations.

\section{Definitions and basic facts}
\label{sec:stable-triples}

\subsection{Holomorphic triples and their moduli spaces}
\label{sec:triples-definitions}

Let $X$ be a compact Riemann surface (some of what follows is also true
also for a compact K\"ahler manifold 
\cite{garcia-prada:1994,alvarez-garcia-prada:2001}). 
Recall (\cite{bradlow-garcia-prada:1996} and \cite{garcia-prada:1994}) that
a \emph{holomorphic triple}  $T = (E_{1},E_{2},\phi)$ on $X$ consists of two 
holomorphic
vector bundles $E_{1}$ and $E_{2}$ on $X$ and a
holomorphic map $\phi \colon E_{2} \to E_{1}$.
A homomorphism from $T' = (E_1',E_2',\phi')$ 
to $T = (E_1,E_2,\phi)$   is a commutative diagram 
\begin{displaymath}
  \begin{CD}
    E_2' @>\phi'>> E_1' \\
    @VVV @VVV  \\
    E_2 @>\phi>> E_1,
  \end{CD}
\end{displaymath}
where the vertical arrows are holomorphic maps.
A triple $T'=(E_1',E_2',\phi')$ is a subtriple of $T = (E_1,E_2,\phi)$   
if the sheaf homomorphims  $E_1'\to E_1$ and $E_2'\to E_2$ 
are injective. A subtriple $T'\subset T$  is called {\em proper} if 
 $T'\neq 0 $ and $T'\neq T$.

\begin{definition}
For any $\alpha \in \R$ the \emph{$\alpha$-degree} and
\emph{$\alpha$-slope} of $T$ are 
defined to be
\begin{align*}
  \deg_{\alpha}(T)
  &= \deg(E_{1}) + \deg(E_{2}) + \alpha
  \rk(E_{2}), \\ 
  \mu_{\alpha}(T)
  &=
  \frac{\deg_{\alpha}(T)}
  {\rk(E_{1})+\rk(E_{2})} \\ 
  &= \mu(E_{1} \oplus E_{2}) +
  \alpha\frac{\rk(E_{2})}{\rk(E_{1})+
    \rk(E_{2})},
\end{align*}
where $\deg(E)$,  $\rk(E)$  and $\mu(E)=\deg(E)/\rk(E)$ are the 
degree, rank and slope of $E$, respectively.

We say  $T = (E_{1},E_{2},\phi)$ is
\emph{$\alpha$-stable} if
\begin{displaymath}
  \mu_{\alpha}(T')
  < \mu_{\alpha}(T)
\end{displaymath}
for any proper subtriple $T' = (E_{1}',E_{2}',\phi')$. 
Sometimes it is convenient to use
\begin{equation}\label{delta-stability}
\Delta_\alpha(T')=\mu_\alpha(T')-\mu_\alpha(T),
\end{equation}
in terms of which the $\alpha$-stability of $T$ is equivalent
to $\Delta_\alpha(T')<0$ for any proper subtriple $T'$.
We define  \emph{$\alpha$-semistability} by replacing the above 
strict inequality with a weak inequality. A triple is called
\emph{$\alpha$-polystable} if it is the direct sum of $\alpha$-stable
triples of the same $\alpha$-slope.
\end{definition}

Write $\mathbf{n}=(n_1,n_2)$ and
$\mathbf{d} = (d_1,d_2)$.  We denote by
\begin{displaymath}
  \mathcal{N}_\alpha
  = \mathcal{N}_\alpha(\mathbf{n},\mathbf{d})
  = \mathcal{N}_\alpha(n_1,n_2,d_1,d_2)
\end{displaymath}
the moduli space of $\alpha$-polystable triples $T =
(E_{1},E_{2},\phi)$ which have $\rk(E_i)=n_i$ and $\deg(E_i) = d_i$ for
$i=1,2$. The subspace of  $\alpha$-stable triples is denoted by
$ \mathcal{N}_\alpha^s$. We refer to 
$(\mathbf{n},\mathbf{d})=(n_1,n_2,d_1,d_2)$ as the {\em type} of the triple.

There are certain necessary conditions in order for   $\alpha$-semistable 
triples to exist.
Let $\mu_i=d_i/n_i$ for $i=1,2$. We define
\begin{align} 
 \alpha_m= &\mu_1-\mu_2, \label{alpha-bounds-m} \\
      \alpha_M = & (1+ \frac{n_1+n_2}{|n_1 - n_2|})(\mu_1 - \mu_2), \;\;
  n_1\neq n_2. \label{alpha-bounds-bigM} 
\end{align}
\begin{proposition}\cite[Theorem 6.1]{bradlow-garcia-prada:1996}
  \label{prop:alpha-range} 
The moduli space   $\mathcal{N}_\alpha(n_1,n_2,d_1,d_2)$ is a complex analytic
variety, which is projective when $\alpha$ is rational.
A necessary condition for $\mathcal{N}_\alpha(n_1,n_2,d_1,d_2)$ 
to be non-empty is 
\begin{enumerate}
\item[]
$0\leq \alpha_m \leq \alpha \leq \alpha_M$ \ \  if \ \  $n_1\neq n_2$,
\item[]
$0\leq \alpha_m \leq \alpha$ \ \    if  \ \ $n_1= n_2$.
\end{enumerate}
\end{proposition}

\begin{remark}
If  $\alpha_m=0$ and  $n_1\neq n_2$ then
$\alpha_m=\alpha_M=0$ and the moduli space
of $\alpha$ stable triples is empty unless $\alpha=0$.
\end{remark}

A direct  construction of these   moduli spaces 
has been given by Schmitt \cite{Schmitt} using geometric invariant theory.

Given a triple  $T=(E_1,E_2,\phi)$ one has the dual triple
$T^*=(E_2^*,E_1^*,\phi^*)$, where $E_i^*$ is the dual of $E_i$ and
$\phi^*$ is the transpose of $\phi$. The following  is not difficult 
to prove (\cite[Proposition 3.16]{bradlow-garcia-prada:1996}).
\begin{proposition}
\label{prop:duality}
The  $\alpha$-(semi)stability of $T$   is equivalent to the 
$\alpha$-(semi)stability of $T^*$.
The map $T\mapsto T^*$ defines a bijection
$$
\mathcal{N}_\alpha(n_1,n_2,d_1,d_2) = \mathcal{N}_\alpha(n_2,n_1,-d_2,-d_1),
$$
which is moreover an isomorphism.
\end{proposition}

This can be used to restrict  our study to $n_1\geq n_2$ and appeal to duality to deal with the
case $n_1<n_2$.

\subsection{Critical values}\label{sec:critical-values}

A holomorphic triple $T=(E_1,E_2,\phi)$ of type $(n_1,n_2,d_1,d_2)$ is strictly
$\alpha$-semistable if and only if it has a proper subtriple 
$T'=(E_1',E_2',\phi')$ such that $\mu_{\alpha}(T')= \mu_{\alpha}(T)$,
i.e.
\begin{equation}
  \label{eq:strict-alpha-ss}
  \mu(E'_1 \oplus E'_2) + \alpha \frac{n'_2}{n_1'+n_2'} 
  =  \mu(E_1 \oplus E_2) + \alpha \frac{n_2}{n_1+n_2}.
\end{equation}
There are two ways in which this can happen: The first one is if there
exists a subtriple $T'$ such that 
\begin{align*}
  \frac{n'_2}{n_1'+n_2'} &= \frac{n_2}{n_1+n_2},\;\; \mbox{and} \\
  \mu(E'_1 \oplus E'_2) &= \mu(E_1 \oplus E_2).
\end{align*}
In this case the terms containing $\alpha$ drop from
\eqref{eq:strict-alpha-ss} and $T$ is strictly
$\alpha$-semistable for all values of $\alpha$.  We refer to this
phenomenon as \emph{$\alpha$-independent semistability}. This cannot 
happen if $\GCD(n_2,n_1+n_2,d_1+d_2)=1$. The other way in which 
strict $\alpha$-semistability can happen is if equality holds in 
\eqref{eq:strict-alpha-ss} but 
\begin{equation}
  \label{eq:not-alpha-indep}
  \frac{n'_2}{n_1'+n_2'} \neq \frac{n_2}{n_1+n_2}.
\end{equation}
The values of 
$\alpha$ for which this
happens are called critical values. 

\begin{definition}
We say that $\alpha\in [\alpha_m,\infty)$ is a  \emph{critical value} if
 there exist integers $n'_1$, $n'_2$, $d'_1$ and $d'_2$
such that
\begin{displaymath}
  \frac{d'_1+d'_2}{n'_1+n'_2} + \alpha\frac{n'_2}{n'_1+n'_2}
  = \frac{d_1+d_2}{n_1+n_2} + \alpha\frac{n_2}{n_1+n_2},
\end{displaymath}
that is,
$$
\alpha=\frac{(n_1+n_2)(d_1'+d_2')-(n_1'+n_2')(d_1+d_2)}{n_1'n_2-n_1n_2'},
$$
with 
$0 \le n'_i \leq n_i$,  $(n'_1,n'_2,d'_1,d'_2) \neq (n_1,n_2,d_1,d_2)$,
$(n'_1,n'_2) \neq (0,0)$ and  $n_1'n_2\neq n_1n_2'$.
We say that $\alpha$ is {\em generic} if it is not critical. 
\end{definition}

\begin{proposition}\cite{bradlow-garcia-prada:1996}
\label{triples-critical-range}
Fix $(n_1,n_2,d_1,d_2)$.
\begin{itemize}
\item[(1)]
The critical values of $\alpha$
form a discrete subset of  $\alpha\in [\alpha_m,\infty)$,
where $\alpha_m$ is as in (\ref{alpha-bounds-m}).
\item[(2)] If $n_1\neq n_2$ the number of critical values is  finite and
  lies in the interval $[\alpha_m,\alpha_M]$, where $\alpha_M$ is as in 
(\ref{alpha-bounds-bigM}).
\item[(3)] The stability
criteria  for two values of  $\alpha$ lying between two consecutive 
critical values
are equivalent; thus the corresponding moduli spaces are isomorphic.
\item[(4)]  If $\alpha$ is generic and  $\GCD(n_2,n_1+n_2,d_1+d_2) = 1$,
 then $\alpha$-semistability is equivalent to $\alpha$-stability.
\end{itemize}
\end{proposition}

For the application of triples to   $\U(p,q)$-Higgs bundles 
(\cite{bradlow-garcia-gothen:2002}; see also the Introduction),   
it is important to have criteria to rule out
strict $\alpha$-semistability when $\alpha=2g-2$, where $g$ is the genus
of the surface. 
One such criterion, dealing actually with any  integral values of $\alpha$,
  is given by the following.
\begin{lemma}
\label{critical-integer}
Let $m$ be an integer such that $\GCD(n_1+n_2,d_1+d_2-mn_1)=1$. Then 
  \begin{itemize}
  \item[(1)] $\alpha=m$ is not a critical value,
  \item[(2)] there are no $\alpha$-independent semistable triples.
  \end{itemize}
\end{lemma}
\begin{proof}
To prove (1), suppose that $\alpha=m$ is a critical value. There exist
then a triple $T$ and a proper subtriple $T'$ so that 
$$
(d_1'+d_2'+mn_2')(n_1+n_2)=(d_1+d_2+m n_2)(n_1'+n_2').
$$
Thus $n_1+n_2$ divides 
$(d_1+d_2+ m n_2)(n_1'+n_2')$. But
 $n_1+n_2>n_1'+n_2'$, so  we get that $\GCD(n_1+n_2,d_1+d_2+m 
n_2)>1$. Writing $d_1+d_2+mn_2=d_1+d_2-mn_1+m(n_1+n_2)$, we see that
$\GCD(n_1+n_2,d_1+d_2-m n_1)>1$, in contradiction with the 
hypothesis. To prove (2), we  show that 
$\GCD(n_2,n_1+n_2,d_1+d_2)=1$, from which the  result follows by  (4) 
in Proposition \ref{triples-critical-range}. Suppose that 
$\GCD(n_2,n_1+n_2,d_1+d_2)\neq 1$. Then there is 
$(n_1',n_2',d_1',d_2')$ such that $ \frac{n_2'}{n_1'+n_2'}= 
\frac{n_2}{n_1+n_2}$ and $\frac{d_1'+d_2'}{n_1'+n_2'} 
= \frac{d_1+d_2}{n_1+n_2}$.  
It follows that 
$$ 
\frac{d_1'+d_2'-mn_1'}{n_1'+n_2'}=\frac{d_1+d_2-mn_1}{n_1+n_2},
$$
and hence $\GCD(n_1+n_2,d_1+d_2-mn_1)\neq 1$, in contradiction with 
the hypothesis. 
\end{proof}

\subsection{Vortex equations}\label{sec:vortex-equations}

There is a correspondence between stability and the existence of
solutions to certain  gauge-theoretic equations on a triple 
$T=(E_1,E_2,\phi)$, known as the {\em vortex equations} 
(\cite{bradlow-garcia-prada:1996} and \cite{garcia-prada:1994}).
The vortex equations
  \begin{equation}
    \label{eq:coupled-vortex}
  \begin{split}
    \sqrt{-1}\Lambda F(E_1) + \phi\phi^{*}
      &= \tau_{1} \Id_{E_1},  \\
    \sqrt{-1}\Lambda F(E_2)  - \phi^{*}\phi
      &= \tau_{2} \Id_{E_2}, 
  \end{split}
  \end{equation}
are equations for Hermitian metrics on $E_1$ and $E_2$.
Here $\Lambda$ is contraction by the K\"ahler form of a metric on $X$  
(normalized so that $\vol(X)=2\pi$), $F(E_i)$ is the curvature of the 
unique connection on $E_i$ compatible
with the Hermitian metric and the holomorphic structure of $E_i$, 
and $\tau_1$ and $\tau_2$ are real parameters satisfying 
$d_1+d_2=n_1\tau_1+n_2\tau_2$.  Here $\phi^\ast$ is the adjoint
of $\phi$ with respect to the Hermitian metrics.
One has the following.

\begin{theorem}\cite[Theorem 5.1]{bradlow-garcia-prada:1996}
\label{thm:vortices-hitchin-kobayashi}
A  solution to \eqref{eq:coupled-vortex} exists if and only if $T$ is
$\alpha$-polystable for $\alpha=\tau_1-\tau_2$.
\end{theorem}

Using the  vortex interpretation of the moduli space of triples  
one can easily  identify the  moduli space of triples for $\alpha=\alpha_m$.

\begin{proposition}\label{moduli-alpha_m}
A triple $T=(E_1,E_2,\phi)$ is $\alpha_m$-polystable
if and only if $\phi=0$ and $E_1$ and $E_2$ are polystable. We thus have
$$
\mathcal{N}_{\alpha_m}(n_1,n_2,d_1,d_2)\cong M(n_1,d_1) \times M(n_2,d_2),
$$
where $M(n_i,d_i)$ is the moduli space of semistable bundles of rank $n_i$
and degree $d_i$. 
\end{proposition}
\begin{proof}
Consider  equations \eqref{eq:coupled-vortex} on $T$.
If $\alpha=\alpha_m$ then $\tau_1=\mu_1$ and $\tau_2=\mu_2$ and hence 
in order to have solutions of \eqref{eq:coupled-vortex} 
we must have  $\phi=0$. In this 
case,  \eqref{eq:coupled-vortex} say that the  Hermitian metrics on 
$E_1$ and $E_2$ have constant central curvature. But this is 
equivalent to the polystability of $E_1$ and $E_2$ by the theorem of 
Narasimhan and Seshadri 
\cite{narasimhan-seshadri:1965}.
\end{proof}

\section{Extensions and deformations of triples}
\label{sec:extensions-of-triples}

In order to analyse the differences between the  moduli spaces 
$\mathcal{N}_\alpha$ as  $\alpha$ changes, as well as the   
smoothness properties of the moduli space for a given value of $\alpha$, 
we need to study the homological algebra of triples. This is done by 
considering the hypercohomology of  a certain complex of sheaves, in 
a similar way  to what is done in the study of infinitesimal deformations
by  Biswas and Ramanan \cite{biswas-ramanan:1994}. 
In fact, it is  a special case of the more general situation
considered  in \cite{gothen-king:2002}.

\subsection{Extensions}

Let $T'=(E'_1,E'_2,\phi')$ and $T''=(E''_1,E''_2,\phi'')$ be two triples
and, as usual, let
\begin{align*}
  (\mathbf{n}',\mathbf{d}') &= (n_{1}',n_{2}',d_{1}',d_{2}'), \\
  (\mathbf{n}'',\mathbf{d}'') &= (n_{1}'',n_{2}'',d_{1}'',d_{2}''),
\end{align*}
where $n'_i = \rk(E'_i)$, $n''_i = \rk(E''_i)$, $d'_i = \deg(E'_i)$
and $d''_i = \deg(E''_i)$.  Let $\Hom(T'',T')$ denote the linear  
space of homomorphisms from  $T''$ to
$T'$, and let $\Ext^1(T'',T')$  denote the linear space of equivalence
classes of extensions of the form
\begin{displaymath}
  0 \lto T' \lto T \lto T'' \lto 0,
\end{displaymath}
where by this we mean a commutative  diagram
  $$
  \begin{CD}
  0@>>>E_1'@>>>E_1@>>> E_1''@>>>0\\
  @.@A\phi' AA@A \phi AA@A \phi'' AA\\
  0@>>>E'_2@>>>E_2@>>>E_2''@>>>0.
  \end{CD}
$$
Hence, to analyse $\Ext^1(T'',T')$ one considers the
complex of sheaves 
\begin{equation}
  \label{eq:extension-complex}
    C^{\bullet}(T'',T') \colon {E_{1}''}^{*} \otimes E_{1}' \oplus
  {E_{2}''}^{*} \otimes E_{2}'
  \overset{c}{\lto}
  {E_{2}''}^{*} \otimes E_{1}',
\end{equation} 
where the map $c$ is defined by
\begin{displaymath}
  c(\psi_{1},\psi_{2}) = \phi'\psi_{2} - \psi_{1}\phi''.
\end{displaymath}
\begin{proposition}
  \label{prop:hyper-equals-hom}
  There are natural isomorphisms
  \begin{align*}
    \Hom(T'',T') &\cong \HH^{0}(C^{\bullet}(T'',T')), \\
    \Ext^{1}(T'',T') &\cong \HH^{1}(C^{\bullet}(T'',T')),
  \end{align*}
and a long exact sequence associated to the complex
$C^{\bullet}(T'',T')$:
\begin{equation}
  \label{eq:long-exact-extension-complex}
\begin{array}{ccccccc}
  0 &\lto \mathbb{H}^0(C^{\bullet}(T'',T')) &
  \lto & H^0({E_{1}''}^{*} \otimes E_{1}' \oplus {E_{2}''}^{*} \otimes E_{2}')
  & \lto &  H^0({E_{2}''}^{*} \otimes E_{1}') \\
    &  \lto \mathbb{H}^1(C^{\bullet}(T'',T')) &
  \lto &  H^1({E_{1}''}^{*} \otimes E_{1}' \oplus {E_{2}''}^{*} \otimes E_{2}')
&  \lto & H^1({E_{2}''}^{*} \otimes E_{1}') \\
&   \lto \mathbb{H}^2(C^{\bullet}(T'',T')) & \lto & 0. & & 
\end{array}
\end{equation}
\end{proposition}
\begin{proof}
The proof is omitted since it is very similar to that given in 
\cite{biswas-ramanan:1994} in the  study of  deformations, and
it is a  special case of a much 
more general result proved in  \cite{gothen-king:2002}. 
\end{proof}
We introduce the following notation: 
\begin{align}
  h^{i}(T'',T') &= \dim\HH^{i}(C^{\bullet}(T'',T')), \notag \\
  \chi(T'',T') &= h^0(T'',T') - h^1(T'',T') + h^2(T'',T'). \label{euler}
\end{align}

\begin{proposition}
  \label{prop:chi(T'',T')}
  For any holomorphic triples $T'$ and $T''$ we have
  \begin{align*}
    \chi(T'',T') &= \chi({E_{1}''}^{*} \otimes E_{1}')
    + \chi({E_{2}''}^{*} \otimes E_{2}')
    - \chi({E_{2}''}^{*} \otimes E_{1}')  \\
    &= (1-g)(n''_1 n'_1 + n''_2 n'_2 - n''_2 n'_1) \\
    & \quad + n''_1 d'_1 - n'_1 d''_1
    + n''_2 d'_2 - n'_2 d''_2
    - n''_2 d'_1 + n'_1 d''_2,
  \end{align*}
where $\chi(E)=\dim H^0(E) - \dim H^1(E)$ is the Euler characteristic
of $E$.
\end{proposition}
\begin{proof}
  Immediate from the long exact sequence
  \eqref{eq:long-exact-extension-complex} and the Riemann--Roch formula.
\end{proof}
\begin{corollary}
  \label{cor:chi-relation}
  For any extension $0 \to T' \to T \to T'' \to 0$ of triples,
  \begin{displaymath}
    \chi(T,T) = \chi(T',T') + \chi(T'',T'') + \chi(T'',T') +
    \chi(T',T''). 
  \end{displaymath}
  \qed
\end{corollary}
\begin{remark} Proposition \ref{prop:chi(T'',T')} shows that 
$\chi(T'',T')$ depends only  on the topological invariants 
$(\mathbf{n}',\mathbf{d}')$ and 
$(\mathbf{n}'',\mathbf{d}'')$ of $T'$ and $T''$.  Whenever convenient
we shall therefore use the notation
\begin{displaymath}
  \chi(\mathbf{n}'',\mathbf{d}'',\mathbf{n}',\mathbf{d}')
  = \chi(T'',T').
\end{displaymath}
\end{remark}

\subsection{Vanishing of $\HH^{0}$ and $\HH^{2}$}

 The following vanishing results play a central role in our study.
 \begin{proposition}
   \label{prop:h0-vanishing}
   Suppose that $T'$ and $T''$ are $\alpha$-semistable. 
   
   \begin{itemize}
\item[(1)]
   If $\mu_\alpha(T')<\mu_\alpha (T'')$ then  
 $\HH^{0}(C^{\bullet}(T'',T')) \cong 0$.
\item[(2)] 
If $\mu_\alpha(T')=\mu_\alpha (T'')$ and $T''$ is
   $\alpha$-stable, then
  $$
     \HH^{0}(C^{\bullet}(T'',T')) \cong 
     \begin{cases}
       \C \quad &\text{if $T' \cong T''$} \\
       0 \quad &\text{if $T' \not\cong T''$}. 
     \end{cases}
 $$
 \end{itemize}
   \end{proposition}
 \begin{proof} By Proposition \ref{prop:hyper-equals-hom} we can 
 identify $\HH^{0}(C^{\bullet}(T'',T'))$ with $\Hom(T'',T')$. The 
 statements (1) and (2) are thus the direct analogs for triples of the 
 same results for semistable bundles. The proof is identical. Suppose 
 that $h:T''\rightarrow T'$ is a non-trivial homomorphism 
 of triples. If 
 $T'=(E_1',E_2',\Phi')$ and $T''=(E_1'',E_2'',\Phi'')$ then $h$ is 
 given by a pair of holomorphic  maps $u_i:E_i''\rightarrow E_i''$ for 
 $i=1,2$ such that $\Phi'\circ u_2=u_1\circ\Phi''$. We can thus define
 subtriples of $T''$ and $T'$ respectively by 
 $T_N=(\ker(u_1),\ker(u_2),\Phi'')$ and $T_I=(\im(u_1),\im(u_2),\Phi')$,  where in $T_I$, it is in general necessary to take the saturations of 
 the image $\im(u_1)$ and $\im(u_2)$. By the semistability conditions, 
 we get
 $$
 \mu_\alpha(T_N)\le\mu_{\alpha}(T'')\le\mu_{\alpha}(T_I)\le\mu_\alpha(T').
 $$ 
 The conclusions follow directly from this.
 \end{proof}
\begin{proposition}
  \label{prop:h2-vanishing}
Suppose that the  triples $T'$ and $T''$ are $\alpha$-semistable and 
satisfy
$\mu_\alpha(T')=\mu_\alpha (T'')$. Then 
\begin{itemize}
\item[(1)] $\HH^{2}(C^{\bullet}(T'',T')) = 0$ whenever $\alpha > 2g-2$.
\item[(2)] If one of $T'$, $T''$ is   $\alpha+\epsilon$-stable for some
$\epsilon\ge 0$, then $\HH^{2}(C^{\bullet}(T'',T')) = 0$ whenever
$\alpha \geq2g-2$.
\end{itemize}
\end{proposition}
\begin{proof}
{}From (\ref{eq:long-exact-extension-complex}) it is clear that the 
vanishing of 
$\HH^{2}(C^{\bullet}(T'',T'))$ is equivalent to the surjectivity of the map
$$
  H^1({E_{1}''}^{*} \otimes E_{1}' \oplus {E_{2}''}^{*} \otimes E_{2}')
  \lto  H^1({E_{2}''}^{*} \otimes E_{1}').
$$
By Serre duality this is equivalent to the injectivity of the map
\begin{equation}
  \label{eq:petri}
\begin{array}{ccc}
H^0({E_1'}^{*} \otimes E_{2}''\otimes K) &   \overset{P}{\lto} &
H^0({E_{1}'}^{*} \otimes E_{1}''\otimes K) \oplus H^0({E_{2}'}^{*}
\otimes E_{2}''\otimes K)\\
\psi &\longmapsto & ((\phi''\otimes \Id)\circ \psi, \psi\circ \phi').
\end{array}
\end{equation}
\noindent {\it Proof of (1):} Suppose that $P$ is not injective. Then there is 
a non-trivial homomorphism $\psi:E_1' \to E_2''\otimes K$ in 
$\ker P$. Let $I=\im\psi$ and $N=\ker\psi$. Since 
$(\phi''\otimes \Id_K)\circ\psi=0$, $I\subset\ker\phi''$ and hence 
$T_I''=(0, I\otimes K^*, 0)$ is a proper subtriple of  $T''$. Similarly,
the fact that $\psi\circ \phi'=0$ implies that $\im\phi'\subset N$ and thus
 $T_N'=(\ker\psi,E_2',\phi')$ is a proper subtriple of $T'$.  Let $k=\rk(N)$ and 
$l=\deg(N)$. Then, from the exact sequence
\begin{displaymath}
  0 \lto N \lto E_1' \lto I \lto 0
\end{displaymath}
we see that $\rk(I)=n_1'-k$ and $\deg(I)=d_1'-l$. Hence
\begin{align*}
 \mu_\alpha(T_N') &= \frac{l+d_2'}{k+n_2'}+\alpha \frac{n_2'}{k+n_2'}, \\
\mu_\alpha(T_I'') &= \frac{d_1'-l}{n_1'-k}+ 2-2g +\alpha.
\end{align*}
Adding these two expressions, and clearing denominators we see that
\begin{displaymath}
d_1'+d_2'+(n_1'-k)(2-2g)+\alpha(n_1'+n_2'-k)= 
(k+n_2')\mu_\alpha(T_N')+ (n_1'-k)\mu_\alpha (T_I'').
\end{displaymath}
But  $\mu_\alpha(T_N')\leq \mu_\alpha(T')$, $\mu_\alpha(T_I'')\leq
 \mu_\alpha(T'')$ and $\mu_\alpha(T')=\mu_\alpha(T'')$.
{}From this we obtain that 
\begin{equation}\label{eqn:alpha-est}
d_1'+d_2'+(n_1'-k)(2-2g)+\alpha(n_1'+n_2'-k)\leq d_1'+d_2'+\alpha 
n_2', 
\end{equation}
and hence
\begin{displaymath}
\alpha(n_1'-k)\leq(n_1'-k)(2g-2).
\end{displaymath}
Since  $n_1'-k>0$ we get  that  $\alpha\leq 2g-2$. Hence $P$ must be 
injective if the hypotheses of the part (1) of the proposition are 
satisfied. 

\noindent {\it Proof of (2):} Suppose that $T''$ is $\alpha+\epsilon$-stable 
for some $\epsilon\ge 0$. It follows that 
$\mu_{\alpha+\epsilon}(T_I'')< 
\mu_{\alpha+\epsilon}(T'')$, i.e
$$
\mu_{\alpha}(T_I'')-
\mu_{\alpha}(T'')<\epsilon(\frac{n_2''}{n_1''+n_2''}-1)\le 0.
$$
Thus, following exactly the same argument as in the proof of (1), we 
get a strict inequality in (\ref{eqn:alpha-est}).  We conclude that 
that if P is not injective then $\alpha< 2g-2$, i.e.\ if $\alpha\ge 
2g-2$ then $P$ must be injective.
If  $T'$ is $\alpha+\epsilon$-stable for some $\epsilon\ge 0$ then we 
get that
$$
\mu_{\alpha}(T_N')-
\mu_{\alpha}(T')<\epsilon(\frac{n_2'}{n_1'+n_2'}-\frac{n_2'}{k+n_2'})\le 0.
$$
The rest of the argument is the same as in the case that $T''$ is 
$\alpha+\epsilon$-stable.
\end{proof}
\begin{corollary}\label{dimension-ext1}
Let  $T'$ and $T''$ be  $\alpha$-semistable triples  
with  $\mu_\alpha(T')=\mu_\alpha (T'')$, and $\alpha> 2g-2$.
Then
 $$
\dim \Ext^1(T'',T')=    h^0(T'',T')  - \chi(T'',T').
$$
The same holds for $\alpha\geq 2g-2$ if in addition $T'$ or $T''$ is 
$\alpha+\epsilon$-stable for some $\epsilon\ge 0$.
\end{corollary}
\begin{proof}
It follows from Proposition \ref{prop:hyper-equals-hom} and (\ref{euler})
that
\begin{equation}
\dim \Ext^1(T'',T')=h^0(T'',T')+h^2(T'',T')- \chi(T'',T').
\label{dim-ext1}
\end{equation}
The result follows immediately from this and the vanishing of 
$h^2(T'',T')$ given by  Proposition \ref{prop:h2-vanishing}.
\end{proof}

\subsection{Deformation theory for triples}

Since the  space of infinitesimal deformations of $T$ is
isomorphic to $\HH^{1}(C^{\bullet}(T,T))$, the
considerations of the previous sections also apply to studying deformations of a holomorphic
triple $T$.  To be precise, one has the following.
\begin{theorem}\label{thm:smoothdim}
Let $T=(E_1,E_2,\phi)$ be an $\alpha$-stable triple of type 
$(n_1,n_2,d_1,d_2)$.
\begin{itemize}
\item[(1)] The Zariski tangent space at the point defined by $T$
in the moduli space of stable triples  is isomorphic
to $\HH^{1}(C^{\bullet}(T,T))$.
\item[(2)]
If $\HH^{2}(C^{\bullet}(T,T))= 0$, then the moduli space of $\alpha$-stable
triples is smooth in  a neighbourhood of the point defined by $T$.
\item[(3)]
$\HH^{2}(C^{\bullet}(T,T))= 0$ if and only if the homomorphism
$$
  H^1(E_1^* \otimes E_1 \oplus E_2^* \otimes E_2)
  \lto  H^1(E_2^* \otimes E_1)
$$
in the corresponding long exact sequence is surjective.
\item[(4)] At a smooth point $T\in \mathcal{N}^s_\alpha(n_1,n_2,d_1,d_2)$
the dimension of the moduli space of $\alpha$-stable triples is
\begin{align}
  \dim \mathcal{N}^s_\alpha(n_1,n_2,d_1,d_2)
  &= h^{1}(T,T) = 1 - \chi(T,T) \notag \\
  &= (g-1)(n_1^2 + n_2^2 - n_1 n_2) - n_1 d_2 + n_2 d_1 + 1.
    \label{eq:dim-triples}
\end{align}
\item[(5)] If $\phi$ is injective or surjective then $T=(E_1,E_2,\phi)$
defines a smooth point in the moduli space.
\item[(6)] If $\alpha\geq 2g-2$, then $T$ defines a smooth point
in the moduli space, and hence $\mathcal{N}^s_\alpha(n_1,n_2,d_1,d_2)$ is 
smooth.
\end{itemize}
\end{theorem}
\begin{proof}
Statements (1) and (2) follow from  Theorems 2.3 and 
3.1 in \cite{biswas-ramanan:1994}, respectively.
An indirect proof of (1) and (2),  exploiting the correspondence between
triples on $X$ and stable bundles on $X\times \mathbb{P}^1$ (see
Section \ref{sec:reduction}) also follows  from  
\cite{bradlow-garcia-prada:1996}.
Statement (3) follows from the long exact sequence (\ref{eq:long-exact-extension-complex}) 
with $T=T'=T''$. (4) follows from (1), (2) and 
Propositions~\ref{prop:chi(T'',T')} and \ref{dimension-ext1}.
(5) is proved in \cite[Proposition 6.3]{bradlow-garcia-prada:1996}. (6) is a 
consequence of Proposition~\ref{prop:h2-vanishing}.
\end{proof}

\section{Bounds for $\chi$}
\label{sec:chi-bounds}
In our approach to the study of how the moduli spaces of triples 
vary with the parameter, it is of crucial importance to be able 
to estimate the
Euler characteristics $\chi(T'',T')=
\chi(\mathbf{n}'',\mathbf{d}'',\mathbf{n}',\mathbf{d}')$ when $T'$ and
$T''$ are polystable triples with the same $\alpha$-slope.  The basic
idea is to identify $\chi(T'',T')$ as a hypercohomology Euler
characteristic for the complex $C^{\bullet}(T'',T')$ defined in
(\ref{eq:extension-complex}) and to notice that the complex is itself a
holomorphic triple. As such it ought to satisfy a stability condition
induced from the stability condition of $T'$ and $T''$.  In principle, 
a way to obtain the stability condition for
$C^{\bullet}(T'',T')$ should be provided by the correspondence between the
stability of the holomorphic triples and the existence of solutions to
the vortex equations given by Theorem
\ref{thm:vortices-hitchin-kobayashi}.  However, there seem to be no
simple way to construct a solution to the vortex
equations for $C^{\bullet}(T'',T')$ from solutions on $T'$ and $T''$.
Instead we consider slightly more general objects than
triples, known as \emph{holomorphic chains}.  These are studied in
\cite{alvarez-garcia-prada:2001}.

\subsection{Holomorphic chains}\label{sec:holomorphic-chains}

A \emph{holomorphic chain} is a diagram 
\begin{displaymath}
  \mathcal{C} \colon E_{m} \overset{\phi_m}{\lto} E_{m-1}
  \overset{\phi_{m-1}}{\lto} \cdots
  \overset{\phi_{1}}{\lto} E_{0},
\end{displaymath}
where each $E_i$ is a holomorphic vector bundle and $\phi_{i} \colon E_{i} \to 
E_{i-1}$ is a holomorphic map.  Let 
\begin{align*}
  \mu(\mathcal{C}) &= \mu(E_{0} \oplus \cdots \oplus E_{m}), \\
  \lambda_{i}(\mathcal{C}) &=   \frac{\rk(E_{i})}
       {\sum_{i=0}^{m}\rk(E_{i})}, \quad i=0,\ldots, m.
\end{align*}
For $\balpha = (\alpha_1,\ldots,\alpha_m) \in \R^{m}$, the 
\emph{$\balpha$-slope} of $\mathcal{C}$ is defined to be
\begin{displaymath}
  \mu_{\balpha}(\mathcal{C}) = \mu(\mathcal{C}) +
  \sum_{i=1}^{m}\alpha_{i}\lambda_{i}(\mathcal{C}).
\end{displaymath}
The notion of $\balpha$-stability 
is defined via the standard $\balpha$-slope condition on 
subchains, that is, for any holomorphic subchain 
$\mathcal{C}'\subset \mathcal{C}$ we must have 
$\mu_{\balpha}(\mathcal{C}') <\mu_{\balpha}(\mathcal{C})$.
Semistability and polystability are defined as usual.
A holomorphic triple is a holomorphic chain of length 2, and 
the stability notions coincide, taking $\balpha=(\alpha)$. 
As for triples, there are natural gauge-theoretic equations for holomorphic chains, 
which we now describe.  Define $\btau = (\tau_{0}, \ldots, \tau_{m}) 
\in \R^{m+1}$ by
\begin{equation}
  \label{eq:tau-alpha}
  \tau_{i} = \mu_{\alpha}(\mathcal{C}) - \alpha_{i}, \quad i=0,\ldots, m,
\end{equation}
where we make the convention $\alpha_0 = 0$.  Then $\balpha$ can be 
recovered from $\btau$ by 
\begin{equation}
  \label{eq:alpha-tau}
  \alpha_{i} = \tau_{0} - \tau_{i}, \quad i=0,\ldots, m.
\end{equation}
The 
\emph{$\btau$-vortex equations} 
\begin{displaymath}
  \sqrt{-1}\Lambda F(E_{i}) + \phi_{i+1} \phi_{i+1}^* - \phi_{i}^*\phi_{i}
  = \tau_{i} \Id_{E_{i}}, \quad i=0,\ldots, m,
\end{displaymath}
are equations for Hermitian metrics on $E_{0}, \ldots, E_{m}$.  Here,
as in (\ref{eq:coupled-vortex}), 
$F(E_{i})$ is the curvature of the Hermitian connection on $E_{i}$,
$\Lambda$ is contraction with the K\"ahler form and $\vol(X) = 2\pi$.
By convention $\phi_0 = \phi_{m+1}=0$.
One has the generalization of Theorem \ref{thm:vortices-hitchin-kobayashi}
to the  case of holomorphic chains.
\begin{theorem} \cite[Theorem 3.4]{alvarez-garcia-prada:2001}
  \label{thm:chain-hitchin-kobayashi}
  A holomorphic chain $\mathcal{C}$ is $\balpha$-polystable if and
  only if the $\btau$-vortex equations have a solution, where
  $\balpha$ and $\btau$ are related by \eqref{eq:tau-alpha}.
\end{theorem}

\subsection{A length 3 holomorphic chain}
Let $T'=(E'_1,E'_2,\phi')$ and $T''=(E''_1,E''_2,\phi'')$ be two triples.
Let us  consider the length $3$ holomorphic chain 
\begin{equation}
  \label{eq:H-chain}
  \widetilde{C^{\bullet}}(T'',T') \colon
  {E_{1}''}^{*} \otimes E_{2}'
  \xrightarrow{a_2}
  {E_{1}''}^{*} \otimes E_{1}' \oplus {E_{2}''}^{*} \otimes E_{2}'
  \xrightarrow{a_1}
  {E_{2}''}^{*} \otimes E_{1}' ,
\end{equation}
where 
\begin{align*}
  a_{2}(\psi) &= (\phi'\psi, -\psi\phi''), \\
  a_{1}(\psi_{1},\psi_{2}) &= \phi'\psi_{2} - \psi_{1}\phi''.
\end{align*}
We shall sometimes write this chain briefly as 
\begin{displaymath}
  \widetilde{C^{\bullet}}(T'',T') \colon
  C_{2} \xrightarrow{a_{2}} C_{1} \xrightarrow{a_{1}} C_{0}.
\end{displaymath}
Note that the last two terms of $\widetilde{C^{\bullet}}(T'',T')$ coincide
with the complex $C^{\bullet}(T'',T')$.  Note also that 
$\widetilde{C^{\bullet}}(T'',T')$ is not in general a complex.  Our goal
in this section is to prove, using Theorem \ref{thm:chain-hitchin-kobayashi}, that if $T'$ and $T''$ are 
$\alpha$-polystable then $\widetilde{C^{\bullet}}(T'',T')$ is
$\balpha$-polystable for a suitable choice of $\balpha$.
\begin{lemma}
  \label{lemma:chain-vortex-solution}
  Let $T'$ and $T''$ be holomorphic triples and suppose we have
  solutions to the $(\tau_{1}',\tau_{2}')$-vortex equations on $T'$
  and the $(\tau_{1}'',\tau_{2}'')$-vortex equations on $T''$, such
  that $\tau_{1}' - \tau_{1}''= \tau_{2}' - \tau_{2}''$.
  Then the induced Hermitian metric on $\widetilde{C^{\bullet}}(T'',T')$
  satisfies the chain vortex equations
  \begin{align}
    \label{eq:chain-vortex-1}
    \sqrt{-1}\Lambda F(C_{0}) + a_{1}{a_{1}}^{*}
      &= \tilde{\tau}_{0} \Id_{C_{0}}, \\
    \label{eq:chain-vortex-2}
    \sqrt{-1}\Lambda F(C_{1}) + a_{2}{a_{2}}^{*} - {a_{1}}^{*}a_{1}
      &= \tilde{\tau}_{1} \Id_{C_{1}}, \\
    \label{eq:chain-vortex-3}
    \sqrt{-1}\Lambda F(C_{2}) - {a_{2}}^{*}a_{2}
      &= \tilde\tau_{2} \Id_{C_{2}},
  \end{align}
  for $\btau = (\tilde{\tau}_{0},\tilde{\tau}_{1},\tilde{\tau}_{2})$
  given by
  \begin{align*}
    \tilde{\tau}_{0} &= \tau_{1}' - \tau_{2}'', \\
    \tilde{\tau}_{1} &= \tau_{1}' - \tau_{1}''
                      = \tau_{2}' - \tau_{2}'', \\
    \tilde{\tau}_{2} &= \tau_{2}' - \tau_{1}''. \\
  \end{align*}
\end{lemma}
\begin{proof}
  We shall only show that the induced Hermitian metric satisfies
  \eqref{eq:chain-vortex-2}, since the proofs that it satisfies the
  two remaining equations are similar (but simpler).
  
  The  vortex equations for $T'$ and $T''$ are
  \begin{align*}
    \sqrt{-1}\Lambda F(E_{1}')
      + \phi'{\phi'}^{*}&=\tau_{1}' \Id_{E_{1}'}, &
      \sqrt{-1}\Lambda F(E_{1}'')
      + \phi''{\phi''}^{*}&=\tau_{1}'' \Id_{E_{1}''}, \\
      \sqrt{-1}\Lambda F(E_{2}')
      - {\phi'}^{*}\phi'&=\tau_{2}' \Id_{E_{2}'}, &
      \sqrt{-1}\Lambda F(E_{2}'')
      - {\phi''}^{*}\phi''&= \tau_{2}'' \Id_{E_{2}''}. \\
  \end{align*}
  We shall write the left hand side of \eqref{eq:chain-vortex-2} in
  terms of these known data of the triples $T'$ and $T''$.
  First, we note that
  \begin{displaymath}
    F({E_{i}'}^{*}) = -F(E_{i}')^{t}, \quad i=1,2,
  \end{displaymath}
  and similarly for $F({E_{i}''}^{*})$.  Hence
  \begin{align}
    F(C_{1}) &= F({E_{1}''}^{*} \otimes E_{1}'
      \oplus {E_{2}''}^{*} \otimes E_{2}') \notag \\
    &= \bigl(
      F({E_{1}''}^{*}) \otimes \Id + \Id \otimes F(E_{1}'),\
      F({E_{2}''}^{*}) \otimes \Id + \Id \otimes F(E_{2}')
      \bigr) \notag \\
    &= \bigl(
      -F({E_{1}''})^{t} \otimes \Id + \Id \otimes F(E_{1}'),\
      -F({E_{2}''})^{t} \otimes \Id + \Id \otimes F(E_{2}')
      \bigr). \label{eq:F(C_1)}
  \end{align}
  Next we calculate $a_{1}^{*}$: note that for $\xi \otimes x \in
  C_{0}$ and $(\eta_1 \otimes y_1, \eta_2 \otimes y_2) \in C_{1}$ we
  have
  \begin{align*}
    \bigl\langle a_{1}^{*}(\xi \otimes x),\
    (\eta_1 \otimes y_1, \eta_2 \otimes y_2) \bigr\rangle_{C_{1}}
    &=\bigl\langle \xi \otimes x,\
    a_{1}(\eta_1 \otimes y_1,\eta_2 \otimes y_2) \bigr\rangle_{C_{0}} \\
    &\hspace{-1cm}=\bigl\langle \xi \otimes x,\ - \eta_1\phi'' \otimes y_1
    + \eta_2 \otimes \phi'(y_2) \bigr\rangle_{C_{0}} \\
    &\hspace{-1cm}=\bigl\langle \xi \otimes x,\ - {\phi''}^{t}(\eta_1)
    \otimes y_1 + \eta_2 \otimes \phi'(y_2) \bigr\rangle_{C_{0}} \\
    &\hspace{-1cm}= -\bigl\langle \xi, 
    {\phi''}^{t}(\eta_1) \bigr\rangle_{{E_{2}''}^{*}}
    \ \bigl\langle x, y_1 \bigr\rangle_{E_{1}'} 
    \ +\ \bigl\langle \xi,\eta_2 \bigr\rangle_{{E_{2}''}^{*}}
    \ \bigl\langle x,\phi'(y_2) \bigr\rangle_{E_{1}'} \\
    &\hspace{-1cm}= -\bigl\langle {\phi''}^{t*}(\xi), 
    \eta_1 \bigr\rangle_{{E_{1}''}^{*}}
    \ \bigl\langle x, y_1 \bigr\rangle_{E_{1}'} 
    \ +\ \bigl\langle \xi,\eta_2 \bigr\rangle_{{E_{2}''}^{*}}
    \ \bigl\langle {\phi'}^{*}(x),y_2 \bigr\rangle_{E_{2}'} \\
    &\hspace{-1cm}=\bigl\langle (-{\phi''}^{t*}(\xi) \otimes x,
    \xi \otimes {\phi'}^{*}(x)),\
    (\eta_1 \otimes y_1, \eta_2 \otimes y_2) \bigr\rangle_{C_{1}}.
   \end{align*}
  Hence,
  \begin{equation}
    \label{eq:a_1^*}
    a_{1}^{*}(\xi \otimes x) = \bigl(-{\phi''}^{t*}(\xi) \otimes x,\
    \xi \otimes {\phi'}^{*}(x)\bigr).
  \end{equation}
  Similarly, to calculate $a_{2}^{*}$ consider $\xi \otimes x \in
  C_{2}$ and $(\eta_1 \otimes y_1, \eta_2 \otimes y_2) \in C_{1}$.
  Then
  \begin{equation}
    \label{eq:a_2^*}
    a_{2}^{*}(\eta_1 \otimes y_1, \eta_2 \otimes y_2)
    = \eta_1 \otimes {\phi'}^{*}(y_1)
    -{\phi''}^{t*}(\eta_2) \otimes y_2.
  \end{equation}
  Using \eqref{eq:a_2^*} and \eqref{eq:a_1^*} we can now calculate for
  $(\eta_1 \otimes y_1, \eta_2 \otimes y_2) \in C_{1}$:
  \begin{multline}
    \label{eq:a_2a_2^*}
    a_{2}a_{2}^{*}(\eta_1 \otimes y_1, \eta_2 \otimes y_2) \\
    = \bigl(\eta_1 \otimes \phi'{\phi'}^{*}(y_1)
    -{\phi''}^{t*}(\eta_2) \otimes \phi'(y_2),\
    -{\phi''}^{t}(\eta_1) \otimes {\phi'}^{*}(y_1)
    +{\phi''}^{t}{\phi''}^{t*}(\eta_2) \otimes y_2 \bigr),
  \end{multline}
  and
  \begin{multline}
    \label{eq:a_1^*a_1}
    a_{1}^{*}a_{1}(\eta_1 \otimes y_1, \eta_2 \otimes y_2) \\
    =\bigl({\phi''}^{t*}{\phi''}^{t}(\eta_1) \otimes y_1
    -{\phi''}^{t*}(\eta_2) \otimes \phi'(y_2),\
    -{\phi''}^{t}(\eta_1) \otimes {\phi'}^{*}(y_1)
    +\eta_2 \otimes {\phi'}^{*}\phi'(y_2) \bigr).
  \end{multline}
  Putting together \eqref{eq:F(C_1)}, \eqref{eq:a_2a_2^*} and
  \eqref{eq:a_1^*a_1} we finally obtain
  \begin{multline}
    \label{eq:chain-1-lhs}
    \bigl(\sqrt{-1}\Lambda F(C_{1}) + a_{2}{a_{2}}^{*} -
    {a_{1}}^{*}a_{1}\bigr)
    (\eta_1 \otimes y_1, \eta_2 \otimes y_2) \\
    = \Bigl( \eta_1 \otimes \bigl(\sqrt{-1}\Lambda F(E_{1}') +
    \phi'{\phi'}^{*} \bigr)(y_1)
    + \bigl(-\sqrt{-1}\Lambda F(E_{1}'')^{t} +
    {\phi''}^{t*}{\phi''}^{t} \bigr)(\eta_1) \otimes y_1, \\
    \eta_2 \otimes \bigl(\sqrt{-1}\Lambda F(E_{2}')
    - {\phi'}^{*}\phi'\bigr)(y_2)
    - \bigl(\sqrt{-1}\Lambda F(E_{2}'')^{t} -
    {\phi''}^{t}{\phi''}^{t*} \bigr)(\eta_2) \otimes y_2 \Bigr).
  \end{multline}
  Notice that the unpleasant mixed term
  $\bigl(-{\phi''}^{t*}(\eta_2) \otimes \phi'(y_2),
  -{\phi''}^{t}(\eta_1) \otimes {\phi'}^{*}(y_1) \bigr)$ appears both
  in $a_{1}^{*}a_{1}$ and $a_{2}a_{2}^{*}$ and therefore cancels.
  This would not have been the case if we had considered the 
  vortex equations on the triple $C^{\bullet}(T'',T')$ and is the
  reason why we must consider the chain $\widetilde{C^{\bullet}}(T'',T')$.
  Combining \eqref{eq:chain-1-lhs} with the  vortex equations
  (or their transposes) for the triples $T'$ and $T''$ we get
  \begin{multline}
    \bigl(\sqrt{-1}\Lambda F(C_{1}) + a_{2}{a_{2}}^{*} -
    {a_{1}}^{*}a_{1}\bigr)
    (\eta_1 \otimes y_1, \eta_2 \otimes y_2) \\
    = \bigl((\tau_{1}'-\tau_{1}'')\eta_1 \otimes y_1,
    (\tau_{2}'-\tau_{2}'')\eta_2 \otimes y_2 \bigr).
  \end{multline}
  Since $\tau_{1}'-\tau_{1}'' = \tau_{2}'-\tau_{2}''$ this
  concludes the proof.
\end{proof}
\begin{proposition}
  \label{prop:H-alpha-stable}
  Let $T'$ and $T''$ be $\alpha$-polystable triples.  Then the
  holomorphic chain $\widetilde{C^{\bullet}}(T'',T')$ is
  $\balpha$-polystable for
  \begin{math}
    \balpha = (\alpha_1,\alpha_2) = (\alpha,2\alpha).
  \end{math}
\end{proposition}
\begin{proof}
  Since the triples $T'$ and $T''$ are $\alpha$-polystable, it follows
  from Theorem \ref{thm:vortices-hitchin-kobayashi} that they
  support solutions to the $(\tau_{1}',\tau_{2}')$- and
  $(\tau_{1}'',\tau_{2}'')$-vortex equations, respectively, where
  $\alpha = \tau_{1}' - \tau_{2}' = \tau_{1}'' - \tau_{2}''$.  
  Notice
  that $\tau_{1}' - \tau_{1}''= \tau_{2}' - \tau_{2}''$.  Thus it 
  follows 
  from Lemma~\ref{lemma:chain-vortex-solution} that the holomorphic chain
  $\widetilde{C^{\bullet}}(T'',T')$ supports a solutions to the chain
  vortex equations for $\btau = (\tau_{1}' - \tau_{2}'',\tau_{1}' -
  \tau_{1}'',\tau_{2}' - \tau_{1}'')$.  
  Now
  Theorem~\ref{thm:chain-hitchin-kobayashi} and
  \eqref{eq:alpha-tau} imply that $\widetilde{C^{\bullet}}(T'',T')$ is
  $\balpha$-polystable for 
  \begin{align*}
    \alpha_1 &= \tau_{1}' - \tau_{2}'' - \tau_{2}' + \tau_{2}'' =
    \alpha, \\
    \alpha_2 &= \tau_{1}' - \tau_{2}'' - \tau_{2}' + \tau_{1}'' =
    2\alpha. 
  \end{align*}
\end{proof}

\subsection{Bounds  for $\chi(T'',T')$}

We start with some technical lemmas needed to estimate the Euler 
characteristic $\chi(T'',T')$.

\begin{lemma}\label{lemma:stability-for-H_*}
  Let $T'=(E'_1,E'_2,\phi')$ and $T''=(E''_1,E''_2,\phi'')$ be triples
  for which the chain $\widetilde{C^{\bullet}}(T'',T')$ is 
  $\balpha=(\alpha,2\alpha)$-polystable.  Let
  \begin{align*}
    C_{1} &= {E_{1}''}^{*} \otimes E_{1}' \oplus
    {E_{2}''}^{*} \otimes E_{2}', \\
    C_{0} &= {E_{2}''}^{*} \otimes E_{1}',
  \end{align*}
  and $a_{1} \colon C_{1} \to C_{0}$ be defined as in 
  (\ref{eq:H-chain}).
  Then the following inequalities hold.
  \begin{align}
    \label{eq:ker-a_1}
    \deg(\ker(a_{1})) &\leq \rk(\ker(a_{1}))(\mu_{\alpha}(T') -
    \mu_{\alpha}(T'')), \\
    \label{eq:im-a_1}
    \deg(\im(a_{1})) &\leq \bigl(\rk(C_{0}) - \rk(\im(a_{1}))\bigr)
    (\mu_{\alpha}(T'') - \mu_{\alpha}(T') - \alpha) + \deg(C_{0}).
  \end{align}
\end{lemma}
\begin{proof} If $\rk(\ker(a_{1}))=0$ then (\ref{eq:ker-a_1}) is obvious.
Assume therefore that $\rk(\ker(a_{1}))>0$.  Using $\ker(a_{1})$, we 
can then define a quotient of the chain 
  $\widetilde{C^{\bullet}}(T'',T')$ by 
  \begin{displaymath}
      \mathcal{K}\colon  0 \to \ker(a_{1}) \to 0.
  \end{displaymath}
  Thus, since
  $\mu_{\balpha}(\mathcal{K}) = \mu(\ker(a_{1})) + \alpha$, it
  follows  from the definition of  $\balpha$-polystability
that
  \begin{align*}
    \mu(\ker(a_{1})) + \alpha &\leq
    \mu_{\balpha}(\widetilde{C^{\bullet}}(T'',T'))
    = \mu_{\alpha}(T') - \mu_{\alpha}(T'') + \alpha.
  \end{align*}
  We therefore have
  \begin{displaymath}
    \mu(\ker(a_{1})) \leq \mu_{\alpha}(T') - \mu_{\alpha}(T''),
  \end{displaymath}
  which is equivalent to \eqref{eq:ker-a_1}.
The second inequality, i.e.\ (\ref{eq:im-a_1}), is obvious when
$\rk(\im(a_{1}))=\rk(C_0)$. We thus assume $\rk(\im(a_{1})<\rk(C_0)$.
Using the cokernel $\coker(a_{1})$ (or its saturation if it is not 
torsion free), we can define a subchain of the chain 
$\widetilde{C^{\bullet}}(T'',T')$ by 
  \begin{displaymath}
      \mathcal{Q}\colon  0 \to 0 \to \coker(a_{1}).
  \end{displaymath}
By the $\balpha$-polystability of    
$\widetilde{C^{\bullet}}(T'',T')$ we have  
$\mu_{\balpha}(\mathcal{Q})\ge 
\mu_{\balpha}(\widetilde{C^{\bullet}}(T'',T'))$. This, together with
the fact that 
$$\mu(\coker(a_1))\le\frac{\deg(C_{0})-\deg(\im(a_1))}{\rk(C_0)-\rk(\im(a_1))}\ 
,$$ 
\noindent leads directly to \ref{eq:im-a_1}. 
\end{proof}

\begin{lemma}
  \label{lem:linear-algebra}
  Let $c' \colon V'_2 \to V'_1$ and $c'' \colon V''_2 \to V''_1$ be
  linear maps between finite dimensional vector spaces.
  Assume that $V'_1 \oplus V'_2 \neq 0$ and
  $V''_1 \oplus V''_2 \neq 0$.
  Define
  \begin{align*}
    f \colon \Hom(V''_1,V'_1) \oplus \Hom(V''_2,V'_2) &\lto
    \Hom(V''_2,V'_1) \\
    (\psi_1,\psi_2) &\longmapsto c'\psi_2 - \psi_1c''.
  \end{align*}
  If $f$ is an isomorphism, then exactly one of the following
  alternatives must occur:
  \begin{itemize}
  \item[$(1)$] $V'_1 = V''_2 =0$ and $c'=c''=0$. 
  \item[$(2)$] $V''_1=0$, $V'_1, V'_2,
  V''_2 \neq 0$ and $c'\colon V'_2
  \overset{\cong}{\longrightarrow} V'_1$.  \item[$(3)$]
  $V'_2=0$, $V'_1, V''_1, V''_2
  \neq 0$ and $c''\colon V''_2
  \overset{\cong}{\longrightarrow} V''_1$.
  \end{itemize}
  In particular, if $V'_1$, $V'_2$, $V''_1$
  and $V''_2$ are all non-zero then $f$ cannot be an
  isomorphism.
\end{lemma}
\begin{proof}
  If $(c',c'') = (0,0)$ then $f = 0$ and therefore 
  \begin{displaymath}
    \Hom(V''_2,V'_1) 
    = \Hom(V''_1,V'_1) 
    = \Hom(V''_2,V'_2) = 0.
  \end{displaymath}
  If $V'_1 \neq 0$ then $V''_1 = V''_2 = 0$ i.e.
$V''_1 \oplus  V''_2 = 0$. Hence $V'_1 = 0$.  Similarly one sees that
  $V''_2 =0$ and thus alternative $(1)$ occurs.
  Henceforth assume that $(c',c'') \neq (0,0)$.  Let $r'_i = \dim
  V'_i$ and $r''_i = \dim V''_i$ for $i=1,2$.  If
  $f$ is an isomorphism then $r''_1 r'_1 + r''_2 r'_2 = r''_2 r'_1$
  from which it follows that $r''_2(r'_1 - r'_2) = r''_1 r'_1$ and
  $r'_1(r''_2-r''_1)  = r''_2 r'_2$.
  Hence
  \begin{align}
    r'_1 &\geq r'_2, \label{eq:r'_1-geq-r'_2} \\
    r''_2 &\geq r''_1. \label{eq:r''_2-geq-r''_1}
  \end{align}
  Assume that we have strict inequality in \eqref{eq:r'_1-geq-r'_2}
  and \eqref{eq:r''_2-geq-r''_1}.  Then, in particular, $\coker(c')$
  and $\ker(c'')$ must both be non-zero.  Choose a complement to
  $\im(c')$ in $V'_1$ so that 
  \begin{displaymath}
    V'_1 = \im(c') \oplus \im(c')^{\perp}.
  \end{displaymath}
  We then have an inclusion
  \begin{displaymath}
    \Hom(\ker(c''),\im(c')^{\perp}) \into \Hom(V''_2,V'_1).
  \end{displaymath}
  Let $\psi = (\psi_1,\psi_2) \in \Hom(V''_1,V'_1) \oplus
  \Hom(V''_2,V'_2)$ and $x \in \ker(c'')$, then
  \begin{displaymath}
    f(\psi)(x) = c' \psi_2 (x) - \psi_1 c'(x) = c' \psi_2 (x),
  \end{displaymath}
  which belongs to $\im(c')$.  Hence $\im(f)$ and
  $\Hom(\ker(c''),\im(c')^{\perp})$ have trivial intersection and,
  therefore, $f$ cannot be an isomorphism, which is absurd.  It
  follows that equality must hold in at least one of the inequalities
  \eqref{eq:r'_1-geq-r'_2} and \eqref{eq:r''_2-geq-r''_1}.
  Suppose that equality holds in \eqref{eq:r'_1-geq-r'_2}, i.e. $r'_1
  =r'_2 = 0$. Then $r_1'r_1''=0$,
 i.e. $r_1'=0$ or $r_1''=0$.
Suppose first that $r'_1 =0$, then $r_2'=r_1'=0$, which contradicts
our assumption that $V_1'\oplus V_2'\neq 0$. Thus we must have 
$r_1''=0$ and $r_1'\neq 0$. We thus also have 
  $V'_2 \neq 0$ (since $r_2'=r_1'$) and $V_2'\neq 0$ 
(since $r_2''+r_1''\neq 0)$. Furthermore, since $(c',c'')\neq (0,0)$
we can assume that $c''\neq 0$. In this case  
$f(\psi_1,\psi_2) = f(\psi_1,0)= -\psi_1 c''$.
In particular, if $f$ is an isomorphism then so is $c''$.
Thus  alternative $(2)$ occurs.
  In a similar manner one sees that if equality holds in
  \eqref{eq:r''_2-geq-r''_1} then alternative $(3)$ occurs.  Obviously
  the three alternatives are mutually exclusive.
\end{proof}

\begin{lemma} \label{lem:linear-algebra-triples}
Suppose that $T'$ and $T''$ are non-zero triples of types
$(n_1',n_2',d_1',d_2')$ and $(n_1'',n_2'',d_1'',d_2'')$ respectively.
Let $n_1=n_1'+n_1''$, $n_2=n_2'+n_2''$,
$d_1=d_1'+d_1''$, $d_2=d_2'+d_2''$, $\mu_1=d_1/n_1$, and $\mu_2=d_2/n_2$.
  Let $\alpha_m$ and $\alpha_M$ be the extreme $\alpha$ values for the
   triples of  type $(n_1,n_2,d_1,d_2)$, as defined in 
(\ref{alpha-bounds-m}) and 
  (\ref{alpha-bounds-bigM}), with the convention that $\alpha_M=\infty$ 
  if $n_1=n_2$.
Let $\alpha_m<\alpha<\alpha_M$ and suppose 
that $\mu_{\alpha}(T')=\mu_{\alpha}(T'')$, then the  map
  \begin{displaymath}
    a_{1}\colon {E_{1}''}^{*} \otimes E_{1}' \oplus {E_{2}''}^{*}
    \otimes E_{2}' \to {E_{2}''}^{*} \otimes E_{1}'
  \end{displaymath}
  cannot be  an isomorphism.
  \end{lemma}
\begin{proof}
Let us consider the triple $T=T'\oplus T''$. It is clear that
$\mu_\alpha(T)=\mu_\alpha(T')=\mu_\alpha(T'')$.

 If $a_1$ is an isomorphism then, applying 
 Lemma~\ref{lem:linear-algebra} fibrewise, it follows that one
  of the following alternatives must occur:
  \begin{itemize}
  \item[(a)] $E'_1 = E''_2 =0$ and $\phi'=\phi''=0$.
  \item[(b)] $E''_1=0$, $E'_1, E'_2, E''_2 \neq 0$ and $\phi'\colon E'_2
    \overset{\cong}{\longrightarrow} E'_1$.
  \item[(c)] $E'_2=0$, $E'_1, E''_1, E''_2 \neq 0$ and $\phi''\colon E''_2
    \overset{\cong}{\longrightarrow} E''_1$.
  \end{itemize}
  We shall consider each case in turn.
  \emph{Case (a).}  In this case we have $T'=(0,E_2,0)$,
  $T''=(E_1,0,0)$ and $T=(E_1,E_2,0)$. It follows from 
  $\mu_{\alpha}(T')=\mu_{\alpha}(T)$ that
  $\alpha=\mu(E_1)-\mu(E_2)=\alpha_m$.
  \emph{Case (b).}  In this case we have 
  $n_1 = n'_1$ and $n_2 = n'_2 + n''_2 = n'_1 + n''_2$.  Hence $n_2 > 
  n_1$. Furthermore, from $\mu_{\alpha}(T')=\mu_{\alpha}(T)$ we get
  $\mu(E_1)+\frac{\alpha}{2}=\mu(E_1\oplus E_2)+\frac{n_2}{n_1+n_2}$, 
  i.e.\ $\alpha=\frac{2n_2}{n_2-n_1}\alpha_m=\alpha_M$.
  \emph{Case (c).}  In this case we have 
 $n_2 = n''_2 $ and  $n_1 = n'_1+n_1'' = n'_1 + n''_2$.  Hence $n_1 > 
  n_2$. Furthermore, from $\mu_{\alpha}(T'')=\mu_{\alpha}(T)$ we get
$\alpha=\frac{2n_1}{n_1-n_2}\alpha_m=\alpha_M$.
If $n_1=n_2$ then case (a) is the only possibility, so  
$\alpha=\alpha_m$. If  $n_1\ne n_2$, then (a) or exactly one of (b) and (c) 
are the only possibilities, depending on whether $n_1<n_2$ or 
$n_1>n_2$. In both cases we see that $\alpha=\alpha_m$ or $\alpha=\alpha_M$. 
\end{proof}

\begin{proposition}\label{prop:hyper-1-vanishing-criterion}
Suppose that $T'$ and $T''$ are non-zero triples of types
$(n_1',n_2',d_1',d_2')$ and $(n_1'',n_2'',d_1'',d_2'')$ respectively.
Let $n_1=n_1'+n_1''$, $n_2=n_2'+n_2''$,
$d_1=d_1'+d_1''$, $d_2=d_2'+d_2''$, $\mu_1=d_1/n_1$, and $\mu_2=d_2/n_2$.
  Let $\alpha_m$ and $\alpha_M$ be the extreme $\alpha$ values for the
   triples of  type $(n_1,n_2,d_1,d_2)$, as defined in 
(\ref{alpha-bounds-m}) and 
  (\ref{alpha-bounds-bigM}), with the convention that $\alpha_M=\infty$ 
  if $n_1=n_2$.
Let $\alpha_m<\alpha<\alpha_M$. Suppose 
that $\mu_{\alpha}(T')=\mu_{\alpha}(T'')$ and  that 
the chain  $\widetilde{C^{\bullet}}(T'',T')$, as defined in  
(\ref{eq:H-chain}),
is $(\alpha,2\alpha)$-stable. Then 
  \begin{displaymath}
    \chi(T'',T')\le 1-g 
  \end{displaymath}
\noindent if $\alpha \geq  2g-2$. In particular, if $g\geq 2$ then
$\chi(T'',T')\le 0$. 
\end{proposition}

\begin{proof}
{}From the long exact sequence \eqref{eq:long-exact-extension-complex}
  and the Rie\-mann-Roch formula we obtain
  \begin{equation}
    \label{eq:chi-T''-T'}
    \chi(T'',T')
    = (1-g)\bigl(\rk(C_{1}) - \rk(C_{0})\bigr) + \deg(C_{1}) -
    \deg(C_{0}), 
  \end{equation}
   \noindent where $C_1$ and $C_0$ are as in (\ref{eq:H-chain}).  
   We can apply Lemma \ref{lemma:stability-for-H_*}, and then use the
   estimates (\ref{eq:ker-a_1}) and (\ref{eq:im-a_1}). Together with
  \begin{align}
    \deg(C_{1}) &= \deg(\ker(a_{1})) + \deg(\im(a_{1})), 
    \label{eq:deg(H_1)}\\
    \rk(C_{1})&=\rk(\ker(a_{1})) + \rk(\im(a_{1})) ,\label{eq:rk(H_1)}
  \end{align}
  
  \noindent these yield
    \begin{multline*}
      \deg(C_{1}) \leq (\mu_{\alpha}(T') - \mu_{\alpha}(T''))
    \bigl(\rk(C_{1}) - \rk(C_{0})\bigr) \\
    - \alpha\bigl(\rk(C_{0}) - \rk(\im(a_{1}))\bigr) + \deg(C_{0}).
    \end{multline*}
  
  \noindent Using that $\mu_{\alpha}(T') = \mu_{\alpha}(T'')$, we can then 
  deduce that
  \begin{displaymath}
      \deg(C_{1}) - \deg(C_{0})
      \leq  -\alpha\bigl(\rk(C_{0}) - \rk(\im(a_{1}))\bigr).
  \end{displaymath}
  Combining this with \eqref{eq:chi-T''-T'} we get
  \begin{equation}\label{eqn:chi-est}
    \chi(T'',T')
    \leq (1-g)\bigl(\rk(C_{1}) - \rk(C_{0})\bigr)
      - \alpha\bigl(\rk(C_{0}) - \rk(\im(a_{1}))\bigr).
  \end{equation}
  If $\alpha\ge 2g-2$ then we get 
 $$
 \chi(T'',T')\le (1-g)\bigl(\rk(C_{0})+\rk(C_{1}) - 2\rk(\im(a_{1}))\bigr),
 $$
with equality if and only if $\alpha=2g-2$. Furthermore $\rk(\im(a_{1})) 
\leq \rk(C_{0})$ and $\rk(a_{1}) \leq \rk(C_{1})$, with equality in both if and only 
if $a_1$ is an isomorphism. Thus in all cases we get 
$\chi(T'',T')\le 0$, with equality if and only if $\alpha=2g-2$ and
$a_1$ is an isomorphism. But by Lemma \ref{lem:linear-algebra-triples}, since 
$\alpha_m<\alpha<\alpha_M$, then $a_1$ cannot be an isomorphism. 
Thus in all cases we get 
$\rk(C_{0})+\rk(C_{1})-2\rk(\im(a_{1}))\ge 1$ and hence $\chi(T'',T')\le 1-g$. 
 \end{proof}
\begin{remark}
Since the roles of $T'$ and $T''$ in Proposition 
\ref{prop:hyper-1-vanishing-criterion}  are symmetric,
we obtain the same bound for $\chi(T',T'')$.
\end{remark}

\section{Crossing critical values}
\label{sec:crossing-critical-values}

In this section we study the differences between the stable loci 
$\mathcal{N}^s_\alpha(\mathbf{n},\mathbf{d})$ in the moduli spaces 
$\mathcal{N}_\alpha(\mathbf{n},\mathbf{d})$, for fixed
values of $\mathbf{n}=(n_1,n_2)$ and 
$\mathbf{d} = (d_1,d_2)$ but different values of $\alpha$. Since
in this section $\mathbf{n}$ and $\mathbf{d}$ are fixed, we use the 
abbreviated notation 
$$\mathcal{N}^s_{\alpha}=\mathcal{N}^s_\alpha(\mathbf{n},\mathbf{d})
\ \quad\ \mathrm{and}\ \quad\ 
\mathcal{N}_{\alpha}=\mathcal{N}_\alpha(\mathbf{n},\mathbf{d}).$$
\noindent Our main result is that for all $\alpha\ge 2g-2$ 
any differences between the $\mathcal{N}^s_{\alpha}$ are confined to 
subvarieties of positive codimension. In particular, the number of 
irreducible components of the spaces  $\mathcal{N}^s_{\alpha}$ are the 
same for all $\alpha$ satisfying 
$\alpha\ge 2g-2$ and $\alpha_m<\alpha <\alpha_M$\footnotemark .
\footnotetext{When $n_1\ne n_2$ the bounds $\alpha_m$ and $\alpha_M$ 
are as in (\ref{alpha-bounds-m}) and (\ref{alpha-bounds-bigM}). When 
$n_1=n_2$ we adopt the convention that $\alpha_M=\infty$}If 
the coprimality condition $\GCD(n_2,n_1+n_2,d_1+d_2)=1$ is satisfied, 
then $\mathcal{N}^s_{\alpha}=\mathcal{N}_{\alpha}$ at all 
non-critical vales of $\alpha$, so the results apply to   
$\mathcal{N}_{\alpha}$ for all non-critical $\alpha\ge 2g-2$.

We begin with a set theoretic description of the differences between 
two spaces $\mathcal{N}^s_{\alpha_1}$ and  $\mathcal{N}^s_{\alpha_2}$ 
when $\alpha_1$ and $\alpha_2$ are separated by a critical value (as 
defined in section \ref{sec:critical-values}). For the rest of 
this section we adopt the following notation: Let $\alpha_c$ be a 
critical value such that 
\begin{equation}\label{eqn:alphac-range}
\alpha_m <\alpha_c <\alpha_M.
\end{equation}
\noindent Set 
\begin{equation}\label{eqtn:alpha-c-pm}
\alpha_c^+ = \alpha_c + \epsilon,\quad 
\alpha_c^- = \alpha_c - \epsilon,
\end{equation}
\noindent where $\epsilon > 0$ is small enough 
so that $\alpha_c$ is the only critical value in the interval 
$(\alpha_c^-,\alpha_c^+)$. 

\subsection{Flip Loci}\label{subs:fliploci}

\begin{definition}Let $\alpha_c\in (\alpha_m,\alpha_M)$ 
be a critical value for triples of type $(\mathbf{n},\mathbf{d})$.  
We define {\it flip loci} 
$\mathcal{S}_{\alpha_c^{\pm}}\subset\mathcal{N}^s_{\alpha_c^{\pm}}$ 
by the conditions that the points in $\mathcal{S}_{\alpha_c^+}$ represent 
triples which are $\alpha_c^+$-stable  but $\alpha_c^-$-unstable, 
while the points in $\mathcal{S}_{\alpha_c^-}$ represent triples 
which are 
$\alpha_c^-$-stable  but $\alpha_c^+$-unstable.
\end{definition}

\begin{remark}The definition of $\mathcal{S}_{\alpha_c^+}$ can be 
extended to the extreme case $\alpha_c=\alpha_m$. However, since all 
$\alpha_m^+$-stable triples must be $\alpha_m^-$-unstable, we see that
$\mathcal{S}_{\alpha_m^+}=\mathcal{N}^s_{\alpha_m^+}$. Similarly, when 
$n_1\ne n_2$ we get $\mathcal{S}_{\alpha_M^-}=\mathcal{N}^s_{\alpha_M^-}$.
 The only interesting cases are thus those those for which 
 $\alpha_m<\alpha_c <\alpha_M$. 
\end{remark}

\begin{lemma}\label{lemma:fliploci}In the above notation:
\begin{equation}\label{eqn:Nsalpha}
\mathcal{N}^s_{\alpha_c^+}-\mathcal{S}_{\alpha_c^+}=
\mathcal{N}^s_{\alpha_c}=
\mathcal{N}^s_{\alpha_c^-}-\mathcal{S}_{\alpha_c^-}.
\end{equation}
\end{lemma}

\begin{proof} By definition we can identify 
$\mathcal{N}^s_{\alpha_c^+}-\mathcal{S}_{\alpha_c^+}=
\mathcal{N}^s_{\alpha_c^-}-\mathcal{S}_{\alpha_c^-}$.

Suppose now that $t$ is a point in 
$\mathcal{N}^s_{\alpha_c^+}-\mathcal{S}_{\alpha_c^+}=
\mathcal{N}^s_{\alpha_c^-}-\mathcal{S}_{\alpha_c^-}$, but that $t$ is not
in $\mathcal{N}^s_{\alpha_c}$. Let $T$ be a triple representing $t$. 
Then $T$ has a subtriple $T'\subseteq T$ for which 
$\mu_{\alpha_c}(T')\ge\mu_{\alpha_c}(T)$, and also 
$\mu_{\alpha_c^{\pm}}(T')<\mu_{\alpha_c^{\pm}}(T)$.
This is not possible, and hence $t\in\mathcal{N}^s_{\alpha_c}$. 

Finally, suppose that $t\in\mathcal{N}^s_{\alpha_c}$ and let $T$ be a 
triple representing $t$. Then $\mu_{\alpha_c}(T')<\mu_{\alpha_c}(T)$ 
for all subtriples $T'\subset T$. But since the set of possible 
values for 
$\mu_{\alpha_c}(T')$ is a discrete subset of $\R$, we can find a 
$\delta>0$ such that $\mu_{\alpha_c}(T')-\mu_{\alpha_c}(T)\le -\delta$ 
for all subtriples $T'\subset T$. Thus 
$\mu_{\alpha_c^{\pm}}(T')-\mu_{\alpha_c^{\pm}}(T)<0$. That is, $t$ 
is in $\mathcal{N}^s_{\alpha^{\pm}}$, and hence 
$\mathcal{N}^s_{\alpha} \subseteq
\mathcal{N}^s_{\alpha^{\pm}}-\mathcal{S}_{\alpha^{\pm}}$. 
\end{proof}

\noindent Our goal is to show that the flip loci 
$\mathcal{S}_{\alpha_c^{\pm}}$ are contained in subvarieties of 
positive codimension in $\mathcal{N}^s_{\alpha_c^{\pm}}$ 
respectively. 

\begin{proposition}\label{prop:vicente}
Let $\alpha_c\in (\alpha_m,\alpha_M)$  be a critical value for 
triples of type 
$(\mathbf{n},\mathbf{d})=(n_1,n_2,d_1,d_2)$. 
Let $T=(E_1,E_2,\phi)$ be a triple of this type. 
\begin{itemize}
\item[(1)] Suppose that $T$ represents a point in 
$\mathcal{S}_{\alpha_c^+}$, i.e.\ suppose that $T$ is
$\alpha_c^+$-stable but $\alpha_c^-$-unstable.  Then $T$ has a 
description as the middle term in an extension 
 \begin{equation}\label{destab}
 0\to T'\to T\to T'' \to 0
 \end{equation}
in which 
 \begin{itemize}
 \item[(a)]  $T'$ and $T''$ are both $\alpha_c^+$-stable, with 
 $\mu_{\alpha_c^+}(T')<\mu_{\alpha_c^+}(T)$,
 \item[(b)] $T'$ and $T''$ are both $\alpha_c$-semistable with
 $\mu_{\alpha_c}(T')=\mu_{\alpha_c}(T)$.
 \end{itemize}
 
\item[(2)] Similarly, if $T$ represents a point in 
$\mathcal{S}_{\alpha_c^-}$, i.e.\ if $T$ is $\alpha_c^-$-stable but
$\alpha_c^+$-unstable,  then $T$ has a 
description as the middle term in an extension 
$\mathrm{(\ref{destab})}$
in which 
 \begin{itemize}
 \item[(a)]  $T'$ and $T''$ are both $\alpha_c^-$-stable with 
 $\mu_{\alpha_c^-}(T')<\mu_{\alpha_c^-}(T)$,
 \item[(b)] $T'$ and $T''$ are both $\alpha_c$-semistable with
 $\mu_{\alpha_c}(T')=\mu_{\alpha_c}(T)$.
 \end{itemize}
\end{itemize}
\end{proposition}
{\em Proof.\/} In both cases (i.e.\ (1) and (2)), since its stability 
property changes at 
$\alpha_c$, the triple $T$ must be strictly $\alpha_c$-semistable, i.e.\ 
it must have a proper subtriple $T'$ with 
$\mu_{\alpha_c}(T')=\mu_{\alpha_c}(T).$  We can thus consider the (non-empty) set
 $$
 {\mathcal F}_1=\{T'\subsetneq T\ |\
 \mu_{\alpha_c}(T')= \mu_{\alpha_c}(T)\ \}.
 $$
\noindent {\em Proof of (1)}  
Suppose first that $T$ is 
$\alpha_c^+$-stable but $\alpha_c^-$-unstable. We observe that if
$T'\in {\mathcal F}_1$, then $\frac{n_2'}{n_1'+n'_2}< \frac{n_2}{n_1+n_2}$, since
otherwise $T$ could not be $\alpha_c^+$-stable. But the allowed 
values for 
$\frac{n'_2}{n_1'+n'_2}$ are limited by the constraints 
$0\le n_1'\le n_1$, $0\le n'_2\le n_2$ and $n_1'+n'_2\ne 0$. We can thus define
 $$
  \lambda_0=\max\left\{\frac{n'_2}{n_1'+n'_2}\ \biggm|\ T'\in {\mathcal F}_1\ \right\}\
 $$
and set 
 $$
 {\mathcal F}_2=\left\{T_1\subset {\mathcal F}_1\ \biggm|\\
 \frac{n'_2}{n_1'+n'_2}=\lambda_0\ \right\}.
 $$
\noindent Now let $T'$ be any triple in $\mathcal F_2$. Since $T'$ has maximal 
$\alpha_c$-slope, we can assume that $T''=T/T'$ is a locally free 
triple, i.e.\ if 
$T''=(E_2'',E_1'',\Phi)$ then $E_2''$ and $E_1''$ are both locally 
free.  Furthermore, since $T$ is 
$\alpha_c$-semistable and $\mu_{\alpha_c}(T')=
\mu_{\alpha_c}(T)=\mu_{\alpha_c}(T'')$, it follows that both
$T'$ and $T''$ are $\alpha_c$-semistable and of the same $\alpha_c$-slope.
We now show that $T''$ is $\alpha_c^+$-stable. Suppose not. Then 
there is a proper subtriple $\tilde{T}''\subset T''$ with 
$\mu_{\alpha_c^+}(\tilde{T}'')\ge\mu_{\alpha_c^+}(T'')$. 
However, since we can assume that $\alpha_c^+$ is not a critical 
value for triples of type $(\tilde{T}'')$, we must have 
$$\mu_{\alpha_c^+}(\tilde{T}'')>\mu_{\alpha_c^+}(T'').$$  
Thus, since $(T'')$ is $\alpha_c$-semistable, we must have
$\mu_{\alpha_c}(\tilde{T}'')\le\mu_{\alpha_c}(T'')$ 
and also 
$$
\frac{\tilde{n}_2''}{\tilde{n}_1''+\tilde{n}_2''}>\frac{n''_2}{n_1''+n''_2}.
$$
If $\mu_{\alpha_c}(\tilde{T}'')<\mu_{\alpha_c}(T'')$, say 
$\mu_{\alpha_c}(\tilde{T}'')=\mu_{\alpha_c}(T'')-\delta$, then in order to have
$\mu_{\alpha_c^+}(\tilde{T}'')>\mu_{\alpha_c^+}(T'')$ we must have
$$\frac{\tilde{n}_2''}{\tilde{n}_1''+\tilde{n}_2''}>
\frac{n''_2}{n_1''+n''_2}+\frac{\delta}{\epsilon}.$$
\noindent Letting $\epsilon$ approach zero, we see that 
$\frac{\tilde{n}_2''}{\tilde{n}_1''+\tilde{n}_2''}$ must be arbitrarily large. 
This cannot be if $0\le \tilde{n}_1''\le n_1''$ and 
$0\le \tilde{n}_2''\le n''_2$ (and $\tilde{n}_1''+\tilde{n}_2''>0$).  
We may thus assume that 
$\mu_{\alpha_c}(\tilde{T}'')=\mu_{\alpha_c}(T'').$ 
Consider now the subtriple $\tilde{T}'\subset T$ defined by the 
pull-back diagram 
 $$
 0\to T'\to \tilde{T}'\to \tilde{T}'' \to 0 .
 $$
\noindent This has $\mu_{\alpha_c}(\tilde{T}')=\mu_{\alpha_c}(T'')=\mu_{\alpha_c}(T)$ and 
thus 
$$\frac{\tilde{n}_2'}{\tilde{n}_1'+\tilde{n}_2'}\le\lambda_0=
\frac{n'_2}{n_1'+n'_2}.$$
\noindent It follows from this and the above extension that
$$\frac{\tilde{n}_2''}{\tilde{n}_1''+\tilde{n}_2''}\le\lambda_0=
\frac{n'_2}{n_1'+n'_2}.$$
\noindent However, since $\mu_{\alpha_c}(T')=\mu_{\alpha_c}(T)$ but 
$\mu_{\alpha_c^+}(T')<\mu_{\alpha_c^+}(T)$, we have that
$$\frac{n'_2}{n_1'+n'_2}<\frac{n''_2}{n_1''+n''_2}.$$
\noindent Combining the previous two inequalities we get
$$\frac{\tilde{n}_2''}{\tilde{n}_1''+\tilde{n}_2''}<\frac{n''_2}{n_1''+n''_2}\ $$
\noindent which is a contradiction. 

Now take $T'\in\mathcal F_2$ with minimum rank (i.e.\ minimum 
$n_1'+n'_2$) in  $\mathcal F_2$. We claim that $T'$ is 
$\alpha_c^+$-stable. If not, then as before
it has a proper subtriple $\tilde{T}'$ with 
 $\mu_{\alpha_c}(\tilde{T}')\le\mu_{\alpha_c}(T')$ and
$\frac{\tilde{n}_2'}{\tilde{n}_1'+\tilde{n}_2'}>\frac{n'_2}{n_1'+n'_2}$.
Then $\tilde{n}_1'+\tilde{n}_2'<n_1'+n'_2$, which contradicts 
the minimality of $n_1'+n'_2$. Thus $T'$ is $\alpha_c^+$-stable. 
Moreover, since $T$ is $\alpha_c^+$-stable it follows that 
$\mu_{\alpha_c^+}(T')<\mu_{\alpha_c^+}(T)$.
Thus taking $T'\in\mathcal F_2$ with minimum rank, and $T''=T/T'$, we 
get a description of $T$ as an extension in which (a)-(b) are 
satisfied. 
\bigskip

\noindent {\em Proof of (2).}
If 
$T$ is $\alpha_c^-$-stable but 
$\alpha_c^+$-unstable, then $\frac{n'_2}{n_1'+n'_2}> \frac{n_2}{n_1+n_2}$ for all
$T'\in\mathcal F_1$. The proof of (a) must thus be modified as
follows. With 
 $$
 \lambda_0=\min\left\{\frac{n'_2}{n_1'+n'_2}\ \biggm|\ T'\in {\mathcal F}_1\ \right\}\
 $$
\noindent we can define 
 $$
 {\mathcal F}_2=\left\{T'\subset {\mathcal F}_1\ \biggm|\\
 \frac{n'_2}{n_1'+n'_2}=\lambda_0\ \right\}
 $$
\noindent and select $T'\in\mathcal F_2$ such that $T'$ has minimal 
rank in 
$\mathcal F_2$. It follows in a similar fashion to that above that
$T$ has a description as
 $$
 0\to T'\to T\to T'' \to 0
 $$
in which all the requirements of the proposition are satisfied. 
\qed
\begin{remark}Unlike for Jordan-H\"older filtrations for semistable 
objects, the filtrations produced by the above proposition are always 
of length two, i.e.\ always yield a description of the semistable 
object as an extension of stable objects. This is achieved by 
exploiting the extra `degree of freedom' provided by the parameter 
$\alpha_c$. The true advantage of never having to consider extensions 
of length greater than two is that it removes the need for inductive 
procedures in the analysis of the flip loci. 
\end{remark}
\begin{definition}
\label{def:S-plusminus} 
Let $\alpha_c\in (\alpha_m,\alpha_M)$ 
be a critical value for triples of type 
$(\mathbf{n},\mathbf{d})$.  Let
$(\mathbf{n}',\mathbf{d}') = (n'_1,n'_2,d'_1,d'_2)$ and
$(\mathbf{n}'',\mathbf{d}'') = (n''_1,n''_2,d''_1,d''_2)$ be such that
\begin{equation}\label{eqn:C1}
(\mathbf{n},\mathbf{d})=(\mathbf{n}',\mathbf{d}')+ 
(\mathbf{n}'',\mathbf{d}''), 
\end{equation}
(i.e.\ $n_1 = n'_1 + n''_1$, $n_2 = n'_2 + n''_2$, $d_1 = d'_1 + 
d''_1$, and $d_2 = d'_2 + d''_2$), and also 
\begin{equation}\label{eqn:C2}
 \frac{d'_1+d'_2}{n'_1+n'_2}+\alpha_c\frac{n'_2}{n'_1+n'_2}=
\frac{d''_1+d''_2}{n''_1+n''_2}+\alpha_c\frac{n''_2}{n''_1+n''_2}.
\end{equation}
\begin{itemize}
\item[(1)] Define 
$\tilde{\mathcal{S}}_{\alpha_c^+}(\mathbf{n}'',\mathbf{d}'',\mathbf{n}',
\mathbf{d}')$ to be the set of all isomorphism classes of extensions 
\begin{displaymath}
  0 \lto T' \lto T \lto T'' \lto 0,
\end{displaymath}
where $T'$ and $T''$ are $\alpha_c^+$-stable triples with topological 
invariants $(\mathbf{n}',\mathbf{d}')$ and 
$(\mathbf{n}'',\mathbf{d}'')$ respectively, and the isomorphism is
on the triple $T$. 
\item[(2)] Define 
$\tilde{\mathcal{S}}^{0}_{\alpha_c^+}(\mathbf{n}'',\mathbf{d}''
,\mathbf{n}',\mathbf{d}') 
\subset \tilde{\mathcal{S}}_{\alpha_c^+}(\mathbf{n}'',\mathbf{d}'',
\mathbf{n}',\mathbf{d}')$
to be the set of all extensions for which moreover  $T$ is 
$\alpha_c^+$-stable.  In an analogous manner, define 
$\tilde{\mathcal{S}}_{\alpha_c^-}(\mathbf{n}'',\mathbf{d}'', 
\mathbf{n}',\mathbf{d}')$ and
$\tilde{\mathcal{S}}^{0}_{\alpha_c^-}(\mathbf{n}'',\mathbf{d}'' 
,\mathbf{n}',\mathbf{d}') 
\subset \tilde{\mathcal{S}}_{\alpha_c^+}(\mathbf{n}'',\mathbf{d}'',
\mathbf{n}',\mathbf{d}')$.
\item[(3)]  Define
\begin{align*}
\tilde{\mathcal{S}}_{\alpha_c^+}= \bigcup
\tilde{\mathcal{S}}_{\alpha_c^+}(\mathbf{n}'',\mathbf{d}'',
\mathbf{n}',\mathbf{d}')\ &,\ 
\tilde{\mathcal{S}}^0_{\alpha_c^+} = \bigcup
\tilde{\mathcal{S}}^0_{\alpha_c^+}(\mathbf{n}'',\mathbf{d}'',
\mathbf{n}',\mathbf{d}')\\
\end{align*}
\noindent where the union is over all $(n'_1,n'_2,d'_1,d'_2)$ and
$(n''_1,n''_2,d''_1,d''_2)$ such that the above conditions apply, and 
also $\frac{n'_2}{n'_1+n'_2}<\frac{n''_2}{n''_1+n''_2}$. 
\item[(4)]  Similarly, define 
\begin{align*}
\tilde{\mathcal{S}}_{\alpha_c^-}= \bigcup
\tilde{\mathcal{S}}_{\alpha_c^-}(\mathbf{n}'',\mathbf{d}'',\mathbf{n}',\mathbf{d}')\ 
&,\ 
\tilde{\mathcal{S}}^0_{\alpha_c^-}= \bigcup
\tilde{\mathcal{S}}^0_{\alpha_c^-}(\mathbf{n}'',\mathbf{d}'',\mathbf{n}',\mathbf{d}')
\end{align*}
\noindent where the union is over all $(n'_1,n'_2,d'_1,d'_2)$ and
$(n''_1,n''_2,d''_1,d''_2)$ such that the above conditions apply, and 
also $\frac{n'_2}{n'_1+n'_2}>\frac{n''_2}{n''_1+n''_2}$. 
\end{itemize}
\end{definition}
\begin{remark}It can happen that $\tilde{\mathcal{S}}^0_{\alpha_c^-}$ 
or $\tilde{\mathcal{S}}^0_{\alpha_c^-}$ is empty. For instance there 
may be no possible choices of  $(n'_1,n'_2,d'_1,d'_2)$ and 
$(n''_1,n''_2,d''_1,d''_2)$ which satisfy all the required 
conditions. In this case, the implication of the next lemma is that 
one or both of the flip loci $\mathcal{S}_{\alpha_c^{\pm}}$ is empty. 
\end{remark}
\begin{lemma}\label{lemma:vmaps}
There are maps, say 
$v^{\pm}:\tilde{\mathcal{S}}^0_{\alpha_c^{\pm}}
\longrightarrow\mathcal{N}^s_{\alpha_c^{\pm}}$, 
which map triples to their equivalence classes. The images contain 
the flip loci $\mathcal{S}_{\alpha_c^{\pm}}$. 
\end{lemma}
\begin{proof}The existence of the maps is clear. The second statement,
about the images of the maps, follows by Proposition 
\ref{prop:vicente}. Indeed, suppose that $T$ represents a point in 
$\mathcal{S}_{\alpha_c^{+}}$ and that
 $$
 0\to T'\to T\to T'' \to 0
 $$
is an extension of the type described in proposition 
\ref{prop:vicente}, with 
$T'$ a triple of type $(\mathbf{n}',\mathbf{d}')$  and $T''$
a triple of type $(\mathbf{n}'',\mathbf{d}'')$. Then  
$(\mathbf{n}',\mathbf{d}')$ and $(\mathbf{n}',\mathbf{d}')$ satisfy
conditions (\ref{eqn:C1}) and (\ref{eqn:C2}). Furthermore, since 
$\mu_{\alpha^+_c}(T')< \mu_{\alpha^+_c}(T'')$, we must have
$\frac{n'_2}{n'_1+n'_2}<\frac{n''_2}{n''_1+n''_2}$. Thus $T$ 
is contained in  $v^+(\tilde{\mathcal{S}}^0_{\alpha_c^+})$. A similar 
argument shows that $\mathcal{S}_{\alpha_c^{+}}$ is contained in  
$v^-(\tilde{\mathcal{S}}^0_{\alpha_c^-})$.
\end{proof}

\subsection{Codimension estimates and comparison of moduli spaces}\label{subs:codim}

Consider a critical value $\alpha_c\in (\alpha_m,\alpha_M)$ for 
triples of type $(\mathbf{n},\mathbf{d})$.  Fix 
$(\mathbf{n}',\mathbf{d}') = (n'_1,n'_2,d'_1,d'_2)$ and 
$(\mathbf{n}'',\mathbf{d}'') = (n''_1,n''_2,d''_1,d''_2)$ as in 
Definition \ref{def:S-plusminus}. For simplicity we shall denote the 
moduli spaces of $\alpha_c^{\pm}$-semistable triples of 
type 
$(\mathbf{n}',\mathbf{d}')$, respectively 
$(\mathbf{n}'',\mathbf{d}'')$, by 
\begin{displaymath}
  \mathcal{N}_{\alpha_c^{\pm}}' =
  \mathcal{N}_{\alpha_c^{\pm}}(\mathbf{n}',\mathbf{d}')
  \quad\text{and}\quad
  \mathcal{N}_{\alpha_c^{\pm}}'' =
  \mathcal{N}_{\alpha_c^{\pm}}( \mathbf{n}'',\mathbf{d}'').
\end{displaymath}
\begin{proposition}\label{prop:Sprojbundle} If $\alpha_c> 2g-2$ then 
$\tilde{\mathcal{S}}_{\alpha_c^\pm}(\mathbf{n}'',\mathbf{d}'',\mathbf{n}',\mathbf{d}')$
is a locally trivial fibration over $\mathcal{N}_{\alpha_c^{\pm}}' 
\times\mathcal{N}_{\alpha_c^{\pm}}''$, with projective fibers of dimension 
$$-\chi(\mathbf{n}'',\mathbf{d}'',\mathbf{n}',\mathbf{d}')-1.
$$
In particular, 
$\tilde{\mathcal{S}}_{\alpha_c^\pm}
(\mathbf{n}'',\mathbf{d}'',\mathbf{n}',\mathbf{d}')$ has dimension 
$$
1-\chi(\mathbf{n}',\mathbf{d}',\mathbf{n}',\mathbf{d}') 
 -\chi(\mathbf{n}'',\mathbf{d}'',\mathbf{n}'',\mathbf{d}'')
 -\chi(\mathbf{n}'',\mathbf{d}'',\mathbf{n}',\mathbf{d}'),
$$
\noindent where 
$\chi(\mathbf{n}',\mathbf{d}',\mathbf{n}',\mathbf{d}')$ etc. are 
as in section \ref{sec:extensions-of-triples}. The same is true for 
$\tilde{\mathcal{S}}_{\alpha_c^+}
(\mathbf{n}'',\mathbf{d}'',\mathbf{n}',\mathbf{d}')$ when $\alpha_c=2g-2$.

\end{proposition}
\begin{proof} 
{}From the defining properties of $\tilde{\mathcal{S}}_{\alpha_c^\pm} 
(\mathbf{n}'',\mathbf{d}'',\mathbf{n}',\mathbf{d}')$ there is map 
\begin{equation}\label{eqn:fibration}
\tilde{\mathcal{S}}_{\alpha_c^\pm}
(\mathbf{n}'',\mathbf{d}'',\mathbf{n}',\mathbf{d}')\longrightarrow 
 \mathcal{N}_{\alpha_c^{\pm}}' \times\mathcal{N}_{\alpha_c^{\pm}}''
\end{equation}
\noindent which sends an extension 
  $$
 0\to T'\to T\to T'' \to 0
 $$
 \noindent to the pair $([T'],[T''])$, where $[T']$ denotes the class 
 represented by $T'$ and similarly for $[T'']$. We first examine the 
 fibers of this map.  
 
 Notice that $T'$ and $T''$
 satisfy the hypothesis of Proposition \ref{prop:h2-vanishing}
 and therefore of Corollary \ref{dimension-ext1}. Notice moreover that, since 
 $\mu_{\alpha_c^{\pm}}(T')<\mu_{\alpha_c^{\pm}}(T'')$, it is not possible
 to have $T'\cong T''$. Thus  (cf.\ Corollary \ref{dimension-ext1} and 
 Proposition \ref{prop:h0-vanishing}(2)) we 
 have
 \begin{align}\label{eq:fiberdim}
 \dim \mathbb{P}(\mathrm{Ext}^1(T'',T'))&=
 \dim {\Ext}^1(T'',T')-1 \notag\\
 &=-\chi(T'',T')-1 \notag\\
 &= -\chi(\mathbf{n}'',\mathbf{d}'',\mathbf{n}',\mathbf{d}')-1,
 \end{align}
 \noindent which is independent of $T'$ and $T''$. Note that
if $\alpha_c=2g-2$, $T'$ and $T''$ satisfy the hypothesis of
Proposition \ref{prop:h0-vanishing}(2)) for $\alpha_c^+$, but not for
$\alpha_c^-$.

It remains to establish that the fibration (\ref{eqn:fibration}) is 
locally trivial.  If the coprimality conditions 
$\GCD(n_1',n_2',d_1'+d_2')=1=\GCD(n_1'',n_2'',d_1''+d_2'')$ hold then the moduli
spaces $ \mathcal{N}_{\alpha_c^{\pm}}'$ and 
$\mathcal{N}_{\alpha_c^{\pm}}''$ are fine moduli spaces (cf.\ 
\cite{Schmitt}). That is, there are universal objects, say $\mathcal{U}'$ and 
$\mathcal{U}''$, defined over
$\mathcal{N}_{\alpha_c^{\pm}}'\times X$ and 
$\mathcal{N}_{\alpha_c^{\pm}}''\times X$. These can be viewed as
coherent sheaves of algebras (cf.\ \cite{AG2}), or more precisely as 
examples of the $Q$-bundles considered in \cite{gothen-king:2002}. 
Pulling these back to 
$\mathcal{N}_{\alpha_c^{\pm}}'\times \mathcal{N}_{\alpha_c^{\pm}}''\times X$
we can construct $Hom(\mathcal{U}'',\mathcal{U}')$ (where we have 
abused notation for the sake of clarity). Taking the projection from 
$\mathcal{N}_{\alpha_c^{\pm}}'\times 
\mathcal{N}_{\alpha_c^{\pm}}''\times X$ onto
$\mathcal{N}_{\alpha_c^{\pm}}'\times \mathcal{N}_{\alpha_c^{\pm}}''$, 
we can then construct the first direct image sheaf. By the results in 
\cite{gothen-king:2002}, we can identify the fibers as hypercohomology groups 
which, in this case, parameterize extensions of triples. We thus 
obtain 
$\tilde{\mathcal{S}}_{\alpha_c^\pm}$ as the projectivization of the 
first direct image of 
$Hom(\mathcal{U}'',\mathcal{U}')$. 
If the coprimality conditions fail, then the universal objects do not 
exist globally. However they still exist locally over (analytic) open 
sets in the stable locus in the base 
$\mathcal{N}_{\alpha_c^{\pm}}'\times 
\mathcal{N}_{\alpha_c^{\pm}}''$. This is sufficient for our purpose since
by construction the image of the map in \eqref{eqn:fibration} lies in 
the stable locus. The result now follows from (\ref{eq:fiberdim}) and 
formula (\ref{eq:dim-triples}) (in Theorem 
\ref{thm:smoothdim}) as applied to 
 $\mathcal{N}_{\alpha_c^{\pm}}'$ and 
 $\mathcal{N}_{\alpha_c^{\pm}}''$.
\end{proof}

\begin{proposition}\label{prop:codim-est} 
If $\alpha_c> 2g-2$ then the loci 
$\mathcal{S}_{\alpha_c^{\pm}}\subset\mathcal{N}^s_{\alpha_c^{\pm}}$ are 
locally contained in subvarieties of  codimension bounded below by 
$$
\mathrm{min}\{-\chi(\mathbf{n}'',\mathbf{d}'',\mathbf{n}',\mathbf{d}')\},
$$ 
\noindent where the minimum is over all $(\mathbf{n}',\mathbf{d}')$ and 
$(\mathbf{n}'',\mathbf{d}'')$ which satisfy
(\ref{eqn:C1}) and (\ref{eqn:C2}) and also 
$\frac{n'_2}{n'_1+n'_2}<\frac{n''_2}{n''_1+n''_2}$ (in the case of
$\mathcal{S}_{\alpha_c^{+}}$) or
$\frac{n'_2}{n'_1+n'_2}>\frac{n''_2}{n''_1+n''_2}$
(in the case of 
$\mathcal{S}_{\alpha_c^{-}}$). The same is true for 
$\mathcal{S}_{\alpha_c^{+}}$ when $\alpha_c= 2g-2$.
\end{proposition}

\begin{proof} If $\alpha_c > 2g-2$ then we can assume $\alpha_c^{\pm}\ge 2g-2$. 
Clearly also, $\alpha_c^{+}\ge 2g-2$ when $\alpha_c=2g-2$.  Thus by 
Theorem \ref{thm:smoothdim} the moduli spaces 
$\mathcal{N}^s_{\alpha^{\pm}}$ are smooth and have dimension 
$1-\chi(\mathbf{n},\mathbf{d},\mathbf{n},\mathbf{d})$.  By Corollary \ref{cor:chi-relation} and 
Proposition \ref{prop:Sprojbundle} we obtain
\begin{align*}
\dim\mathcal{N}^s_{\alpha^{\pm}}=1 
&- \chi(\mathbf{n}',\mathbf{d}',\mathbf{n}',\mathbf{d}')
- \chi(\mathbf{n}'',\mathbf{d}'',\mathbf{n}'',\mathbf{d}'')\\ 
& - \chi(\mathbf{n}'',\mathbf{d}'',\mathbf{n}',\mathbf{d}') 
- \chi(\mathbf{n}',\mathbf{d}',\mathbf{n}'',\mathbf{d}'') \\ =
&\dim \tilde{\mathcal{S}}_{\alpha_c^\pm}(\mathbf{n}'',\mathbf{d}'',
\mathbf{n}',\mathbf{d}')- \chi(\mathbf{n}',\mathbf{d}',\mathbf{n}'',\mathbf{d}'').
\end{align*}
\end{proof}

\begin{proposition}\label{prop:chain-poly-2} 
  Let $T'$ and $T''$ be $\alpha_c^\pm$-polystable triples. Then the
  holomorphic chain $\widetilde{C^{\bullet}}(T'',T')$ (as defined in
  (\ref{eq:H-chain})) is $(\alpha_c,2\alpha_c)$-polystable.
\end{proposition}
\begin{proof} 
  {}From Proposition \ref{prop:H-alpha-stable}, we have that the
  $\alpha_c^\pm$-polystability of $T'$ and $T''$ implies the
  $(\alpha_c^\pm,2\alpha_c^\pm)$-polystability of
  $\widetilde{C^{\bullet}}(T'',T')$.
  
  Now, the critical values for the chain form a discrete set of points
  in the $(\alpha_1,\alpha_2)$ plane. We can thus pick $\epsilon>0$ so
  that, with $\alpha_c^{\pm}=\alpha_c\pm \epsilon$, the point
  $(\alpha_c^{\pm},2\alpha_c^{\pm})$ is not a critical point. We can
  in fact assume that there are no critical points in
  $B_{\epsilon}^{\circ}(\alpha_c,2\alpha_c)$, i.e.\ in the punctured
  ball of radius $\epsilon$ centered at $(\alpha_c,2\alpha_c)$. Thus
  $(\alpha_c^{\pm},2\alpha_c^{\pm})$-polystability is equivalent to
  $(\alpha_c,2\alpha_c)$-polystability.
\end{proof}

 \begin{theorem}\label{thm:codim}
Let $\alpha_c\in (\alpha_m,\alpha_M)$  be a critical value for 
triples of type $(\mathbf{n},\mathbf{d})$.  If $\alpha_c> 2g-2$ then 
the loci $\mathcal{S}_{\alpha_c^{\pm}} 
\subset\mathcal{N}^s_{\alpha_c^{\pm}}$ are contained in subvarieties of 
codimension at least $g-1$. In particular, they are contained in 
subvarieties of strictly positive codimension if $g\ge 2$. If 
$\alpha_c = 2g-2$ then the same is true for 
$\mathcal{S}_{\alpha_c^+}$.
 \end{theorem}
\begin{proof}
 Combining Propositions \ref{prop:hyper-1-vanishing-criterion} and 
\ref{prop:chain-poly-2} we have that
$$
- \chi(\mathbf{n}',\mathbf{d}',\mathbf{n}'',\mathbf{d}'')=
-\chi(T',T'')\geq g-1
$$ 
(notice that the order of $T'$ and $T''$ in
these Propositions is irrelevant). The result follows now from
Proposition \ref{prop:codim-est}.
\end{proof}

\begin{theorem}\label{thm:birationality}Let $\alpha_1$ and 
$\alpha_2$ be any two values in $(\alpha_m,\alpha_M)$ 
such that $2g-2 \le\alpha_1$ and 
$\alpha_m <\alpha_1 <\alpha_2 <\alpha_M$. 
Then the  moduli spaces $\mathcal{N}^s_{\alpha_1}$ and 
$\mathcal{N}^s_{\alpha_2}$ have the same number of irreducible 
components, in particular, $\mathcal{N}^s_{\alpha_1}$ is 
irreducible if and only if $\mathcal{N}^s_{\alpha_2}$ is.
\end{theorem}

\begin{proof} This follows immediately from Theorem \ref{thm:codim} 
if $\alpha_1$ and $\alpha_2$ are non-critical, and from Theorem 
\ref{thm:codim} together with Lemma \ref{lemma:fliploci}
if either of them is critical.
\end{proof}

\section{Special values of $\alpha$}\label{sect:special-alpha}

In 
this section we identify some critical values and special
subintervals  in the range 
$(\alpha_m,\alpha_M)$. We describe their significance for the structure of 
$\alpha$-stable triples.

\subsection{Small $\alpha$}

Let $\alpha_m^+=\alpha_m+\epsilon$, with
$\epsilon$ such that the interval $(\alpha_m, \alpha_m^+]$ does not
contain any critical value (sometimes we refer to this value of 
$\alpha$ as {\em small}. The following is important in the construction
of the moduli space for small $\alpha$ (\cite{bradlow-garcia-prada:1996}).
\begin{proposition}
\label{moduli-small}If a triple $T=(E_1,E_2,\phi)$ is  $\alpha_m^+$-semistable
triple, $E_1$ and $E_2$ are semistable. In the converse direction, if one of $E_1$ and
$E_2$ is stable and the other is semistable, $T=(E_1,E_2,\phi)$ is  
$\alpha_m^+$-stable.
\end{proposition}
\begin{corollary}
If $\GCD(n_1,d_1)=1$ and $\GCD(n_2,d_2)=1$,
then $\mathcal{N}_{\alpha_m^+}^s(n_1,n_2,d_1,d_2)$
is isomorphic to the projectivization of a Picard sheaf
over $M(n_1,d_1) \times M(n_2,d_2)$, where $M(n,d)$
is the moduli space of stable bundles of rank $n$ and 
degree $d$.
\end{corollary}
\begin{proof}
Let $\mathbb{E}_1$ and $\mathbb{E}_2$ the universal bundles over
$X\times M(n_1,d_1)$ and $X \times M(n_2,d_2)$, respectively.
Consider the canonical projections 
\begin{align*}
  \pi&\colon X \times M(n_1,d_1)\times M(n_2,d_2)\to 
  M(n_1,d_1)\times M(n_2,d_2)\ ; \\
  \hat{\pi}&\colon X \times M(n_1,d_1)\times M(n_2,d_2)\to X\ ; \\
  \pi_1&\colon X \times M(n_1,d_1)\times M(n_2,d_2)\to 
  X\times M(n_1,d_1)\ ; \\
\pi_2&\colon X \times M(n_1,d_1)\times M(n_2,d_2)\to 
  X\times M(n_2,d_2)\ . 
\end{align*}
{}From Proposition \ref{moduli-small} we deduce that
$$
\mathcal{N}_{\alpha_m^+}^s(n_1,n_2,d_1,d_2)=
\mathbb{P}
(R^1\pi_*(\pi_1^*\mathbb{E}_1\otimes \pi_2^*\mathbb{E}_2^*
\otimes\hat{\pi}^*K)^*).
$$
\end{proof}

\subsection{Critical values determined by the kernel}

Throughout this section we  assume that the triple 
$(E_1,E_2,\phi)$ has type $(n_1,n_2,d_1,d_2)$, with $n_1\geq n_2$. 
The case $n_1<n_2$ can be dealt with  via duality of  triples. 

\begin{definition}
  For each integer $0\leq j < n_2$\ set
  \begin{equation}
    \label{eq:1.30}
    \alpha_j=\frac{2n_1n_2}{n_2(n_1-n_2)+
      (j+1)(n_1+n_2)}(\mu_1-\mu_2).
  \end{equation}
\end{definition}
\begin{proposition}
  \label{prop:kernel-rank}
  Let $T=(E_1,E_2,\phi)$\ be a triple in which $n_1\geq n_2$. Let
  $N\subset E_2$\ be the kernel of $\phi: E_2\longrightarrow
  E_1$. Suppose that $T$\ is $\alpha$-semistable for some
  $\alpha>\alpha_j$.  Then $N$\ has rank at most $j$. In
  particular, if $T$\ is $\alpha$-semistable for some
  $\alpha>\alpha_0$\ then $N = 0$, i.e.\ $\phi$\ is
  injective.  
\end{proposition}
\begin{proof} Suppose that
\begin{equation}
\rk(N)=k>0.\label{eq:1.40}
\end{equation}
We consider the subtriples $T_N=(0,N,0)$\ and
$T_I=(I,E_2,\phi)$, where $I$\ denotes the
image sheaf $\im(\phi)$. If $N \neq 0$, then the triple
$T_N$\ is a proper subtriple, and so is $T_I$\ since $n_1\geq n_2$.
The $\alpha$-semistability condition applied to $T_N$\ yields
\begin{equation}
\alpha n_1\leq(n_1+n_2)(\mu -\mu_N),\label{eq:1.41}
\end{equation}
where $\mu_N$\ denotes the slope of $N$\ and $\mu$\ is the slope of 
$E_1\oplus E_2$.

The $\alpha$-semistability condition applied to $T_I$\ yields
\begin{equation}
\mu(E_2\oplus I)+\alpha\frac{n_2}{i+n_2}\leq
\mu+\alpha\frac{n_2}{n_1+n_2},\label{eq:1.42}
\end{equation}
where $i=\rk(I)$.  Furthermore, from the exact sequence
\begin{equation}
0 \longrightarrow N \longrightarrow E_2 \to I \longrightarrow 0,
\label{eq:1.43}
\end{equation}
we get
\begin{align}
k+i&=n_2,\label{eq:1.44} \\
k\mu_N+i\mu_I&=n_2\mu_2.\label{eq:1.45}
\end{align}
Using \eqref{eq:1.45} we can write
\begin{equation}
\mu(E_2\oplus I)=
\frac{2n_2\mu_2-k\mu_N}{2n_2-k} \label{eq:1.46}
\end{equation}
and hence \eqref{eq:1.42} yields
\begin{equation}
-k(n_1+n_2)\mu_N\leq(n_1+n_2)((2n_2-k)\mu-2n_2\mu_2)
+\alpha n_2(n_2-k-n_1).\label{eq:1.47}
\end{equation}
Combining $k$ times \eqref{eq:1.41} and \eqref{eq:1.47} yields
\begin{equation}
\alpha \leq\frac{2n_1n_2}{n_2(n_1-n_2)+k(n_1+n_2)}(\mu_1-\mu_2). 
\label{eq:1.48}
\end{equation}
We have thus shown that if $\rk(N)=k$\ and the triple is
$\alpha$-stable, then
$$
\alpha\leq\alpha_{k-1}
$$
where $\alpha_{k-1}$\ is given by \eqref{eq:1.30} with $j=k-1$. 
Since 
$$
\alpha_{k-1}>\alpha_{k}>\dots > \alpha_{n_2-1},
$$
we can  conclude that if the triple is $\alpha$-semistable with
$\alpha>\alpha_{k-1}$, then the rank of $N$\ is strictly less than $k$.
In particular, if $\alpha>\alpha_0$, where
\begin{equation}\label{alpha0}
\alpha_0 =\frac{2n_1n_2}{n_2(n_1-n_2)+
      (n_1+n_2)}(\mu_1-\mu_2), 
\end{equation}
then $T$ is injective.
\end{proof}

As an immediate consequence we obtain the following.
\begin{proposition}\label{injective-triples}
Let $\alpha>\alpha_0$, where $\alpha_0$ is given by (\ref{alpha0}).
\begin{itemize}
\item[(1)]
An $\alpha$-semistable  triple $(E_1,E_2,\phi)$  defines a sequence  of
the form
\begin{equation}
0 \lto E_2 \overset{\phi}{\lto} E_1 \lto F\oplus S \lto 0,
\label{eq:3.20}
\end{equation}
where $F$\ is locally free and $S$\ is a torsion sheaf.
\item[(2)] If $n_1=n_2$ then
an $\alpha$-semistable  triple $(E_1,E_2,\phi)$  defines a sequence  of
the form
\begin{equation}
0 \lto E_2 \overset{\phi}{\lto} E_1 \lto S \lto 0,
\end{equation}
 where $S$ is a torsion sheaf of degree  $d_1 -d_2$.
\end{itemize}
\end{proposition}

\begin{lemma}\label{lemma:alpha-m-form} Let $\alpha_0$ be given by (\ref{alpha0}).
\begin{itemize}
\item[(1)] If $n_1>n_2$ then
\begin{equation}
\alpha_0=\frac{n_2(n_1-n_2)}{n_2(n_1-n_2)+n_1+n_2}\alpha_M
=\frac{2n_1n_2}{n_2(n_1-n_2)+n_1+n_2}\alpha_m ,
\end{equation}
where $\alpha_m$ and $\alpha_M$ are given by (\ref{alpha-bounds-m})
and (\ref{alpha-bounds-bigM}), respectively.
 
\item[(2)] 
If $n_1=n_2=n$ then 
\begin{equation}
\label{alpha0-equal-ranks}
\alpha_0=n\alpha_m=n(\mu_1-\mu_2)=d_1-d_2.
\end{equation}
\item[(3)] If $n_1\ge n_2$ then $\alpha_0\ge\alpha_m$,
with equality if and only if $\alpha_m=0$ or $n_2=1$. 
\end{itemize}
\end{lemma}
\begin{proof} Parts (1) and (2) are immediate. Using (1)
we compute 
$$\alpha_0-\alpha_m 
=\frac{n_1+n_2}{n_2(n_1-n_2)+n_1+n_2}(n_2-1)\alpha_m,$$ from which 
(3) follows. 
\end{proof}

\subsection{Critical values determined by the cokernel}

In this section we  assume that $n_1>n_2$. The range for $\alpha$
is then $[\alpha_m,\alpha_M]$,  where $\alpha_m$ and $\alpha_M$ are given 
by (\ref{alpha-bounds-m})
and (\ref{alpha-bounds-bigM}). Let us define 
\begin{equation}
\alpha_t:=\alpha_M-\frac{n_1+n_2}{n_2(n_1-n_2)}.\label{alpha-torsion}
\end{equation}
\begin{proposition}
  \label{prop:torsion-degree-bound}
  Suppose that a triple $T=(E_1,E_2,\phi)$  of the 
form (\ref{eq:3.20}) with $n_1>n_2$ 
  is $\alpha$-semistable for some $\alpha>\alpha_m$.  Then
  $$
  s\leq\frac{n_2(n_1-n_2)}{(n_1+n_2)}(\alpha_M-\alpha),
  $$
  where $s$\ is the degree of $S$.  In particular, if
  $ \alpha >\alpha_t$, then $S=0$, i.e. the quotient sheaf
$E_1/E_2$ is locally free.
\end{proposition}
\begin{proof} If $T=(E_1,E_2,\phi)$\ is of the form
  in \eqref{eq:3.20}, with $S \neq 0$, then we can find a proper
  subtriple $T'=(E'_1,E_2,\phi)$\ of the form
\begin{equation}
0 \lto E_2 \overset{\phi}{\lto} E'_1 \lto S   \lto 0.
\label{eq:3.30}
\end{equation}
Indeed, $E'_1$\ is the kernel of the sheaf map
$E_1\longrightarrow F\oplus S \longrightarrow S$. Notice that
$n'_1=n_2$\ and $d'_1=d_2+s$, where $n'_1,\ d'_1$\ denote the rank
and degree of $E'_1$,  etc. We compute
\begin{equation}
\Delta_{\alpha}(T')=
\frac{n_1}{n_1+n_2}(\mu_2-\mu_1)+
\frac{\alpha}{2}\left(\frac{n_1-n_2}{n_1+n_2}\right)+
\frac{s}{2n_2},\label{eq:3.31}
\end{equation}
where,  as in (\ref{delta-stability}), 
$\Delta_\alpha(T')=\mu_\alpha(T')-\mu_\alpha(T)$.
But
$$
\frac{n_1}{n_1+n_2}(\mu_1-\mu_2)=
\frac{\alpha_M}{2}\left(\frac{n_1-n_2}{n_1+n_2}\right)
$$
and hence
\begin{equation}
\Delta_{\alpha}(T')=
\frac{n_1-n_2}{2(n_1+n_2)}\left( \alpha-\alpha_M +
\frac{n_1+n_2}{n_2(n_1-n_2)}s \right). \label{eq:3.32}
\end{equation}
If the triple is $\alpha$-semistable then
$\Delta_{\alpha}(T')\leq 0$\ and the
result follows.
\end{proof}

Let us define
 \begin{equation}\label{alpha_e}
  \alpha_e=
\max \{\alpha_m, \alpha_0, \alpha_t\}.
 \end{equation}

The following is an  immediate consequence of Proposition 
\ref{prop:torsion-degree-bound}.

\begin{proposition} \label{prop:alpha_e}
Let $\alpha>\alpha_e$. An $\alpha$-semistable triple $(E_1,E_2,\phi)$
defines an extension
\begin{equation}\label{free-extension}
0 \lto E_2 \overset{\phi}{\lto} E_1 \lto F \lto 0,
\end{equation}
with $F$ locally free.
\end{proposition}
It turns out that for  extension like (\ref{free-extension}),
arising from semistable triples the dimension of $H^1(E_2\otimes F^*)$
does not depend on the given triple. More precisely:
\begin{proposition} \label{vanishing-morphisms}
Let  $(E_1,E_2,\phi)$ be an $\alpha$-semistable triple
in which $\ker \phi=0$ and $\coker \phi$ is locally free,
 defining an extension like (\ref{free-extension}). Then 
$H^0(E_2\otimes F^*)=0$ and hence 
\begin{equation}\label{dimension-extensions}
\dim H^1(E_2\otimes F^*)=n_2d_1-n_1d_2+n_1(n_1-n_2)(g-1).
\end{equation}
\end{proposition}
\begin{proof}
{}From \cite[Lemma 4.5]{bradlow-garcia-prada:1996} we 
have that the 
$\alpha$-semistability of $(E_1,E_2,\phi)$  for arbitrary $\alpha$  
implies that  $H^0(E_1\otimes E_2^*)=0$.
{}From  (\ref{free-extension}), we have an injective homomorphism 
$F^*\to E_1^*$, which after  tensoring with $E_2$
gives that   $H^0(E_2\otimes F^*)$ injects in 
$H^0(E_1\otimes E_2^*)$, and hence the desired vanishing.
By Riemann--Roch we obtain (\ref{dimension-extensions}).
\end{proof}

\section{Moduli space of triples  with  $n_1\neq n_2$}
\label{sec:n_1-neq-n_2}

Throughout this section we assume that  $n_1>n_2$.
The case $n_1<n_2$ can be dealt with by triples duality. 
Recall that the allowed range for the stability parameter is
$\alpha_m\leq\alpha\leq\alpha_M$,
where $\alpha_m=\mu_1-\mu_2$ and 
$\alpha_M=\frac{2n_1}{n_1-n_2}\alpha_m$, and we assume  that 
$\mu_1-\mu_2>0$.   We describe the moduli space $\mathcal{N}_\alpha$
for $2g-2 \leq \alpha \leq \alpha_M$, beginning with $\alpha=\alpha_M$.  

\subsection{Moduli space  for   $\alpha=\alpha_M$}

\begin{proposition}\label{moduli-alpha-M}
Let $T=(E_1,E_2,\phi)$ be an $\alpha_M$-polystable triple. Then  
$E_1=\im\phi\oplus F$, where $F=\coker \phi$, and $T$ decomposes as 
the direct sum of two 
$\alpha_M$-polystable triples of the same $\alpha_M$-slope, $T'$ and 
$T''$, where $T'=(\im\phi,E_2,\phi)$, and  $T''=(F,0,0)$.
In particular, $T$ is never $\alpha_M$-stable. Moreover,
$E_2\cong\im \phi$ and $E_2$  and  $F$ are polystable.
\end{proposition}
\begin{proof}
 By Proposition, \ref{prop:alpha_e}, $T$ defines an extension
\begin{equation}
0 \lto E_2 \overset{\phi}{\lto} E_1 \lto F \lto 0,
\end{equation}
with $F$ locally free. Let $T'=(\im\phi,E_2,\phi)$. Of course
$\phi:E_2\to\im\phi$ is an isomorphism, and
$$
\mu_{\alpha_M}(T')=\mu(E_2)+\frac{\alpha_M}{2},
$$
but this is equal to $\mu_{\alpha_M}(T)$ and hence
$T$ cannot be $\alpha_M$-stable.
Since we assume that $T$ is $\alpha_M$-polystable, it  must decompose
as $T'\oplus T''$, where $T''=(F,0,0)$. It is clear from
the polystability of $T$  that
$T'$ and $T''$ are $\alpha_M$-polystable with the same 
$\alpha_M$-slope. Applying the $\alpha_M$-semistability 
condition to the
subtriples $(E_2',\phi(E_2'),\phi)\subset T'$ and 
$(F',0,0)\subset T''$,  we obtain that $\mu(E_2')\leq \mu(E_2)$
and $\mu(F')\leq \mu(F)$, and hence $E_2$ and $F$ are semistable.
In fact the polystability of $T'$ and $T''$ imply
the polystability of $E_2$ and $F$, respectively.
\end{proof}

As a consequence of Proposition \ref{moduli-alpha-M}, we obtain the
following.

\begin{corollary}
\label{cor:N-alpha-M}
Suppose that $n_1>n_2$ and $\mu_1-\mu_2>0$. Then 
\begin{equation}
  \label{eqn:alpha-Mmoduli}
\mathcal{N}_{\alpha_M}(n_1,n_2,d_1,d_2)\cong
\mathcal{N}_{\alpha_M}(n_2,n_2,d_2,d_2) \times M(n_1-n_2,d_1-d_2),
\end{equation}
where $M(n_1-n_2,d_1-d_2)$ is the moduli space of polystable bundles of rank 
$n_1-n_2$ and degree $d_1-d_2$. 
\end{corollary}

\subsection{Moduli space for large $\alpha$}

Let $\alpha_L$ be the largest critical value in 
 $(\alpha_m,\alpha_M)$, and let $\mathcal{N}_L$
(respectively $\mathcal{N}^s_L$) denote the moduli space of 
$\alpha$-polystable (respectively 
$\alpha$-stable) triples for
$\alpha_L<\alpha<\alpha_M$. We refer to $\mathcal{N}_L$ as the `large
$\alpha$' moduli space. 
By definition, $\alpha_L$ is at least as big as $\alpha_e$ (where
$\alpha_e$ is as in (\ref{alpha_e})). Thus if $T$ is $\alpha$-stable
for $\alpha>\alpha_L$, then we can assume it is of the form 
(\ref{free-extension}), i.e. it gives rise to an extension 
$$
0 \lto E_2 \overset{\phi}{\lto} E_1 \lto F \lto 0.
$$
In particular,  $I=\im\phi$\ is a subbundle with 
torsion free 
  quotient in $E_1$, and  $\phi:E_2\lto I$ is an
  isomorphism. Thus we get a subtriple $T_I=(I,E_2,\phi)$ in
  which the bundles have the same rank and degree, and $\phi$ is an
  isomorphism.

\begin{proposition}  \label{prop:E2-stability}
Let $T=(E_1,E_2,\phi)$ represent a point in $\mathcal{N}_L$,
i.e.\ suppose that the triple is $\alpha$-semistable for some
$\alpha$\ in the range $\alpha_L<\alpha<\alpha_M$. Then
\begin{itemize}
\item[$(1)$] the triple $T_I=(I,E_2,\phi)$\ is $\alpha_M$-semistable,
\item[$(2)$] the bundle $E_2$\ is semistable.
\end{itemize}
\end{proposition}
\begin{proof}
(1). Let $T'=(E_1',E_2',\phi')$ be any subtriple of
$T_I$.  Since $T'$ is also a subtriple of $T$, we get
\begin{equation}
\mu_{\alpha}(T')\le
\mu_{\alpha}(T).\label{eq:3.80}
\end{equation}

\noindent A direct computation shows that
\begin{equation}
\mu_{\alpha}(T)=\mu_{\alpha}(T_I)+
\frac{n_1-n_2}{2(n_1+n_2)}(\alpha_M-\alpha), \label{eq:3.81}
\end{equation}

\noindent where we have used the fact that $n_1>n_2$ in $T$
and hence
$\alpha_M=\frac{2n_1}{n_1-n_2}(\mu_1-\mu_2)=2(\mu(F)-\mu_2)$.
Thus for all $\alpha_L<\alpha<\alpha_M$ we have

$$\mu_{\alpha}(T')-\mu_{\alpha}(T_I)\le
\frac{n_1-n_2}{2(n_1+n_2)}(\alpha_M-\alpha).$$

\noindent Taking the limit $\alpha\rightarrow\alpha_M$, we get

$$\mu_{\alpha_M}(T')-\mu_{\alpha_M}(T_I)\le 0,$$

\noindent i.e.\ $T_I$ is $\alpha_M$-semistable.

(2). Let $E'_2\subset E_2$\ be any proper subsheaf.  Then
$(\phi(E'_2),E'_2,\phi)$\ is a subtriple of
$T_I$. Since $\phi:E_2\lto\phi(E_2)$ is
an isomorphism, this subtriple has $\mu(\phi(E'_2))=\mu(E'_2)$\ and
$n'_2=n'_1$. The $\alpha_M$-semistability condition of
$T_I$  thus gives
$$\mu(E'_2)+\frac{\alpha_M}{2}\leq\mu_2+\frac{\alpha_M}{2},$$ \noindent
(where we have made use of the fact that
$\mu(\phi(E_2))=\mu(E_2)=\mu_2$).  It follows from this that
$\mu(E'_2)\leq\mu_2$, i.e.\ that $E_2$\ is semistable.
\end{proof}


\begin{proposition}
  \label{prop:triple-stable-implies-quotient-semistable}
  Suppose that the triple $T=(E_1,E_2,\phi)$\ is of the form in 
(\ref{free-extension}), i.e. 
  $$
   0 \lto E_2 \overset{\phi}{\lto} E_1 \lto F \lto 0,
  $$
with $F$\ locally free. Then there is an $\epsilon> 0$ such that $F$
is semistable if the triple is $\alpha$-semistable for any 
$\alpha >\alpha_M-\epsilon$. Indeed the
conclusion holds for any
\begin{equation}
0<\epsilon < \frac{2}{m(m-1)^2},\label{eq:3.60}
\end{equation}
where $m=n_1-n_2=\rk(F)$.
\end{proposition}
\begin{proof} Let $F'\subset F$\ be any proper
  subsheaf. 
Denote the rank and slope of $F$\ (resp.
  $F'$) by $m$\ and $\mu_F$\ (resp.\ $m'$\ and $\mu_{F'}$). We
  can always find $E'_1\subset E_1$\ such that
  $F'=E'_1/E_2$, i.e.\ such that we have
$$
0 \lto E_2 \overset{\phi}{\lto} E'_1 \lto F' \lto 0.
$$
Let $T'=(E_1',E_2,\phi)$. For convenience, define
\begin{equation}
\Delta_{\alpha}\equiv\Delta_{\alpha}(T')=
\mu_{\alpha}(T')-\mu_{\alpha}(T).
\label{eq:3.61} 
\end{equation}
Using
\begin{align}
n_1&=n_2+m, \notag \\
n'_1&=n_2+m', \notag \\
n_1\mu_1 &= n_2\mu_2+m\mu_F,\label{eq:3.62}\\
n'_1\mu'_1 &= n_2\mu_2+m'\mu_{F'}, \notag
\end{align}
we get
\begin{equation}
\mu_{F'}-\mu_F=\frac{(2n_2+m)(2n_2+m')}{2n_2m'}\Delta_{\alpha}-
\left(\frac{m-m'}{2m'}\right)(\alpha -2(\mu_F-\mu_2)).\label{eq:3.63}
\end{equation}
But $2(\mu_F-\mu_2)=\alpha_M$. Thus, setting
\begin{equation}
\alpha=\alpha_M-\epsilon,\label{eq:3.64}
\end{equation}
we get
\begin{equation}
\mu_{F'}-\mu_F=\frac{(2n_2+m)(2n_2+m')}{2n_2m'}\Delta_{\alpha}+
\left(\frac{m-m'}{m'}\right)\frac{\epsilon}{2}.\label{eq:3.65}
\end{equation}
If now we take
$$
\frac{\epsilon}{2}<\frac{1}{m(m-1)^2},
$$
then for all $0<m'<m$ we get
\begin{equation}
\left(\frac{m-m'}{m'}\right)\frac{\epsilon}{2}<
\frac{1}{m(m-1)}.\label{eq:3.66}
\end{equation}
Hence, if the triple is $\alpha$-semistable, so that 
$\Delta_{\alpha}\le 0$, then we get 
\begin{equation}
\mu_{F'}-\mu_F<\frac{1}{m(m-1)}.\label{eq:3.67}
\end{equation}
Since $\mu_F$ and $\mu_{F'}$ are rational numbers, the first with 
denominator $m$, and the second with denominator $m'\le (m-1)$, 
equation (\ref{eq:3.67}) equivalent to the condition 
$\mu_{F'} - \mu_F\leq 0$. 
\end{proof}

We can combine   
Propositions~\ref{prop:alpha_e},
\ref{prop:triple-stable-implies-quotient-semistable} and 
\ref{prop:E2-stability} to obtain  the following.
\begin{proposition}
\label{prop:triple-stable-implies-bundles-semistable} 
Let $T=(E_1,E_2,\phi)$ be an  $\alpha$-semistable triple for
some $\alpha$\ in the range $\alpha_L<\alpha<\alpha_M$. Then
$T$ is of the form
$$
0 \lto E_2 \overset{\phi}{\lto} E_1 \lto F \lto 0, 
$$
with  $F$ locally free, and $E_2$ and $F$ are semistable.
\end{proposition}
In the converse direction we have:
\begin{proposition}
Let $T=(E_1,E_2,\phi)$ be a triple 
of  the form
$$
0 \lto E_2 \overset{\phi}{\lto} E_1 \lto F \lto 0, 
$$
with   $F$  locally free.
If  $E_2$ is semistable  and $F$ is stable then $T$ 
 is $\alpha$-stable for  $\alpha=\alpha_M-\epsilon$ in the range 
$\alpha_L<\alpha<\alpha_M$. 
\end{proposition}
\begin{proof}
Any subtriple $T'=(E_1',E_2',\phi')$ defines a commutative diagram
  $$
  \begin{CD}
  0@>>>E_2@>\phi>>E_1@>>>F@>>>0\\
  @.@AAA@AAA@AAA\\
  0@>>>E'_2@>\phi'>>E'_1@>>>F'@>>>0,
  \end{CD}
$$
where $F'\subset F$. Then
\begin{align}
\Delta_{\alpha}\equiv\Delta_{\alpha}(T') & = 
\mu_{\alpha}(T')-\mu_{\alpha}(T) \notag\\
&=\mu(E_1'\oplus E_2')-\mu(E_1\oplus E_2)+
\alpha(\frac{n_2'}{n_1'+n_2'}- \frac{n_2}{n_1+n_2}).
\end{align}
 Denote the rank and slope of $F$\ (resp.
  $F'$) by $m$\ and $\mu_F$\ (resp.\ $m'$\ and $\mu_{F'}$).
Using
\begin{align}
n_1&=n_2+m, \notag \\
n'_1&=n'_2+m', \notag \\
n_1\mu_1 &= n_2\mu_2+m\mu_F,\notag \\
n'_1\mu'_1 &= n'_2\mu_2'+m'\mu_{F'}, \notag
\end{align}
and the fact that $\alpha_M=2(\mu_F-\mu_2)$, and setting
 $\alpha=\alpha_M-\epsilon$,
we obtain 
\begin{align}
\Delta_{\alpha} & = 
\frac{2n_2'+\mu_2'+m\mu_{F'}}{2n_2'+m'}-\frac{2n_2+\mu_2+m\mu_F}{2n_2+m}
+ 2(\mu_F-\mu_2)(\frac{n_2'}{2n_2'+m'}- \frac{n_2}{2n_2+m}) \notag \\
& -\epsilon(\frac{n_2'}{2n_2'+m'}- \frac{n_2}{2n_2+m})\notag \\
&= \frac{2n_2'}{2n_2'+m'}(\mu_2'-\mu_2) + \frac{m'}{2n_2'+m'}(\mu_{F'}-\mu_F)
-\epsilon(\frac{n_2'm-n_2m'}{(2n_2'+m')(2n_2+m)}).  \label{alpha-epsilon}
\end{align}
Clearing denominators in (\ref{alpha-epsilon}) we obtain 
$$
(2n_2'+m')\Delta_{\alpha}=n_2'(2\Delta_2-\frac{m}{2n_2+m}\epsilon) +
                    m'(\Delta_F+\frac{n_2}{2n_2+m}\epsilon),
$$
where 
\begin{equation}
\Delta_2 =\mu_2'-\mu_2, \quad\ \mathrm{and}\ \quad\ 
\Delta_F=\mu_{F'}-\mu_F. 
\end{equation}
Now suppose that  $E_2$ is semistable and $F$ is stable. The 
semistability of $E_2$ implies that 
$$
2\Delta_2-\frac{m}{2n_2+m}\epsilon<0,
$$
while the  stability of $F$ implies there exists $\delta>0$ such that
$\Delta_F\leq -\delta <0$.
Thus by taking 
$\epsilon <\frac{2n_2 +m}{n_2}\delta$,
we have 
$\Delta_F+\frac{n_2}{2n_2+m}\epsilon<0$, 
and hence $\Delta_\alpha <0$.
\end{proof}
\begin{theorem}\label{thm:largealpha}
Assume that $n_1>n_2$ and $d_1/n_1> d_2/n_2$. Then the moduli space 
$\mathcal{N}^s_L=\mathcal{N}^s_L(n_1,n_2,d_1,d_2)$ 
is  smooth 
of dimension 
$$
(g-1)(n_1^2 + n_2^2 - n_1 n_2) - n_1 d_2 + n_2 d_1 + 1,
$$
and is birationally equivalent to a $\mathbb{P}^N$-fibration over
$M^s(n_1-n_2,d_1-d_2) \times M^s(n_2,d_2)$, where $M^s(n,d)$ is the
moduli space of stable bundles of rank $n$ and degree $d$, and 
$$N=n_2d_1-n_1d_2+n_1(n_1-n_2)(g-1)-1.$$
In particular, $\mathcal{N}_L^s(n_1,n_2,d_1,d_2)$ is non-empty and
irreducible. 

If $\GCD(n_1-n_2,d_1-d_2)=1$ and $\GCD(n_2,d_2)=1$, the birational 
equivalence is an isomorphism. 

Moreover, in all cases,  $\mathcal{N}_L=\mathcal{N}_L(n_1,n_2,d_1,d_2)$ is irreducible and 
hence  birationally equivalent to $\mathcal{N}_L^s$. 
\end{theorem}
\begin{proof}
For every triple $T=(E_1,E_2,\phi)$  in $\mathcal{N}^s_L$,
the homomorphism  $\phi$ is injective and hence,  by (5) in Proposition 
\ref{thm:smoothdim}, $T$ defines a smooth point in the moduli space, whose dimension is 
then given by (4) in Proposition \ref{triples-critical-range}. 

Given 
$F\in M^s(n_1-n_2,d_1-d_2)$ and $E_2\in M^s(n_2,d_2)$, we know
from Proposition \ref{triples-critical-range} that
every extension 
$$
0 \lto E_2 \overset{\phi}{\lto} E_1 \lto F \lto 0, 
$$
defines  a triple $T=(E_1,E_2,\phi)$ in  $\mathcal{N}^s_L$.
These extensions are  classified by
  $H^1(E_2\otimes F^*)$. In fact two classes defining the same element
in the projectivization  $\mathbb{P}H^1(E_2\otimes F^*)$ define equivalent
extensions and therefore equivalent triples.
Now, 
$$
\deg(E_2\otimes F^*)=(n_1-n_2)d_2-n_2(d_1-d_2)= n_1n_2(\mu_2-\mu_1)<0
$$
and, since $E_2\otimes F^*$\ is semistable, then
$H^0(E_2\otimes F^*)=0$.  Hence, by the Riemann--Roch theorem
$$h^1(E_2\otimes F^*)=n_2d_1-n_1d_2+n_1(n_1-n_2)(g-1).$$
In particular this dimension is constant as $F$ and $E_2$ vary in their
corresponding moduli spaces.

We  can describe this globally in terms of Picard sheaves. To do that  
we consider first the case in which $\GCD(n_1-n_2,d_1-d_2)=1$ and  
$\GCD(n_2,d_2)=1$. In this situation there exist universal 
 bundles $\mathbb{F}$  and  $\mathbb{E}_2$ over $X \times M(n_1-n_2,d_1-d_2)$
and $X\times M(n_2,d_2)$, respectively.
Consider the canonical projections 
\begin{align*}
  \pi&\colon X\times M(n_1-n_2,d_1-d_2)\times M(n_2,d_2)\times \to
  M(n_1-n_2,d_1-d_2) \times M(n_2,d_2)\ , \\
  \nu&\colon X \times M(n_1-n_2,d_1-d_2)\times M(n_2,d_2)\to 
  X \times  M(n_1-n_2,d_1-d_2)\ , \\
\intertext{and}  
  \pi_2&\colon X \times M(n_1-n_2,d_1-d_2)\times M(n_2,d_2)\to 
  X\times M(n_2,d_2)\ . 
\end{align*}
The Picard sheaf
$$
\mathcal{S}:=R^1\pi_*(\pi_2^*\mathbb{E}_2\otimes \nu^* \mathbb{F^*}),
$$
is then locally free  and we can identify 
$\mathcal{N}_L=\mathcal{N}_L^s$
with  $\mathcal{P}=\mathbb{P}(\mathcal{S})$. This is indeed a 
$\mathbb{P}^N$ fibration
with $N=n_2d_1-n_1d_2+n_1(n_1-n_2)(g-1)-1$, which in particular
is non-empty since $M(n_1-n_2,d_1-d_2)$ and  $M(n_2,d_2)$ are 
non-empty and $N>0$.

If $\GCD(n_1-n_2,d_1-d_2)\neq 1$ and  $\GCD(n_2,d_2)\neq 1$, there are no
universal bundles and hence the Picard bundle does not exist. However, its 
projectivization over 
$$M^s(n_1-n_2,d_1-d_2) \times M^s(n_2,d_2)$$
does exist. This can be constructed 
by working in  the open set $R$ of the Quot 
scheme corresponding to stable bundles. The point is that an 
appropriate linear group 
$\GL$ acts on $R$, with the centre acting trivially and such that 
$\PGL$ acts freely with the quotient being $M^s(n_1-n_2,d_1-d_2) 
\times M^s(n_2,d_2)$. For the action on the projective bundle 
associated to the universal bundle over $R$, the centre of $\GL$ 
still acts trivially, and  standard descent arguments produce the
required $\mathbb{P}^N$ fibration 
$\mathcal{P}$ over $M^s(n_1-n_2,d_1-d_2) \times M^s(n_2,d_2)$. 

We now show that the complement of $\mathcal{P}$ has strictly positive
codimension in $\mathcal{N}_L^s$.  This follows from two facts. The
first one is that any family of strictly semistable bundles of rank
$n_1-n_2$ and degree $d_1-d_2$ depends on a number of parameters
strictly less than the dimension of $M^s(n_1-n_2,d_1-d_2)$
(cf. e.g. \cite{brambila-grzegorczyk-newstead:1997}).
The same argument applies to
 any family of strictly semistable bundles of rank $n_2$ and degree
$d_2$. The second fact is that the dimension of $H^1(E_2\otimes F^*)$
is fixed by the Riemann--Roch theorem (we use here that $E_2$ and $F$
are semistable).

To prove the last statement, i.e.\ to extend the results to 
$\mathcal{N}_L$, we consider the family 
$\tilde{\mathcal{P}}$ of equivalence classes of extensions 
$$
0 \lto E_2 \overset{\phi}{\lto} E_1 \lto F \lto 0, 
$$
where $F$ and $E_2$ are semistable. Clearly, $\tilde{\mathcal{P}}$ 
contains the family $\mathcal{P}$. The family $\tilde{\mathcal{P}}$ 
is irreducible. This is because  since $F$ and $E_2$ are semistable  
they vary (for fixed ranks and degrees) in irreducible families 
$\mathcal{F}$ and $\mathcal{E}_2$, respectively, and  as shown above
$H^0(E_2\otimes F^*)=0$. Hence  
$\tilde{\mathcal{P}}$ is a projective bundle over $\mathcal{F}\times 
\mathcal{E}_2$.  {}From Proposition 
\ref{prop:triple-stable-implies-bundles-semistable}, 
we know that $\mathcal{N}_L\subset 
\tilde{\mathcal{P}}$, and since 
$\alpha$-semistability is an open condition 
(which follows  from the construction of the moduli space given
in   \cite{bradlow-garcia-prada:1996} and \cite{Schmitt}),
we have that
$\mathcal{N}_L$ is irreducible.
\end{proof}

\begin{remark}
If $n_1<n_2$, we have an analogous theorem for 
$\mathcal{N}^s_L(n_1,n_2,d_1,d_2)$ via the isomorphism
$$
\mathcal{N}^s_\alpha(n_1,n_2,d_1,d_2) = \mathcal{N}^s_\alpha(n_2,n_1,-d_2,-d_1)
$$
given by duality (Proposition \ref{prop:duality}).
\end{remark}

\subsection{Moduli space for $\alpha_m<2g-2\leq \alpha<\alpha_M$}

\begin{theorem}\label{thm:irreducibility-moduli-stable-triples}
Let $\alpha$ be any value in the range $\alpha_m<2g-2\leq\alpha< 
\alpha_M$. Then 
$\mathcal{N}^s_\alpha$ is birationally equivalent to 
$\mathcal{N}^s_L$.  In particular it is non-empty and irreducible. 
\end{theorem}
\begin{proof}
This follows from Theorem \ref{thm:birationality} and Theorem 
\ref{thm:largealpha}. 
\end{proof}

\begin{corollary}\label{cor:gcd=1}
Let $(n_1,n_2,d_1,d_2)$ be such that $\GCD(n_2,n_1+n_2,d_1+d_2)=1$.
If  $\alpha$ is a  generic value satisfying  
$\alpha_m<2g-2\leq\alpha<\alpha_M$, then  $\mathcal{N}_\alpha$ is birationally 
equivalent to $\mathcal{N}_L$, and in particular it is irreducible. 
\end{corollary}
\begin{proof}
{}From (4) in Proposition \ref{triples-critical-range} one has that 
$\mathcal{N}_\alpha=\mathcal{N}_\alpha^s$ if
$\GCD(n_2,n_1+n_2,d_1+d_2)=1$ and $\alpha$ is generic.
In particular, $\mathcal{N}_L=\mathcal{N}_L^s$, and hence the result
follows from Theorem \ref{thm:irreducibility-moduli-stable-triples}.
\end{proof}

\section {Moduli space of triples with $n_1 = n_2$}
\label{sec:n_1=n_2}
Throughout this section we  assume that $n_1=n_2=n$
and $d_1\geq d_2$.

\subsection{Moduli space for $d_1=d_2$}\label{subs:d_1=d_2}

\begin{proposition}\label{same-rank-same-slope}
Suppose  that $n_1=n_2=n$ and $d_1=d_2=d$. Let $T=(E_1,E_2,\phi)$ be 
a triple of type $(n_1,n_2,d_1,d_2)$, and let $\alpha>0$. Then $T$ is 
$\alpha$-(poly)stable if and only if $E_1$ and $E_2$ are (poly)stable and 
$\phi$ is an isomorphism. 
\end{proposition}
\begin{proof}
In this case the injectivity bound $\alpha_0$ given by 
(\ref{alpha0}) is  $\alpha_0=\alpha_m=0$. Hence for  every 
$\alpha$-semistable triple $T=(E_1,E_2,\phi)$ with $\alpha>0$, $\phi$ 
must be injective and therefore an isomorphism. The polystability of 
$E_1$ and $E_2$ is now  straightforward to see. To show the converse, 
suppose that $E_1$ and $E_2$ are both polystable and let 
$T'=(E_1',E_2',\phi')$ be any subtriple of $T$. 

\begin{align*}
\mu_\alpha(T')&= \mu(E_1'\oplus E_2')+\alpha\frac{n_2'}{n_1'+n_2'}\\
              &\leq  \mu(E_1\oplus E_2)+\alpha\frac{n_2'}{n_1'+n_2'}\\
              &\leq  \mu_\alpha(T)+
                \alpha(\frac{n_2'}{n_1'+n_2'}-\frac{1}{2})\\
                            &\leq  \mu_\alpha(T),
\end{align*}
since $n_1'\geq n_2'$ if $\phi$ is injective. 

\end{proof}

\begin{corollary}\label{cor:largealphamoduli}
The moduli space $\mathcal{N}_\alpha(n,n,d,d)$ and the moduli space
$M(n,d)$  of poly\-stable bundles of rank $n$ and degree $d$ 
are isomorphic.  In particular $\mathcal{N}_\alpha(n,n,d,d)$ is
non-empty and irreducible.
\end{corollary}
\begin{proof}
{}From Proposition \ref{same-rank-same-slope} it is clear that 
we have a surjective map, say 
$$\pi:\mathcal{N}_\alpha(n,n,d,d)\to M(n,d).$$
Suppose that 
$\pi([T])=\pi([T'])$, where $[T]$ and $[T']$ are points 
in $\mathcal{N}_\alpha(n,n,d,d)$ represented by triples 
$T=(E,E,\phi)$ and $T'=(E',E',\phi')$ respectively. We may assume 
that $T$ and $T'$ are polystable triples, and hence that $E$ and 
$E'$ are polystable vector bundles.  Thus, since 
$\pi([T])=\pi([T'])$, we can find an isomorphism 
$h_1:E\mapsto E'$. Set $h_2=\phi'\circ h_1\circ\phi^{-1}$  
(remember that $\phi$ and $\phi'$ are bundle isomorphisms). 
Then $(h_1,h_2)$ defines an isomorphism form $T$ to $T'$. 
Thus $\pi$ is injective.
\end{proof}

Combining Proposition \ref{moduli-alpha-M} and  Corollaries
 \ref{cor:N-alpha-M} 
and  \ref{cor:largealphamoduli},   we obtain the following.

\begin{corollary}
Suppose that $n_1>n_2$ and $\mu_1-\mu_2>0$. Then 
$$
\mathcal{N}_{\alpha_M}(n_1,n_2,d_1,d_2) \cong
M(n_2,d_2) \times M(n_1-n_2,d_1-d_2).
$$
In particular,  $\mathcal{N}_{\alpha_M}(n_1,n_2,d_1,d_2)$ is 
non-empty and irreducible. 
\end{corollary}

\subsection{Bounds on $E_1$ and $E_2$ for  $\alpha>\alpha_0$}

\begin{lemma}  \label{lemma:technical}
  Let $(E_1,E_2,\phi)$ be a triple with
  $\ker\phi=0$. Let $(E'_1,E'_2,\phi')$ be a subtriple with
  $n'_1=n'_2=n'$. Thus we get the following diagram, in which $S$ and
  $S'$ are torsion sheaves:
  $$
  \begin{CD}
  0@>>>E_2@>\phi>>E_1@>>>S@>>>0\\
  @.@AAA@AAA@AAA\\
  0@>>>E'_2@>\phi'>>E'_1@>>>S'@>>>0.
  \end{CD}
  $$
  Then
  \begin{align*}
  \Delta_{\alpha}(T')\equiv
  \mu_{\alpha}(T')-\mu_{\alpha}(T)
  &=(\mu(E'_2)-\mu_2)+\frac{1}{2}\left(\frac{s'}{n'}-\frac{s}{n}\right)\\  &=(\mu(E'_1)-\mu_1)-\frac{1}{2}\left(\frac{s'}{n'}-\frac{s}{n}\right).
  \end{align*}
  Here $s$\ and $s'$\ are the degrees of $S$\ and $S'$\ respectively.
\end{lemma}
\begin{proof} {}From the above diagram we get
\begin{align*}
 n\mu_2+s&=n\mu_1,\\
 n\mu(E'_2)+s'&=n\mu(E'_1).
\end{align*}
Thus
\begin{align*}
\mu_{\alpha}(T')=
\frac{1}{2}(\mu(E'_1)+\mu(E'_2))+\frac{\alpha}{2}
&=\frac{1}{2}\left(2\mu(E'_2) + \frac{s'}{n'}\right)+\frac{\alpha}{2}\\
&=\frac{1}{2}\left(2\mu(E'_1)- \frac{s'}{n'}\right)+\frac{\alpha}{2},
\end{align*}
and similarly for $\mu_{\alpha}(T)$.
\end{proof}
\begin{proposition}
  \label{prop:slope-bound-on-subbdls}
  Let $(E_1,E_2,\phi)$\ be an $\alpha$-semistable triple with 
$\ker\phi=0$. Then
\begin{itemize}
\item[$(1)$] For any subsheaf $E'_1\subset E_1$
  $$
  \mu(E'_1)\leq\mu_1+\frac{1}{2}(n-1)(\mu_1-\mu_2).
  $$
\item[$(2)$] For any subsheaf $E'_2\subset E_2$
  $$
  \mu(E'_2)\leq\mu_2+\frac{1}{2}(\mu_1-\mu_2).
  $$
\end{itemize}
\end{proposition}
\begin{proof} Since $\ker\phi=0$\ the results of Lemma
  \ref{lemma:technical} apply. Furthermore, any subsheaf $E'_1\subset
  E_1$\ is part of a subtriple $(E'_1, E'_2,\phi')$\ with
  $n'_1=n'_2=n'$. Likewise, given any subsheaf $E'_2\subset E_2$, we
  can take $E'_1=\phi(E'_2)$. Thus we can use the results of Lemma
  \ref{lemma:technical}, plus the fact that $\alpha$-stability implies
  $\Delta_{\alpha}(T')<0$\ for all subtriples, to conclude
  $$
  \mu(E'_1)-\mu_1-\frac{1}{2}\left(\frac{s'}{n'}-
  \frac{s}{n}\right)<0\
  $$
  for all $E'_1\subset E_1$. Similarly
  $$
  \mu(E'_2)-\mu_2+\frac{1}{2}\left(\frac{s'}{n'}-
    \frac{s}{n}\right) < 0
  $$
  for all $E'_2\subset E_2$.  The results now follow using the fact
  that $0\leq s'\leq s$\ and $1\leq n'< n$.
\end{proof}

\subsection{Stabilization of moduli}
\label{sec:stabilize}

\begin{theorem}[Stabilization Theorem]
  \label{thm:stabilization}
Let $\alpha_0$ be as in  (\ref{alpha0-equal-ranks}).  
\begin{itemize}
\item[$(1)$] Let   $\alpha_1$, $\alpha_2$\ be any real numbers
  such that $\alpha_0<\alpha_1\leq\alpha_2$, then
$$
\mathcal{N}_{\alpha_1}(n,n,d_1,d_2)\subseteq\mathcal{N}_{\alpha_2}(n,n,d_1,d_2). 
$$
\item[$(2)$] There is a real number
  $\alpha_{L}\geq\alpha_0$\ such that
$$
\mathcal{N}_{\alpha_1}(n,n,d_1,d_2)
=\mathcal{N}_{\alpha_2}(n,n,d_1,d_2)
$$
for all $\alpha_1\geq \alpha_2>\alpha_L$.
\end{itemize}
\end{theorem}
\begin{proof}
  (1). Recall
 from Proposition \ref{injective-triples} that if $\alpha>\alpha_0$\ then 
 any triple, 
  $T=(E_1,E_2,\phi)$,
  in $\mathcal{N}_\alpha(n,n,d_1,d_2)$\ has $\rk(\phi)=n$. It follows
  that in any subtriple, say $T'=(E'_1,E'_2,\phi')$, the rank of $E'_1$\ 
  is at least as big as the rank of $E'_2$, i.e.\ $n'_1\geq n'_2$. We
  treat the cases $n'_1> n'_2$\ and $n'_1= n'_2$\ separately. In both
  cases we must show that
$$
\Delta_{\alpha_1}(T')\leq0\ \Rightarrow\
\Delta_{\alpha_2}(T')\leq0\ 
$$
if $\alpha_1\leq\alpha_2$.
If $n'_1=n'_2$\ then for any $\alpha$
\begin{equation}
\Delta_{\alpha}(T')=
\mu(E'_1\oplus E'_2)-\mu(E_1\oplus E_2).\label{eq:2.60}
\end{equation}
In particular,  $\Delta_{\alpha}(T')$\ is independent of
$\alpha$\ and hence
$\Delta_{\alpha_1}(T')=
\Delta_{\alpha_2}(T')$.
If $n'_1>n'_2$, then for any $\alpha$
\begin{align}
\Delta_{\alpha}(T')&=
\mu(E'_1\oplus E'_2)-\mu(E_1\oplus E_2)+
\left(\frac{n'_2}{n'_1+n'_2}-\frac{1}{2}\right)\alpha. \label{eq:2.61}
\end{align}
For each subtriple, $\Delta_{\alpha}(T')$\ is thus a linear
function of $\alpha$, with slope
\begin{equation}
\lambda(T')\ =\ \left(\frac{n'_2}{n'_1+n'_2}-\frac{1}{2}\right)\ =\
\frac{n'_2-n'_1}{2(n'_1+n'_2)}\ \label{eq:2.62}
\end{equation}
and constant term
\begin{equation}
M(T')= \mu(E'_1\oplus E'_2)-\mu(E_1\oplus E_2).\label{eq:2.63}
\end{equation}
We see that if $n'_1>n'_2$\ then $\lambda(T')<0$. It follows
from this that
$$
\Delta_{\alpha_1}(T')\leq 0\ \Longrightarrow\
\Delta_{\alpha_2}(T')\leq 0\ 
$$
if $\alpha_1\leq\alpha_2$.

(2). Consider any $\alpha_1,\alpha_2$\ such that
$\alpha_0<\alpha_1\leq\alpha_2$. By Part (1), the difference (if any)
between $\mathcal{N}_{\alpha_1}$\ and $\mathcal{N}_{\alpha_2}$\ is due
entirely to triples which are $\alpha_2$-stable but not
$\alpha_1$-stable. Any such triple must have a subobject, say
$T'=(E'_1,E'_2,\phi')$, for which
\begin{equation}
\Delta_{\alpha_2}(T')\leq 0 <
\Delta_{\alpha_1}(T').\label{eq:2.64}
\end{equation}
As in (1), we need only consider subobjects for which the rank of
$E'_1$\ is at least as big as the rank of $E'_2$, i.e.\ $n'_1\geq
n'_2$. Clearly \eqref{eq:2.64} is not possible for a subobject with
$n'_1=n'_2$\ (since in that case
$\Delta_{\alpha_1}(T')=\Delta_{\alpha_2}(T')$). Suppose
then that $n'_1>n'_2$. By \eqref{eq:2.61} and the fact that for such a
subobject $\lambda(T')<0$, we get that
\begin{equation}
\Delta_{\alpha}(T') \geq0\  \iff
\alpha \leq \frac{M(T')}{-\lambda(T')}.\label{eq:2.64a}
\end{equation}
We claim that there is a bound, $\alpha_L$, depending
  only on the degrees and ranks of $E_1$\ and $E_2$, such that
  \begin{equation}
    \frac{M(T')}{-\lambda(T')}\leq \alpha_L
    \label{eq:2.65}
  \end{equation}
  for all possible subtriples with $n'_1>n'_2$.
For a triple  $T=(E_1,E_2,\phi)$  in
$\mathcal{N}_{\alpha_2}$ 
Proposition~\ref{prop:slope-bound-on-subbdls} applies, giving
upper bounds on slopes of subsheaves of both $E_1$\ and $E_2$. Using
these bounds we compute
\begin{equation}
M(T')\leq
\frac{nn'_1}{2(n'_1+n'_2)}(\mu_1-\mu_2).\label{eq:2.66}
\end{equation}
Combined with \eqref{eq:2.62}, this gives
\begin{align*}
\frac{M(T')}{-\lambda(T')}&\leq
\frac{nn'_1}{(n'_1-n'_2)}(\mu_1-\mu_2)\\
&\leq n(n-1)(\mu_1-\mu_2).
\end{align*}
We can thus take
\begin{equation}
\alpha_L=n(n-1)(\mu_1-\mu_2).\label{eq:2.67}
\end{equation}
We can now complete the proof of Part (2): if
$\alpha_1>\alpha_L$\ then no triple in
$\mathcal{N}_{\alpha_2}$\ can have a subtriple satisfying
\eqref{eq:2.64a}.  Hence
$\mathcal{N}_{\alpha_2}=\mathcal{N}_{\alpha_1}$.
\end{proof}
\begin{remark}
If $n=2$ then  $\alpha_L=\alpha_0=d_1-d_2$, i.e.\ the stabilization
parameter coincides with the injectivity  parameter.
\end{remark}
It is clear from (\ref{eq:2.67}) that $\alpha_L=0$ correspond to the 
following especial cases. 
\begin{proposition} If $n=1$ or $\alpha_m=0$ then $\alpha_L=0$. 
Hence if $\epsilon$ is any positive real 
number,
\begin{itemize}
\item[(1)]
if $n=1$, then $\mathcal{N}_{\alpha}$ is 
isomorphic to $\mathcal{N}_{\alpha_m+\epsilon}(1,1,d_1,d_2)$ 
for every $\alpha\in (\alpha_m,\infty)$; 
\item[(2)] if $\alpha_m=0$, then $\mathcal{N}_{\alpha}$ is
isomorphic to 
$\mathcal{N}_{\epsilon}(n,n,d_1,d_2)$ for every $\alpha\in (0,\infty)$.
\end{itemize}
\end{proposition}

\subsection{Moduli for large $\alpha$ and $\alpha\geq 2g-2$}
\label{sec:alpha0}

Let $\alpha>\alpha_0$, where $\alpha_0$ is as in (\ref{alpha0-equal-ranks}). 
By Proposition
\ref{injective-triples}, we know that all triples in
$\mathcal{N}_\alpha(n,n,d_1,d_2)$ are of the form
\begin{equation}\label{ext-tor}
0 \lto E_2 \overset{\phi}{\lto} E_1 \lto S \lto 0,
\end{equation}
where $S$ is a torsion sheaf of degree $d=d_1-d_2$.
By analogy with the $n_1\neq n_2$, let us denote by 
$\mathcal{N}_L(n,n,d_1,d_2)$ the `large $\alpha$' moduli space, i.e.\ 
the moduli space of $\alpha$-semistable triples for any 
$\alpha\in (\alpha_L,\infty)$. Since $\alpha_L\geq \alpha_0$ we have 
 that all triples in
$\mathcal{N}_L(n,d_1,d_2)$ are of the form in (\ref{ext-tor})
and that $E_1$\ and $E_2$\ are bounded by the constraints in
Proposition \ref{prop:slope-bound-on-subbdls}.  

In the converse direction we have the following.
\begin{proposition}
  \label{prop:bundles-stable-implies-triple-stable}
Let $T=(E_1,E_2,\phi)$ be a triple such that $\ker\phi=0$
If $E_1$ and $E_2$ are semistable, then $T$ is  $\alpha$-semistable
for large enough  $\alpha$,  i.e.
  $T \in \mathcal{N}_L(n,n,d_1,d_2)$.  If
  either $E_1$\ or $E_2$\ is stable, then $T$\ is
  $\alpha$-stable.  
\end{proposition}
\begin{proof} 
  Since $\ker\phi=0$, it follows (as in the
  proof of Theorem \ref{thm:stabilization}) that in any subtriple, say
  $T'=(E'_1,E'_2,\phi')$, the rank of $E'_1$\ is at least as big as the
  rank of $E'_2$, i.e.\  $n'_1\geq n'_2$.  If $n'_1>n'_2$, then 
  \eqref{eq:2.61}, \eqref{eq:2.62} and  \eqref{eq:2.63} apply, with
  $\lambda(T')<0$\ and
  $\frac{M(T')}{-\lambda(T')}\leq \alpha_L$. 
 Thus $\mu_\alpha (T')-\mu_\alpha (T)<0$\ 
  for $\alpha>\alpha_L$. For subtriples with $n'_1=n'_2$, 
  equation \eqref{eq:2.60} says that
  $$
  \Delta_{\alpha}(T')= \mu(E'_1\oplus E'_2)-\mu(E_1\oplus
  E_2)\ 
  $$
for any $\alpha$. For such subtriples, and for any $\alpha$, it thus 
follows that 
\begin{itemize}
\item[$(1)$] $\Delta_{\alpha}(T')\leq 0$\ if both $E_1$\ and $E_2$\ are
  semistable, and
\item[$(2)$] $\Delta_{\alpha}(T') < 0$\ if at least one of the bundles is 
stable.
\end{itemize}
\end{proof}
\goodbreak

\begin{theorem} \label{thm:existence-n1=n2}
The moduli space $\mathcal{N}^s_L(n,n,d_1,d_2)$ is non-empty.
\end{theorem}
\begin{proof}
Our strategy is to  show that  there exist rank $n$ 
 stable bundles $E_1$ and $E_2$ of degree $d_1$ and $d_2$, respectively,
and a torsion sheaf $S$ of degree $d_1-d_2$,  fitting  in an exact sequence
$$
0\lto E_2\lto  E_1\lto S \lto  0.
$$ 
The result will then  follow from Proposition \ref{prop:bundles-stable-implies-triple-stable}.

To prove this, let $E$ be a vector bundle, and let $\Quot^d(E)$ be 
the Quot scheme of quotients $E\lto S$ where $S$ is a torsion
sheaf of degree $d$. The basic fact we  need  is the  following.

\begin{lemma}\label{hernandez}
Let $L$ be a line bundle and let  $\psi:L^{\oplus n}\lto S$ be an element in
 $\Quot^d(L^{\oplus n})$. Then, if $L$ has  big enough degree
(depending on $n$ and $d$),
 for a  generic  $S$, the vector bundle
$E=\ker\psi$ is stable.
\end{lemma}
\begin{proof}
The proof is implicit in the papers by
Hern\'andez \cite{hernandez} and Maruyama \cite{maruyama}, where they deal 
with the case $L=\mathcal{O}$. There, one needs an extra condition
on $n$ and $d$,  which is not required  in the twisted case when the degree
of $L$ is big enough.
\end{proof}  

Let $L$ be a line bundle of degree $m$ and  $d''>0$ such that
$d_1=nm-d''$.
By Lemma \ref{hernandez}, if  $\psi:L^n\lto S''\in \Quot^{d''}(L^n)$
is generic, then $E_1=\ker\psi$ is a stable bundle of rank $n$ and degree
$d_1$. Let $d=d_1-d_2$ and consider a generic element 
$\eta: E_1\lto S\in \Quot^d(E_1)$. Let $E_2=\ker\eta$, and let $S'$ the
cokernel of the natural inclusion $E_2\lto L^n$. We have the
following commutative diagram:

$$\begin{array}{ccccccccc}
&&0&&0&& 0 &&\\
&&\downarrow&&\downarrow& & \downarrow&  & \\
0&\rightarrow&E_2&\rightarrow&E_1 &\rightarrow& S & \rightarrow & 0 \\
&&\parallel &&\downarrow& &\downarrow&& \\
0&\rightarrow& E_2&\rightarrow& L^n  &\rightarrow& S'  &\rightarrow& 0\\
&&\downarrow&&\downarrow& &\downarrow&& \\
& &0&\rightarrow& S''& = & S'' &\rightarrow& 0\\
&& &&\downarrow& &\downarrow&& \\
 & & & & 0  & & 0.  & &  \\
\end{array}
$$

We see from the diagram that  $E_2$ coincides with the kernel
of $L^n\lto S'$. If $S'$ is general enough we can again apply
Lemma \ref{hernandez}   and conclude that $E_2$ is stable.
To show that this is indeed the case, we observe that the diagram
defines a map
$$
\Quot_0^d(E_1) \times  \Quot_0^{d''}(L^n) \lto \Quot_0^{d+d''}(L^n),   
$$
where $\Quot_0$ denotes an open non-empty subscheme of $\Quot$,  which
is surjective and finite.
\end{proof}

We deal now with the irreducibility of the moduli spaces.

\begin{theorem}[Markman-Xia \cite{markman-xia:2001}]
\label{thm:irreducible-alpha0}
There is an irreducible  family $\mathcal{S}$ parameterizing
quotients $E_1\lto S\lto 0$, where $E_1$ is a rank $n$ and degree
$d_1$ vector bundle such that all its subbundles
have their slope bounded above by a given universal constant, and
$S$ is a torsion sheaf of degree $d>0$.
\end{theorem}

\begin{theorem}\label{thm:irreducibility-alpha0}
If $\alpha>\alpha_0$, then $\mathcal{N}_\alpha(n,n,d_1,d_2)$ is
irreducible.
\end{theorem}
\begin{proof}
Since $\alpha>\alpha_0$, an $\alpha$-semistable triple
$T=(E_1,E_2,\phi)$ defines a sequence as in (\ref{ext-tor}) and hence a
quotient $ E_1 \lto S \lto 0$ in $\mathcal{S}$. 
By (1) in Proposition
\ref{prop:slope-bound-on-subbdls}, the slopes of subbundles of
$E_1$ are  bounded above by the universal constant
$$
\mu_1 +\frac{1}{2}(n-1)(\mu_1-\mu_2).
$$  
Let $\mathcal{S}^0\subset \mathcal{S}$
consist of  elements $E_1\lto S\lto 0$ in $\mathcal{S}$ that come
from an $\alpha$-semistable triple $0\to E_2\to E_1\to\mathcal{S}\to 0$.
Since $\alpha$-semistability is an open
condition $\mathcal{S}^0$ is a Zariski open set of $\mathcal{S}$, which
is  non-empty by Theorem \ref{thm:existence-n1=n2}, and
hence irreducible by Theorem \ref{thm:irreducible-alpha0}.
The irreducibility of
$\mathcal{N}_\alpha(n,n,d_1,d_2)$ follows now from that
of $\mathcal{S}^0$.
\end{proof}

\begin{proposition}\label{prop:nnNL}
The moduli space $\mathcal{N}_L(n,n,d_1,d_2)$ is 
birationally  equivalent to a 
$\mathbb{P}^N$-fibration $\mathcal{P}$ over
$M^s(n,d_2)\times \Sym^d(X)$, where $N=n(d_1-d_2)-1$, $\Sym^d(X)$ is the
$d$-symmetric product of $X$, and $M^s(n_2,d_2)$ is the moduli
space of stable bundles of rank $n_2$ and degree $d_2$.
\end{proposition}

\begin{proof}
Let $E_2$ be a 
rank 
$n$ and degree $d_2$ vector bundle  and let 
$S$ be  a torsion sheaf of degree $d>0$. We construct $E_1$ as an 
extension 
\begin{equation}\label{torsion-extension}
0 \lto E_2 \lto E_1 \lto S \lto 0.
\end{equation}
Such extensions are parameterized by $\Ext^1(S,E_2)$. Suppose
that   $S$ 
is of the form $S=\mathcal{O}_D$, where $D$ is a divisor
in $\Sym^d(X)$. 
Let $L$  be a line bundle. Consider  the  short exact sequence 
$$
0 \lto L^*(-D) \lto L^* \lto \mathcal{O}_D \lto 0,
$$
and apply to it the functor $\Hom(\cdot,E_2)$, to  obtain the long
exact sequence
\begin{equation}
  \label{eq:long-exact-torsion}
\begin{array}{ccccccc}
  0 &\lto &H^0(E_2\otimes L) &\lto& H^0(E_2\otimes L(D)) &\lto&  \\ 
  \Ext^1(\mathcal{O}_D,E_2) & \lto & H^1(E_2\otimes L) & \lto& 
H^1(E_2\otimes L(D)) &\lto &  0.
\end{array}
\end{equation}
We thus have
$$
\dim \Ext^1(\mathcal{O}_D,E_2)=\chi(E_2\otimes L)-\chi(E_2\otimes L(D))=nd,
$$
where $\chi(E)=\dim H^0(E)- \dim H^1 (E)$.
Taking $L$ so that $\deg(L)>>0$, we have that 
$H^1(E_2\otimes L)=0$. If $E_2$ is semistable (or more generally, if
it moves in a  bounded family)  we can take the  same $L$ 
for every $E_2$. Then
$$
\Ext^1(\mathcal{O}_D,E_2)=H^0(E_2\otimes L(D))/H^0(E_2\otimes L).
$$
Let $\mathcal{P}$ be the set of equivalence classes of
extensions \eqref{torsion-extension}, where $E_2$ is stable
then $\mathcal{P}$ is a  $\mathbb{P}^N$-fibration over
$M^s(n,d_2)\times \Sym^d(X)$, where  
$N=nd-1=\dim \mathbb{P}(\Ext^1(\mathcal{O}_D,E_2))$. 
Since we assume that $d$ is positive, $N$ is non-negative and positive
if $n>1$.
Setting $d=d_1-d_2$, a simple 
computation shows that
$$
\dim \mathcal{P}=
(g-1)(n_1^2 + n_2^2 - n_1 n_2) - n_1 d_2 + n_2 d_1 + 1.
$$
Clearly $\mathcal{P}$ is irreducible of the same dimension as $\mathcal{N}_L$,
and since it is contained in $\mathcal{S}$ (like $\mathcal{N}_L$) it
must be  birationally equivalent to $\mathcal{N}_L$.
Notice that  If $\GCD(n,d_2)=1$, then 
$\mathcal{P}$ is the projectivization of a  Picard bundle.
\end{proof}

Combining the results of this section, we arrive at the following 
theorems.

\begin{theorem}\label{thm:moduli-n1=n2}
The moduli space $\mathcal{N}_L(n,n,d_1,d_2)$ is non-empty and irreducible. 
Furthermore, it is
birationally  equivalent to a 
$\mathbb{P}^N$-fibration  over
$M^s(n,d_2)\times \Sym^d(X)$, where the fiber dimension is
$N=n(d_1-d_2)-1$.
\end{theorem}
\begin{proof}
It follows from Theorems \ref{thm:irreducibility-alpha0}
and \ref{thm:existence-n1=n2} and Proposition
\ref{prop:nnNL}.
\end{proof}

\begin{theorem}\label{thm:irreducibility-moduli-stable-triples-n1=n2}
 Let $ \alpha\ge 2g-2>\alpha_m$.  Then 
\begin{itemize}
\item[$(1)$] 
The moduli space  $\mathcal{N}^s_\alpha$ is  
 birationally equivalent 
to $\mathcal{N}_L$  and it is hence non-empty and  irreducible.

\item[$(2)$] 
If in addition either
\begin{itemize}
\item [$\bullet$] $\GCD(n,2n,d_1+d_2)=1$ and $\alpha\ge 2g-2>\alpha_m$ is  generic, or
\item [$\bullet$] $d_1-d_2<\alpha$,
\end{itemize}
then  $\mathcal{N}_{\alpha}(n,n,d_1,d_2)$ is birationally equivalent 
to $\mathcal{N}_L(n,n,d_1,d_2)$ and hence irreducible. 
\end{itemize}
\end{theorem}
\begin{proof}

If $2g-2>\alpha_L$, the result follows from Theorems \ref{thm:stabilization} and \ref{thm:moduli-n1=n2}.
Assume then that  $2g-2\leq\alpha_L$.

(1) {}From Theorem \ref{thm:irreducibility-alpha0} we 
know that $\mathcal{N}_L$ is birationally equivalent to 
$\mathcal{N}_L^s$. The result  follows now from Theorem
\ref{thm:birationality} and Theorem \ref{thm:moduli-n1=n2}.

(2)  For the first  part, we observe that 
 from (4) in Proposition \ref{triples-critical-range} one has that 
$\mathcal{N}_\alpha=\mathcal{N}_\alpha^s$ if
$\GCD(n,2n,d_1+d_2)=1$ and $\alpha$ is generic, and hence the result
follows from (1). The second  part  is a consequence of  Theorem \ref{thm:irreducibility-alpha0}.
\end{proof}

 \section{Triples and dimensional reduction}\label{sec:reduction}

Let $\dP^1$ be the complex projective line. 
The Lie group $\SL(2,\C)$\ acts on  $\xp$
via the trivial action on $X$ and the 
identification $\dP^1=\SL(2,\C)/P$, where $P$ is the subgroup of
lower triangular matrices. 

The theory of holomorphic triples and vortex equations on $X$ is  related 
with the study of stable $\SL(2,\C)$-equivariant bundles on $\xp$ and  
the existence of invariant solutions to the Hermitian--Einstein equations.
In fact, it is in this way (known as dimensional reduction) that  the theory  
originated (see \cite{bradlow-garcia-prada:1996} and 
\cite{garcia-prada:1994} for details).

In this section we recall the basics of this correspondence and apply 
our main results on triples  to the theory of vector bundles on $\xp$.

 \subsection{Existence of stable  bundles on $\xp$ and triples}

 \begin{proposition}\cite[Proposition 2.3]{bradlow-garcia-prada:1996}
\label{prop:triple-extension}
There is a one-to-one correspondence between 
holomorphic triples $(E_1,E_2,\phi)$ 
on $X$ and  holomorphic extensions on $\xp$ of
the form
 \begin{equation}
 0\lra \ps E_1\lra E \lra \ps E_2\otimes\qs\cO(2)\lra 0,  \label{extension}
 \end{equation}
where $p$ and $q$ are the canonical projections from $\xp$ to $X$ and $\dP^1$,
respectively, and $\cO(2)$ is the degree 2 line bundle of $\dP^1$ (the tangent
bundle).
\end{proposition}
\begin{proof}
The proof given  in \cite{garcia-prada:1994} is simply that extensions over $\xp$  
of the form (\ref{extension}) are
parametrized by
 $$
 H^1(\xp, \ps (E_1\otimes E_2^\ast)\otimes \qs\cO(-2)).
 $$
By  the K\"unneth formula,  this is isomorphic to
 $$
  H^0(X, E_1\otimes E_2^\ast)\otimes H^1(\dP^1,\cO(-2))
 \cong H^0(X,E_1\otimes E_2^\ast).
 $$
After fixing an element in $H^1(\dP^1,\cO(-2))$, the
 homomorphism $\phi$ can thus be identified with the extension class defining
$E$.
\end{proof}

Notice that the vector bundles $E$ of the form in (\ref{extension})
come equipped with an action of $\SL(2,\C)$ which lifts 
the action on $\xp$. The action on $E$ is trivial 
on $\ps E_1$ and $\ps  E_2 $, is the standard one on $\cO(2)$,
and leaves invariant the extension class.

To talk about the stability of $E$ one needs a \kahler\ metric on $\xp$.
Let us fix a metric on $X$ and the  Fubini-Study metric on $\dP^1$, both  
normalized to have volume $2\pi$. Let $\alpha>0$ be a real number.
We  consider on $\xp$ the  one-parameter family of $\SU(2)$-invariant \kahler\ metrics
with \kahler\ form
 $$
 \omega_\alpha=\alpha\ps\omega_X\oplus\qs\omega_{\dP^1}.
 $$
Here $\omega_X$ and $\omega_{\dP^1}$ are the Kahler forms on $X$ and 
$\dP^1$, respectively. The degree
of a complex vector  bundle $E$ over $\xp$ with respect to $\omega_\alpha$ is given
by
 $$
 \deg(E)=\int_{\xp} c_1(E)\wedge\omega_\alpha,
 $$
 where $c_1(E)$ is the first Chern class of $E$.
Recall that  $E$  is said to be {\em stable} with respect to $\omega_\alpha$  if for
every non-trivial coherent reflexive subsheaf $E'\subset E$,
 $$
 \mu(E')<\mu(E),
 $$
where $\mu(E)=\deg E /\rk E $\ is the {\em slope} of $E $.
Since we are in complex dimension 2,  $E'$ is locally free.

 \begin{theorem} \cite[Theorem 4.1]{bradlow-garcia-prada:1996}
\label{tsvs}
 Let $T=(E_1,E_2,\phi)$ be a  holomorphic triple over
 $X$ and let $E$  be the holomorphic bundle over $\xp$ defined by $T$ as in
Proposition \ref{prop:triple-extension}.
Then,  if $E_1$\ and $E_2$\ are not isomorphic, $T$ is $\alpha$-stable if and only 
if  $E$ is stable  with respect to $\omega_\alpha$.
If $E_1\cong E_2$, the triple $T$ is
$\alpha$-stable  if and only if  $E$ decomposes as a direct sum
 $$
 E=\ps E_1\otimes \qs\cO(1)\oplus \ps E_2\otimes \qs\cO(1),
 $$
 and $\ps E_i\otimes \qs\cO(1)$
is stable with respect to $\omega_\alpha$.
 \end{theorem}
\begin{remark}
The stability of $\ps E_i\otimes \qs\cO(1)$ is equivalent to the stability
of $E_i$. 

\end{remark}
Let $(n_1,n_2,d_1,d_2)$ be  the type of the triple $T=(E_1,E_2,\phi)$.
Let $\mathcal{M}_\alpha$ be the moduli space of stable bundles on $\xp$
with respect to $\omega_\alpha$, whose topological type is that of $E$ in
 (\ref{extension}).
Combining Theorems \ref{tsvs},  \ref{thm:irreducibility-moduli-stable-triples} 
and \ref{thm:irreducibility-moduli-stable-triples-n1=n2} we can prove existence
of stable bundles on $\xp$. More precisely.

\begin{theorem} \label{thm:existence-stable-vb}
$\mathcal{M}_\alpha$ is non-empty if
\begin{itemize}
\item[(1)] $2g-2\leq \alpha\leq \alpha_M$ if $n_1\neq n_2$, where $\alpha_M$ is
given by (\ref{alpha-bounds-bigM});

\item[(2)] $2g-2\leq \alpha$ if $n_1 = n_2$.
\end{itemize}
\end{theorem}

\begin{remark} 
The moduli space $\mathcal{N}_\alpha$ can be identified with
the $\SL(2,\C)$-invariant part of $\mathcal{M}_\alpha$ 
(\cite{bradlow-garcia-prada:1996}). 
Hence
from Theorems \ref{thm:irreducibility-moduli-stable-triples} and 
\ref{thm:irreducibility-moduli-stable-triples-n1=n2} we can say that
within the range for $\alpha$ in Theorem \ref{thm:existence-stable-vb}
the invariant loci  for different values of $\alpha$ are birationally
equivalent. Whether this is true or not  for the whole moduli spaces
 $\mathcal{M}_\alpha$
for different values of $\alpha$ is something that deserves further  study (see 
\cite{Q} for a discussion
on this in the rank two case).
\end{remark}

\subsection{Existence of $\SU(2)$-invariant Hermitian--Einstein metrics}

By the Hitchin--Kobayashi correspondence proved by Donaldson, Uhlenbeck and Yau
 \cite{D1,D2,UY}, the stability 
of the bundle $E$ on $\xp$ is equivalent to the existence of an  irreducible 
solution to the  Hermitian--Einstein equation. Recall that this is a
Hermitian metric on $E$ such that 
\begin{equation}
 \sqrt{-1}\Lambda F(E)=\mu  \Id_E,
 \label{he}
\end{equation}
where, as usual, $\Lambda$ is contraction with the \kahler\ form
of $\xp$,   $F(E)$ is the curvature of the unique  connection determined by 
the Hermitian metric and the holomorphic structure of $E$,
$\Id_E$  is the identity endomorphism of $E$  and $\mu$ is the slope of $E$.

The action of  $\SL(2,\C)$ on $E$ restricts to an action of the compact subgroup 
$\SU(2)\subset \SL(2,\C)$, and, since
the metric $\omega_\alpha$ on $\xp$ is $\SU(2)$-invariant, one can  
consider  $\SU(2)$-invariant solutions to (\ref{he}). 
The relevant  fact is the 
following.

 \begin{proposition} \cite[Proposition 3.11]{garcia-prada:1994}
\label{dr} Let $T=(E_1,E_2,\phi)$ be a
holomorphic triple of type $(n_1,n_2,d_1,d_2)$ over $X$  and let $E$ be the 
 holomorphic bundle over $\xp$ associated to $T$  by 
Proposition \ref{prop:triple-extension}.
Let  $\tau_1$ and $\tau_2$ be real numbers such that
$d_1+d_2=n_1\tau_1+n_2\tau_2$, and $\tau_1-\tau_2>0$.
 Then $T$ admits a solution to (\ref{eq:coupled-vortex}) if and only
 if $E$ admits an $\SU(2)$-invariant  \he\ metric with respect to $\omega_\alpha$.
 \end{proposition}

Combining the previous results  we have the following.

\begin{corollary} 
The vector bundle $E$ associated to a triple $T$ of type $(n_1,n_2,d_1,d_2)$
 has a Hermitian--Einstein metric, with respect to $\omega_\alpha$ if 
\begin{itemize}
\item[(1)] $2g-2\leq \alpha \leq \alpha_M$ if $n_1\neq n_2$, where $\alpha_M$ is
given by (\ref{alpha-bounds-bigM});

\item[(2)] $2g-2\leq \alpha$ if $n_1 = n_2$.

\end{itemize}

In fact this metric is $\SU(2)$-invariant and it is given by a vortex solution on $T$.
\end{corollary}

\begin{remark}
This is similar in spirit to the instanton solutions of vortex type on $\R^4$
studied by Witten in \cite{W} and Taubes \cite{Ta}.
\end{remark}



\begin{thebibliography}{99}

\bibitem{alvarez-garcia-prada:2001}
{\'A}lvarez-C{\'o}nsul, L., Garc{\'\i}a-Prada, O.: Dimensional reduction,
  $\mathrm{SL}(2,\mathbb{C})$-equivariant bundles and stable holomorphic
  chains.
\newblock Int. J. Math. \textbf{12}, 159--201 (2001).

\bibitem{AG2}{\'A}lvarez-C{\'o}nsul, L., Garc{\'\i}a-Prada, O.:
Hitchin--Kobayashi correspondence, quivers and vortices.
\newblock  Comm. Math. Phys. (to appear). arXiv:math.DG/0112161.

\bibitem{biswas-ramanan:1994}
Biswas, I., Ramanan, S.: An infinitesimal study of the moduli of {H}itchin
  pairs.
\newblock J. London Math. Soc. (2) \textbf{49}, 219--231 (1994).

\bibitem{bradlow-garcia-prada:1996}
Bradlow, S.~B., Garc{\'\i}a-Prada, O.: Stable triples, equivariant bundles
and
  dimensional reduction.
\newblock Math. Ann. \textbf{304}, 225--252 (1996).

\bibitem{bradlow-garcia-gothen:2001}
Bradlow, S.~B., Garc{\'\i}a-Prada, O., Gothen P.~B.:
Representations of the fundamental group of a surface in PU(p,q) and
holomorphic triples.
\newblock C.R. Acad. Sci. Paris, \textbf{333} 347--352 (2001).

\bibitem{bradlow-garcia-gothen:2002:preprint} Bradlow, S.~B.,
Garc{\'\i}a-Prada, O., Gothen P.~B.: Surface group representations,
Higgs bundles and holomorphic triples. \newblock Preprint, 2002.
\texttt{arXiv:math.AG/0206012}.

\bibitem{bradlow-garcia-gothen:2002}
Bradlow, S.~B., Garc{\'\i}a-Prada, O., Gothen P.~B.:
Surface  group representations and  $\U(p,q)$-Higgs bundles.
\newblock Preprint, 2002.

\bibitem{brambila-grzegorczyk-newstead:1997} 
Brambila--Paz, L., Grzegorczyk, I., Newstead, P.E.:
Geography of Brill--Noether loci for small slopes.
\newblock J. Algebraic Geometry \textbf{6} 645--669 (1997).

\bibitem{corlette:1988}
Corlette, K.: Flat ${G}$-bundles with canonical metrics.
\newblock J. Differential Geom. \textbf{28}, 361--382 (1988).

\bibitem{D1}
         Donaldson, S.~K.:
        Anti-self-dual Yang--Mills connections on a complex algebraic
              surface and stable vector bundles,
\newblock Proc. Lond. Math. Soc. {\bf 3} 1--26 (1985).

\bibitem{D2}
        Donaldson, S.~K.:
        Infinite determinants, stable bundles and curvature,
        \newblock Duke Math.
 J. {\bf 54} 231--247 (1987).

\bibitem{donaldson:1987}
Donaldson, S.~K.: Twisted harmonic maps and the self-duality equations.
\newblock Proc. London Math. Soc. (3) \textbf{55}, 127--131 (1987).


\bibitem{garcia-prada:1994}
Garc{\'\i}a-Prada, O.: Dimensional reduction of stable bundles, vortices and
  stable pairs.
\newblock Int. J. Math. \textbf{5}, 1--52 (1994).


\bibitem{G} Gieseker D.: On moduli of vector bundles on an algebraic surface,
         \newblock  Ann. of Math. {\bf 106}  45-60 (1977).

\bibitem{gothen:2001}
Gothen, P.~B.: Components of spaces of representations and stable triples.
\newblock Topology \textbf{40}, 823--850 (2001).

\bibitem{gothen-king:2002} Gothen, P.~B., King, A.~D.: Homological
algebra of quiver bundles.  
\newblock Preprint, 2002.
\newblock \texttt{arXiv: math.AG/0202033}

\bibitem{hernandez}
Hern\'andez, R.: On Harder-Narasimhan stratification over Quot schemes.
\newblock J. Reine Angew. Math. \textbf{371}, 115--124 (1986).



\bibitem{hitchin:1987}
Hitchin, N.~J.: The self-duality equations on a {R}iemann surface.
\newblock Proc. London Math. Soc.  \textbf{55}, 59--126 (1987).

\bibitem{hitchin:1992}
Hitchin, N.~J.: {L}ie groups and {T}eichm\"{u}ller space.
\newblock Topology \textbf{31}, 449--473 (1992).

\bibitem{Ko}  Kobayashi S.:
        \newblock Differential Geometry of Complex Vector Bundles,
        Princeton University Press, New Jersey, 1987.


\bibitem{markman-xia:2001}
Markman, E., Xia, E.~Z.: The moduli of flat $\mathrm{PU}(p,p)$ structures
with
  large {T}oledo invariants.
\newblock Preprint, 2001. \texttt{arXiv:math.AG/0009203 v2}.


\bibitem{maruyama} Maruyama M.:
Elementary transformations in the theory of algebraic vector bundles.
\newblock LNM   {\bf 961}, 241--266 (1982).


\bibitem{Q} Qin, Z.B.: Equivalence classes of polarizations and moduli
spaces of sheaves.
\newblock J. Differential Geom. \textbf{37} 397--415 (1993).


\bibitem{narasimhan-seshadri:1965}
Narasimhan, M.~S., Seshadri, C.~S.: Stable and unitary bundles
on a compact Riemann surface.
\newblock Ann. of Math., \textbf{82}, 540--564 (1965).



\bibitem{Schmitt} Schmitt, A.: A universal construction for the
moduli spaces of decorated vector bundles.
\newblock Habilitationsschrift, University of Essen, (2000).

\bibitem{simpson:1988}
Simpson, C.~T.: Constructing variations of {H}odge structure using
  {Y}ang-{M}ills theory and applications to uniformization.
\newblock J. Amer. Math. Soc. \textbf{1}, 867--918 (1988).

\bibitem{simpson:1992}
Simpson, C.~T.: Higgs bundles and local systems.
\newblock Inst. Hautes {\'E}tudes Sci. Publ. Math. \textbf{75}, 5--95 (1992).

\bibitem{simpson:1994a}
Simpson, C.~T.: Moduli of representations of the fundamental group of a
smooth
  projective variety {I}.
\newblock Inst. Hautes {\'E}tudes Sci. Publ. Math. \textbf{79}, 867--918
(1994).

\bibitem{simpson:1994b}
Simpson, C.~T.: Moduli of representations of the fundamental group of a
smooth
  projective variety {II}.
\newblock Inst. Hautes {\'E}tudes Sci. Publ. Math. \textbf{80}, 5--79 (1994).


\bibitem{Ta}
Taubes, C.H.: On the equivalence of the first and second order equations 
for gauge theories.
\newblock Commun. Math. Phys.\textbf{75} 207--227 (1980).

\bibitem{thaddeus:1994}
Thaddeus, M.: Stable pairs, linear systems and the {V}erlinde formula.
\newblock Invent. Math. \textbf{117}, 317--353 (1994).

\bibitem{UY}
 Uhlenbeck, K.K.  and  Yau S.T.:
         On the existence of Hermitian--Yang--Mills connections
              on stable bundles over compact K\"{a}hler manifolds.
        \newblock Comm. Pure and Appl. Math. {\bf 39--S}, 257--293
(1986).


\bibitem{W}
Witten E. Some exact multipseudoparticle solutions of classical Yang--Mills theory.
\newblock Phys. Rev. Lett. \textbf{38} (1977), 121.
\end{thebibliography}
\end{document}